\documentclass[pdflatex,sn-mathphys-num]{sn-jnl}

\usepackage{amsfonts}
\usepackage{lscape}
\usepackage[dvips]{epsfig}
\usepackage{graphicx}
\usepackage[linesnumbered]{algorithm2e}
\usepackage{amsmath}
\usepackage{url}
\usepackage{multicol}
\usepackage{subcaption}
\usepackage{diagbox}
\usepackage{multirow}
\usepackage{booktabs}
\usepackage{rotating}
\usepackage{xurl}
\usepackage{tikz}

\allowdisplaybreaks

\def\and{\hbox{\sixrm AND }}

\newcommand{\ba}{\begin{array}}
\newcommand{\ea}{\end{array}}
\newcommand{\be}{\begin{equation}}
\newcommand{\ee}{\end{equation}}

\newcommand{\re}{{\rm{I\!R}}}
\newcommand{\tr}{\rm{tr}}
\newcommand{\Id}{\rm{Id}}

\begin{document}

\title[The Ellipsoidal Separation Machine]{The Ellipsoidal Separation Machine}


\author[1]{\fnm{Antonio} \sur{Frangioni}}\email{frangio@di.unipi.it}
\equalcont{These authors contributed equally to this work.}

\author[2]{\fnm{Enrico} \sur{Gorgone}}\email{egorgone@unica.it}
\equalcont{These authors contributed equally to this work.}

\author*[2]{\fnm{Benedetto} \sur{Manca}}\email{bmanca@unica.it}
\equalcont{These authors contributed equally to this work.}

\affil*[1]{\orgdiv{Dipartimento di Informatica}, \orgname{Universita di Pisa}, \orgaddress{\street{Largo B.Pontecorvo 3}, \city{Pisa}, \postcode{56127}, \state{(PI)}, \country{Italy}}}

\affil[2]{\orgdiv{Dipartimento di Matematica e Informatica}, \orgname{Università degli Studi di Cagliari}, \orgaddress{\street{Via Ospedale 72}, \city{Cagliari}, \postcode{09124}, \state{(CA)}, \country{Italy}}}


\abstract{We propose the---to the best of our knowledge---first fully functional implementation of the ``Separation by a Convex Body'' (SCB) approach first outlined in Grzybowski et al.~\cite{grz05} for classification, separating two data sets using an ellipsoid. A training problem is defined that is structurally similar to the Support Vector Machine (SVM) one, thus leading to call our method the Ellipsoidal Separation Machine (ESM). Like SVM, the training problem is convex, and can in particular be formulated as a Semidefinite Program (SDP); however, solving it by means of standard SDP approaches does not scale to the size required by practical classification task. As an alternative, a nonconvex formulation is proposed that is amenable to a Block-Gauss-Seidel approach alternating between a much smaller SDP and a simple separable Second-Order Cone Program (SOCP). For the purpose of the classification approach the reduced SDP can even be solved approximately by relaxing it in a Lagrangian way and updating the multipliers by fast subgradient-type approaches. A characteristic of ESM is that it necessarily defines ``indeterminate points'', i.e., those that cannot be reliably classified as belonging to one of the two sets. This makes it particularly suitable for \emph{Classification with Rejection} (CwR) tasks, whereby the system explicitly indicates that classification of some points as belonging to one of the two sets is too doubtful to be reliable. We show that, in many datasets, ESM is competitive with SVM---with the kernel chosen among the three standard ones and endowed with CwR capabilities using the margin of the classifier---and in general behaves differently; thus, ESM provides another arrow in the quiver when designing CwR approaches.}

\keywords{Data classification, Convex Analysis, Nonsmooth optimization}



\maketitle

\section*{Acknowledgments}
The first author gratefully acknowledges partial funding from the European Union -- Next-GenerationEU -- National Recovery and Resilience Plan (NRRP) –- Mission 4, Component 2, Investment n.~1.1 (call PRIN 2022 D.D. 104 02-02-2022, project title ``Large-scale optimization for sustainable and resilient energy systems'', CUP I53D23002310006). The second and third authors fratefully acknowledges the support to their research by project “FIATLUCS - Favorire l’Inclusione e l’Accessibilità per i Trasferimenti Locali Urbani a Cagliari e Sobborghi” funded by the PNRR RAISE Liguria, Spoke 01, CUP: F23C24000240006.
The first and third authors are members of the group GNAMPA - Gruppo Nazionale per l'Analisi Matematica, la Probabilità e le loro Applicazioni of INDAM - Istituto Italiano di Alta Matematica.

\section{Introduction} 

Pattern classification is a supervised learning approach which recognizes similarities in the data and assigns the samples to one of data sets on the basis of a given classification rule \cite{cri00, svm2, vapnik1}. A hugely popular approach for binary classification is the Support Vector Machine (SVM) one (cf., e.g., \cite{cri00, vapnik2013, scholkopf1999}), where the classification rule is an hyperplane subdividing the space in two half-spaces, each one of which (hopefully) containing only points with the same label. Endowed with the kernel trick, that implicitly nonlinearly maps the input space into a (much larger) feature space where linear separation may be easier, SVM is a staple in Machine Learning (ML) due to its numerous advantages, chiefly among which is that the training problem is a simple convex Quadratic Program (QP). The fundamental reason is that, due to the linearity of the separation hyperplane, membership of a point to each of the two sets is an easy linear constraint.

\smallskip
\noindent
Yet, SVM has some conceptual drawbacks. One of them is related to the fact that the hyperplane has ``no locality''. That is, given the training set $X \subset \re^n$, customarily partitioned into $X^+ \cup X^-$ according to the two labels, a hyperplane is computed, possibly in the feature space, which (approximately) separates $X^+$ from $X^-$. This means that every point $x \in \re^n$ is naturally classified as being in one of the two half-spaces, even if $x$ is ``very far'' from the ``area where the original points belonged'', say an appropriate neighbourhood of $conv(X)$ (possibly after the nonlinear mapping implicit in the kernel). In this case, the new $x$ should arguably be considered an outlier, with a high risk of not belonging to the distribution represented by $X$ and therefore not classifiable. An explicit indication in this sense can be precious in many real-world application; say, in a financial application it may be preferable to remain in cash rather than investing at a loss \cite{DORAZIO, Gunjan}, in a medical application or in product inspection it may be preferable to know that more exams/checks would be needed before reaching a diagnosis/status \cite{Udler, Ozsu, Benbarrad, Sabnis}, and many others. In these Classification with Rejection (CwR) tasks \cite{Chow}, if the risk of misclassifying is too high it is preferable that the classification is not performed and the sample is rejected. This requires a \emph{reject rule}, such that the maximum of a {\it posterior} probability is less than a certain threshold. The trade-off between the error and reject probability is evaluated through a cost function that specifies both the costs of misclassification and rejection. The reject rule can be ``bolted upon'' classical binary classifiers \cite{Wegkamp, bartlett}, and is somehow naturally implemented in SVM by exploiting the margin around the classification hyperplane that is characteristic of that approach \cite{fumera2002, tortorella2004, grandvalet2008, wegkamp2010}. Yet, points $x$ that are ``outliers far away from $conv(X)$'' but still ``well off the classification hyperplane'' would still be confidently classified. Remarkably, while this is true in the feature space, the kernel trick can help address the problem in the original input space due to the corresponding nonlinear mapping, as Figures \ref{fig:polyborders}--\ref{fig:sigmoidborders} illustrate.

\smallskip
\noindent
Yet, we try to tackle the issue at a more fundamental level---directly in the input space---by relying on the concept of separation of two (convex) sets by another (convex) set in the sense of Grzybowski et al.~\cite{grz05}. Given two points $x^+$ and $x^-$, one says that a set $\mathcal{S} \subset \re^n$ separates them if the segment with extremes $x^+$ and $x^-$ touches $\mathcal{S}$ somewhere, i.e.,
\begin{equation}\label{eq:sepcond}
 conv( \, \{ \, x^+ \,,\, x^- \} \, ) \cap \mathcal{S} \neq \emptyset
 \;.
\end{equation}
The arbitrary set $\mathcal{S}$ generalizes the hyperplane which linearly separates $x^+$ and $x^-$ in SVM; indeed, any separating hyperplane satisfies \eqref{eq:sepcond}. Convexity of $\mathcal{S}$ (satisfied by an hyperplane) is convenient since it typically makes checking \eqref{eq:sepcond} an ``easy'' task. The notion is then immediately extended to two arbitrary $X^+$ and $X^-$ sets by requiring that \eqref{eq:sepcond} holds for all pairs $( \, x^+ \,,\, x^- \, ) \in X^+ \times X^-$.

\smallskip
\noindent
Intuitively, using a bounded set $\mathcal{S}$ instead of an hyperplane may help in addressing the issue mentioned above. In fact, in the usual classification task $X^+$ and $X^-$ are given as a finite set of points, i.e., $X^+ = \{ \, x^+_1 \,,\, x^+_2 \,, \ldots ,\,  x^+_{n^+} \, \}$ and $X^- = \{ \, x^-_1 \,,\, x^-_2 \,, \ldots ,\,  x^-_{n^-} \, \}$. Having $\mathcal{S}$ to separate all possible pairs $( \, x^+_i \,,\, x^-_j \, )$ is equivalent to separating the polyhedra $conv(X^+)$ and $conv(X^-)$ \cite{grz05}. A small such separating set should not be ``too much larger'' than these convex hulls, which means that it should not be able to classify a point that is ``too far away'' from the original data.

\smallskip
\noindent
Yet, putting the mechanism in action is nontrivial. Only a few attempts have been done before, mostly focusing on piecewise-linear separators. In particular, it has been proven in \cite{Gaudioso01112010} that $\mathcal{S}$ can be chosen as the Clarke subdifferential at the origin of $\min\{\,p_{conv(X^+)}\,,\,p_{conv(X^-)}\,\}$, where $p_{conv(X^+)}$ and $p_{conv(X^-)}$ denote the support functions of $conv(X^+)$ and $conv(X^-)$, respectively. The result established in \cite{Gaudioso01112010} has been effectively applied in \cite{Astorino01012011} to reduce the number of unlabeled data points in the context of semi-supervised classification. Our choice is rather to have $\mathcal{S}$ an ellipsoid, so as to exploit the powerful SemiDefinite Programming (SDP) techniques that are available to manage the positive semidefinite matrices that characterise them (and that have already been proven useful in the classification context \cite{astorino23}). Somehow similarly to an hyperplane, an ellipsoidal set in $\re^n$ requires only ``a few ingredients'' to be constructed: a positive semidefinite matrix $0 \preceq S \in \re^{n \times n}$, and a vector $\gamma \in \re^n$. Furthermore, the volume of an ellipsoid is a convex function of the parameters that characterise it, which is very convenient for finding ellipsoidal sets with minimum volume, such as the celebrated L\"owner-John ellipsoid. Defining $I^+ = \{ \, 1 \,, \ldots \,, n^+ \, \}$, $I^- = \{ \, 1 \,, \ldots \,, n^- \, \}$, and $P = \{ \, ( \, i \,,\, j \, ) \,:\, i \in I^+ \;\;,\;\; j \in I^- \, \}$ the set of all pairs of (indices of) points in the two sets to separate, we can write a first version of the separation problem, closely mirroring the SVM one (for exact linear separation) as follows:
\begin{subequations}
\begin{align}
 \min \; & size( \, S \, ) & \label{eq:ellips_def1} \\[0.05cm]
 & (x_{ij}-\gamma)^\top S(x_{ij}-\gamma) \leq 1
 & ( \, i \,,\, j \, ) \in P \label{eq:ellips_def2} \\[0.05cm]
 & x_{ij} = \alpha_{ij} x^+_i + ( 1 - \alpha_{ij} ) x^-_j
 & ( \, i \,,\, j \, ) \in P \label{eq:ellips_def3} \\[0.05cm]
 & \alpha_{ij} \in [ \, 0 \,,\, 1 \, ]
 & ( \, i \,,\, j \, ) \in P \label{eq:ellips_def4} \\[0.05cm]
 & S \succeq 0 \label{eq:ellips_def5}
\end{align}
\label{eq:ellips_naive}
\end{subequations}
In \eqref{eq:ellips_def1}, $size( \, S \, )$ means some measure of the size of the ellipsoid (which typically does not depend on $\gamma$) that is convenient from the computational viewpoint; for instance, $size( \, S \, ) = vol( \, S \, ) = \log \det S^{-1}$ can be implemented with SDP constraints and is therefore in principle attractive, but it has a nontrivial impact on the size of the formulation and therefore its solution cost. Some alternative choices will be discussed in \S \ref{sec:reform}. The formulation, as it is written, contains several cross-products between variables and therefore it is not an SDP; we will see in the next section that it is possible to reformulate it as a semidefinite, hence convex, problem. Before doing this we have to face the issue that \eqref{eq:ellips_naive} may be empty if the two sets are not ellipsoidably separable. This can be tackled, analogously to SVM, by slackening the constraints and thereby allowing their violation at a cost, but doing so in the right way requires clarifying first how the obtained ellipsoid $\mathcal{S}$ defined by $( \, S \,,\, \gamma \, )$ is going to be used to classify a new point $z \in \re^n$, a concept has not been previously clarified in the literature.

\smallskip
\noindent
The obvious definition is that $z$ is separated from $X^+$, and therefore classified as belonging to $X^-$, if:
\begin{itemize}
 \item $\mathcal{S} \cap [ \, z \,,\, x^+ \, ] \neq \emptyset$ for all $x^+ \in X^+$;
 \item $\mathcal{S} \cap [ \, z \,,\, x^- \, ] = \emptyset$ for all $x^- \in X^-$.
\end{itemize}
That is, $z$ must be entirely separated from $X^+$, and not separated by $X^-$. Of course, the symmetric relationship defines separation from $X^-$, and therefore inclusion in $X^+$. This leaves many points $z$ as not classified; in particular, for each point in $z \in \mathcal{S}$ all the segments $[ \, z \,,\, x^+ \, ]$ and $[ \, z \,,\, x^- \, ]$ have nonempty intersection with $\mathcal{S}$ and therefore fail the second part of the definition. Besides, there may be many points that are only partly, but not completely, separated from each of the two sets. Because of this we propose the following (to the best of our knowledge, new) classification rule for separation by convex bodies. Given $\mathcal{S}$, we define $\bar{X}^+ = X^+ \setminus \mathcal{S}$ and $\bar{X}^- = X^- \setminus \mathcal{S}$, i.e., we remove from the data all original points that, belonging to the ellipsoid, are not clearly separated by it; this allows to use the classification rule even in cases where \eqref{eq:ellips_naive} has no solution. Then, given an unknown point $z \in \re^n$, we compute
\begin{equation}
\begin{array}{c}
 n^+(z) = | \, \{ \, x^+ \in \bar{X}^+ \,:\,
                    \mathcal{S} \cap [ \, z \,,\, x^+ \, ]
                    \neq \emptyset \, \} \, | \\[4pt]
 n^-(z) = | \, \{ \, x^- \in \bar{X}^- \,:\,
                    \mathcal{S} \cap [ \, z \,,\, x^- \, ]
                    \neq \emptyset \, \} \, |
\end{array}
\label{eq:n+-}
\end{equation}
i.e., the cardinality of the set of points in $X^+$ (respectively, $X^-$) outside of the ellipsoid from which $z$ is properly separated. Computing this is trivial, since the segment $[ \, z \,,\, x \, ]$ intersects the ellipsoid defined by $( \, S \,,\, \gamma \, )$ if
\[
\exists \, \alpha \in [\, 0 \,,\, 1\,] \mbox{ s.t. }  (\alpha x + (1 - \alpha)z - \gamma)^\top S
 (\alpha x + (1 - \alpha)z - \gamma) \leq 1.
 \]
By expanding this expression we obtain the following second degree inequality in the variable $\alpha$
\begin{align*}
\alpha^2 (x^\top S x - 2 x^\top S z + z^\top S z) &+ 2\alpha(x^\top S z - x^\top S \gamma + z^\top S \gamma - z^\top S z)\nonumber \\ 
&+ z^\top S z - 2 z^\top S \gamma + \gamma^\top S \gamma \leq 1 
\end{align*}
which can be solved using elementary techniques.
Then, we define
\begin{equation}
\begin{array}{c}
 r^+(z) = n^-(z) \,/\, | \, \bar{X}^- \, | -
          n^+(z) \,/\, | \, \bar{X}^+ \, |
 \in [ \, -1 \,,\, 1 \, ]\\[4pt]
 r^-(z) = n^+(z) \,/\, | \, \bar{X}^+ \, | -
          n^-(z) \,/\, | \, \bar{X}^- \, |
 \in [ \, -1 \,,\, 1 \, ]
\end{array}
\label{eq:r+-}
\end{equation}
Note that, obviously, $r^+(z) + r^-(z) = 0$. If $r^+(z) = 1$ (hence, $r^-(z) = -1$), then every segment $[ \, z \,,\, x^- \, ]$ with $x^- \in \bar{X}^-$ meets $\mathcal{S}$ and no segment $[ \, z \,,\, x^+ \, ]$ with $x^+ \in \bar{X}^+$ does; hence, $z$ is classified as belonging to $X^+$ according to our definition. Conversely, if $r^-(z) = 1$ ($r^+(z) = -1$) then $z$ is classified as belonging to $X^-$. More in general, for some fixed threshold $1 \geq \varepsilon > 0$ one can classify $z$ to $X^+$ or $X^-$ according to the fact that $r^+(z) \geq \varepsilon$ or $r^-(z) \geq \varepsilon$ (the two conditions being of course incompatible). This leaves the points $z$ where both $r^+(z)$ and $r^-(z)$ are ``small'' (in absolute value) as unclassified, or rejected. A typical example is if $z \in \mathcal{S}$, since then all segments intersect the ellipsoid: $n^+(z) \,/\, | \, \bar{X}^+ \, | = n^-(z) \,/\, | \, \bar{X}^- \, | = 1$, thus $r^+(z) = r^-(z) = 0$. A different case is that of points ``very far on the side'' of $\mathcal{S}$, where $n^+(z) = n^-(z) = 0$; yet, again $r^+(z) = r^-(z) = 0$, i.e., these points cannot be classified as belonging to either set, which is reasonable. In general, points that are ``separated as much from $X^+$ as from $X^-$'' are not classified; the explicit presence of indeterminate points is a specific feature of our approach, which distinguishes it from other classification approaches based on ellipsoids \cite{astorino05, astorino23}. As already mentioned this is not a disadvantage in applications where the false positives are a significant problem: unlike other ML methods, our approach classifies a point as belonging to one of the two sets only if there is ``strong evidence'' that the point actually belongs to the set, leaving many other points of the space as not classified.
Other ML approaches also provide some measure of the confidence in the classification: in SVM, for instance, one can compare the distance from the separating hyperplane with the margin to get a sense of how ``strong'' the classification is. Yet, the separating hyperplane does naturally classify all points in the space, even if they are ``very far'' from those of the original sets. The proposed model has the potential to be more reliable in detecting if a new point can really be classified. 

\begin{figure}[h!]
\centering
    \begin{subfigure}{0.5\textwidth}
    \centering
    \begin{tikzpicture}
        \node at (0,0) {\includegraphics[scale=0.3]{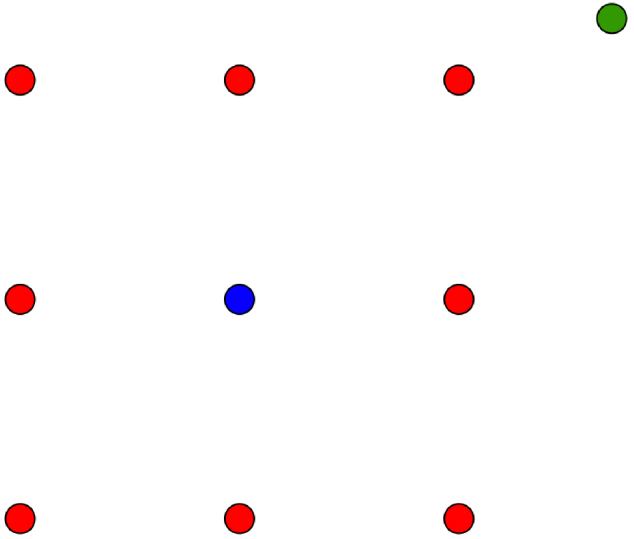}};
        \node at (1.3,1.4) {$z$};
        \end{tikzpicture}
        \caption{}
        \label{fig:classif_a}
    \end{subfigure}
    \hspace{-1cm}
    \begin{subfigure}{0.5\textwidth}
    \centering
        \includegraphics[scale=0.25]{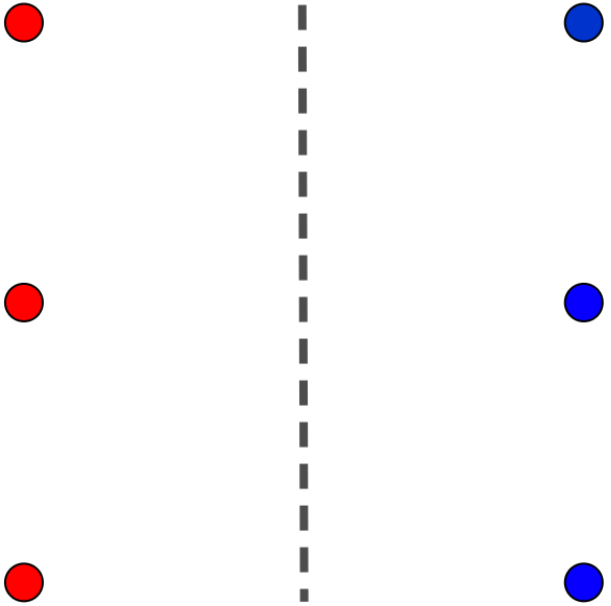}
        \caption{}
        \label{fig:classif_b}
    \end{subfigure}
    \caption{Issues with the classification rule}
    
\end{figure}

\smallskip
\noindent
It is important to note, however, that the given formul{\ae}~break if either $\bar{X}^+ = \emptyset$ or $\bar{X}^- = \emptyset$; that is, if $\mathcal{S}$ is chosen to contain one of the two sets (or both). There is nothing in \eqref{eq:ellips_naive} that helps in avoiding this issue; intuitively, the objective function could, since it leads to ``small'' ellipsoids, but this may well not be enough. Indeed, consider the 2D depicted in Figure \ref{fig:classif_a}, where blue denote(s) the (only) point(s) in $X^+$ and red those in $X^-$. Clearly, the optimal solution to \eqref{eq:ellips_naive} is just the 0-volume ellipsoid centered in blue, as it completely separates the two sets; yet, it leads to the formul{\ae}~breaking. On the other hand, it would not be reasonable to use the whole sets $X^+$ and $X^-$, as opposed to the restricted ones $\bar{X}^+$ and $\bar{X}^-$, to compute the classification indices in \eqref{eq:r+-}. The example clearly shows why: for the given point $z$, any set segment $[ \, z \,,\, x^+ \, ]$ intersects the ellipsoid while none of the segments $[ \, z \,,\, x^+_j \, ]$ do, hence considering $X^+$ for computing \eqref{eq:r+-} would result in $r^+(z) = 1$, i.e., classifying $z$ in $X^+$, which does not seem to be sensible. Besides, the training data cannot be in general assumed to be ellipsoidally separable.

\noindent
For all these reasons, we propose the following form of the training problem:
\begin{subequations}
\begin{align}
 \min \;
 & \textstyle
   size( \, S \, ) + C_1 \sum_{(i,j) \in P} \beta_{ij} +
   C_2 ( \, \sum_{i \in I^+} \nu^+_i + \sum_{j \in I^-} \nu^-_j \, )
 & + \, C_3 \|S\|_F^2 \,/\, 2
   \label{eq:SVMobj2_1}\\[0.05cm]
 & ( x_{ij} - \gamma )^\top S ( x_{ij} - \gamma ) - 1 \leq \beta_{ij}
 & ( \, i \,,\, j \, ) \in P \label{eq:SVMobj2_2} \\[0.1cm]
 & 1- ( x^+_i - \gamma)^\top S ( x^+_i - \gamma)  \leq \nu^+_i
 & i \in I^+ \label{eq:SVMobj2_3}\\[0.05cm]
 & 1- ( x^-_j - \gamma)^\top S ( x^-_j - \gamma)  \leq \nu^-_j
 & j \in I^- \label{eq:SVMobj2_4}\\[0.05cm]
 & x_{ij} = \alpha_{ij} x^+_i + ( 1 - \alpha_{ij} ) x^-_j
 & ( \, i \,,\, j \, ) \in P \label{eq:SVMobj2_5}\\[0.05cm]
 & \beta_{ij} \geq 0
 & ( \, i \,,\, j \, ) \in P \label{eq:SVMobj2_6} \\[0.05cm]
 & \nu^+_i \geq 0 & i \in I^+ \label{eq:SVMobj2_7} \\[0.05cm]
 & \nu^-_j \geq 0 & j \in I^- \label{eq:SVMobj2_8}\\[0.05cm]
 & \alpha_{ij} \in [ \, 0 \,,\, 1 \, ]
 & ( \, i \,,\, j \, ) \in P \label{eq:SVMobj2_9} \\[0.05cm]
 & S \succeq 0 \label{eq:SVMobj2_10}
\end{align}
\label{eq:ellips_nonconv}
\end{subequations}

\noindent
The ``slackened'' form \eqref{eq:SVMobj2_2} of the constraints \eqref{eq:ellips_def2} allow for segments $[ \, x^+_i \,,\, x^-_j \, ]$ not to intersect $\mathcal{S}$, this requiring $\beta_{ij} > 0$ and therefore paying the cost $C_1$ (proportional to how much the constraint is violated) in the objective \eqref{eq:SVMobj2_1}. Conversely, constraints \eqref{eq:SVMobj2_3} and \eqref{eq:SVMobj2_4} try to ensure that all the points in $X^+$ and $X^-$ remain outside of $\mathcal{S}$; those that do not, corresponding to some slack $\nu^\pm > 0$, incur a penalty $C_2$ in the objective. Besides the size term, the objective also contains a regularization term (in particular the Frobenius norm of $S$, although other forms may be chosen) that is intended to play a somewhat analogous role to the margin in the SVM model. The rationale of this last term is to avoid ``too small'' ellipsoids, which may lead to numerical issues in the solution; this point will be discussed in \S \ref{ssec:size}.

\smallskip
\noindent
Similarly to SVM, the (positive) constants $C_1$, $C_2$ and $C_3$ in the objective \eqref{eq:SVMobj2_1} of the model are the \emph{hyperparameters} of the ESM. They control the relative weight of having a ``small'' ellipsoid, w.r.t.~having a ``good separation'' of the original data (separately, of having ``many pairs $( \, i \,,\, j \, )$ separated'' vs.~having ``many original points outside of $\mathcal{S}$''), and having a ``more regular ellipsoid''. As it is customary in ML methods, grid search or related methods will be required to tune the value of the hyperparameters to each specific prediction task.

\section{Semi-definite formulations for the ESM problem}
\label{sec:reform}

The formulation \eqref{eq:ellips_nonconv} of the ESM problem contains bilinear terms, in particular of the form $v^\top S v$ where both $v$ and $S$ are variables. This makes it not obvious whether the problem is convex or not. Furthermore, we still have to specify exactly how the $size()$ term is implemented. We now discuss possible convex (re)formulations of the problem that address all these issues.

\subsection{Choices for the $size$ term}\label{ssec:size}

It is well-known that the axes of an ellipsoid correspond to its $n$ (normalised) eigenvectors, and the length of the $i$-th axis is $1 \,/\, \lambda_i$, where $\lambda_i$ is the corresponding eigenvalue. In fact, if $\lambda_i = 0$ then the axis is ``infinitely long'' and the ellipsoid is degenerate along that direction. The volume of the ellipsoid is therefore (proportional to) the product of the inverse of its eigenvalues
\[
 vol( \, S \, ) = \Pi_{i = 1}^n 1 \,/\, \lambda_i = \det( \, S^{-1} \, )
 \; ,
\]
whose logarithm is SDP-representable \cite{BenNem}; hence, $size( \, S \, ) = \log \det( \, S^{-1} \, )$ is a viable choice. However, actually representing it in a SDP program requires the introduction of $\frac{n}{2} + 2^m + 2$ additional scalar variables and $\frac{n}{2} + 2^m + 1$ semidefinite constraints, where $m$ is the smallest integer such that $2^m \geq n$, thereby significantly increasing the size of the SDP and hence its solution cost.

\smallskip
\noindent
We therefore consider the alternative measure of $S$ given by its \emph{box-size} $B(\,S\,) = \sum_{i = 1}^n 1 \,/\, \lambda_i$, i.e., the sum (rather than the product) of the eigenvalues of $S^{-1}$. The box size can be more cheapily incorporated into an SDP by means of an extra matrix variable $T \in \re^{n \times n}$ and the SDP constraint
\begin{equation}
\begin{pmatrix}
S & I \\
I & T
\end{pmatrix} \succeq 0
\quad\Longleftrightarrow\quad
\tr(\,T\,) \ge \tr(\,S^{-1}\,)
\; :
\label{eq:auxil_volconstr}
\end{equation}
by setting $size( \, S \, ) = \tr(\,T\,)$ one ensures that at optimality $\tr(\,T\,) = B(\,S\,)$, providing a somewhat cheaper measure of the ``extension'' of the ellipsoid.

\smallskip
\noindent
A relevant consequence of these choices in the objective function is that all ellipsoids considered in ESM are non-degenerate, i.e., $\lambda_i > 0$ for all $i = 1, \ldots, n$; in other words, and implicit \emph{strict positive definite constraint} $S \succ 0$ is present in the problem. This makes sense as degenerate ellipsoids are ``infinitely large'' and therefore ``much larger than $conv(\,X\,)$'', negating the basic rationale of separation by (compact) convex sets. However, these choices do nothing for the ``inverse'' form of degeneracy, which is illustrated by the (perfectly separable) case of Figure \ref{fig:classif_b}. Clearly, the optimal solution is just the vertical segment shown, that is, an ellipsoid with an axis of zero length; in other words, this corresponds to the undesirable solution $S = diag( \, +\infty \,,\, 1 \, )$ where one eigenvalue is infinitely large. This would lead the SDP not to have an optimal solution; in practice, one should expect numerical difficulties in the solution of the problem. This justifies the regularization term $\|\,S\,\|_F^2$, or any other that prevents eigenvalues from becoming too large, basically implying that the ellipsoid axes must have a minimal non-zero length. This balancing act of eigenvalues having to be neither too small nor too large, governed by the hyperparmeters, makes the approach more numerically stable and may translate in better generalization capabilities.

\smallskip
\noindent
On the other hand, even the box-size formulation using \eqref{eq:auxil_volconstr} requires an extra dense semidefinite constraint. In order to further reduce the complexity of the formulation, we notice that the box-size, by dint of minimizing the inverse of the eigenvalues of $S$, is in fact maximising the eigenvalues themselves. Thus, a possible alternative is to rather choose $size( \, S \, ) = - \tr( \, S \, )$. This avoids the semidefinite constraint, but it also removes the implicit strict definiteness constraint: it is now possible to have zero eigenvalues of $S$, which means that the ellipsoid has infinitely long axes. Yet, by making the problem significantly easier to solve it allows to perform more exhaustive grid searches for the corresponding hyperparameters, which can then hopefully be sufficient to ensure that no (or, at least, only a limited number of) eigenvalues will turn out to be actually null.

\subsection{Avoiding bilinear terms in $\gamma$ and $S$}
\label{ssec:tildeS}

We now resolve the issue of the bilinear terms in \eqref{eq:SVMobj2_2}--\eqref{eq:SVMobj2_4} due to the center $\gamma$ exploiting the fact that any non-homogeneous ellipsoid $( \, S \,,\, \gamma \, )$ in $\re^n$ with $\gamma \neq 0$ can be viewed as the intersection of an homogeneous ellipsoid $\tilde{S}$ (that is, centered at the origin) in $\re^{n+1}$ with an appropriate hyperplane. Using $\tilde{S}$ instead of $( \, S \,,\, \gamma \, )$ allows to write the condition for a point to belong (or not) to the ellipsoid using Linear Matrix Inequalities (LMIs). The important details are now briefly recalled, with more available in \cite{Glineur98}.

\smallskip
\noindent
An homogeneous ellipsoid $\tilde{S}$ in $\re^{n+1}$ is characterised by a $(n+1) \times (n+1)$ positive semidefinite matrix, that will still be denoted with $\tilde{S}$. This can be seen decomposed as
\begin{equation}
 \tilde{S}=
 \left(\begin{array}{cc} s & t^\top\\ t & F \end{array} \right)
 \;,
\end{equation}
where $s \geq 0 \in \re$, $t \in \re^n$ and $F \succeq 0 \in \re^{n \times n}$. We will ensure that $\tilde{S}$ is a nondegenerate ellipsoid, which corresponds to $s > 0$ and $F \succ 0$. Then, computing the vector $d = - F^{-1} t$, i.e., such that $t = -Fd$, for every $x \in \re^n$ it holds
\begin{equation}
\begin{array}{ll}
0 <
(\,1\,,\,x\,)^\top \tilde{S}(\,1\,,\,x\,)
& = s + 2x^\top t + x^\top Fx
  = s - 2x^\top Fd + x^\top Fx \\
& = (x-d)^\top F (x-d) + s - d^\top F d
\end{array}
\label{eq:deltadef}
\end{equation}
Since \eqref{eq:deltadef} holds for any $x$ it does in particular for $x = d$, which yields $\delta = s - d^\top F d > 0$. Moreover, if $(\,1\,,\,x\,)^\top \tilde{S}(\,1\,,\,x\,) \leq 1$ is satisfied for some $x \in \re^n$, we obtain
\begin{equation}
\label{eq:defndimellips}
(\,1\,,\,x\,)^\top\tilde{S}(\,1\,,\,x\,) \leq 1
\quad\Longleftrightarrow\quad
( x - d )^\top \frac{F}{1-\delta}( x - d ) \leq 1
\;,
\end{equation}
which implies $\delta \leq 1$. We remark that \eqref{eq:defndimellips} relies on the existence of some $x \in \re^n$ such that $(\,1\,,\,x\,) \in \tilde{S}$, which always holds if $s = 1$, i.e., scaling $\tilde{S}$ by $1 \,/\, s$. From \eqref{eq:defndimellips}, the $n$-dimensional non-homogeneous ellipsoid with $( \, S \,,\, \gamma \, ) = ( \, F \,/\, ( 1 - \delta ) \,,\, d \, )$ is characterised as the set of $x \in \re^n$ such that $(\,1\,,\,x\,)$ belongs to the homogeneous $(n+1)$-dimensional (normalised) ellipsoid $\tilde{S}$. This allows to reformulate the constraints \eqref{eq:SVMobj2_2}--\eqref{eq:SVMobj2_4} as
\begin{subequations}
\begin{align}
 &  \tilde{x}_{ij}^\top \tilde{S} \tilde{x}_{ij} - 1 \leq \beta_{ij}
 & ( \, i \,,\, j \, ) \in P \label{eq:SDPconstr_1}\\[0.05cm]
 & 1 - (\tilde{x}^+_i)^\top \tilde{S} \tilde{x}^+_i \leq \nu^+_i 
 & i \in I^+ \label{eq:SDPconstr_2}\\[0.05cm]
 & 1 - ( \tilde{x}^-_j )^\top \tilde{S} \tilde{x}^-_j  \leq \nu^-_j
 & j \in I^- \label{eq:SDPconstr_3}
 \end{align}
 \label{eq:SDPconstr}
\end{subequations}
where $\tilde{x}_{ij} = (\,1\,,\, x_{ij}\,)$, $\tilde{x}_i^+ = (\,1\,,\,x_i^+\,)$, and $\tilde{x}_j^- = (\,1\,,\,x_j^-\,)$. Owing to the fact that $\tilde{x}_i^+$ and $\tilde{x}_j^-$ are constants, \eqref{eq:SDPconstr_2} and \eqref{eq:SDPconstr_2} are already LMIs since
 $v^\top Q v = \langle \, vv^\top \,,\, Q \, \rangle$.
This is not true for \eqref{eq:SDPconstr_1} since $\tilde{x}_{ij}$ are variable, but this will be dealt with in the next subsection.

\smallskip
\noindent
We will therefore use this construction to reformulate the constraints of \eqref{eq:ellips_nonconv} using $\tilde{S}$ instead of $S$ and $\gamma$. However, we need to deal with the fact that the matrix $S$ also appears in the objective function. We propose the following approach:
\begin{itemize}
 \item we replace the term $\|S\|_F^2$ with $\|\tilde{S}\|_F^2$, since $S$ is obtained from the entries of $\tilde{S}$ and thus, the Frobenius norm of the two matrices are related;
 \item owing to the fact that $\log (1-\delta) \le 0$, one has
       \[
        \begin{array}{ll}
         \log\det S^{-1}
         & = -\log\det S = -\log\det (\frac{F}{1-\delta}) \\[0.05cm]
         & = -\log ( \frac{1}{(1-\delta)^n} \det F) =
           n \log (1-\delta) - \log \det F\\[0.05cm]
         & \le - \log \det F = \log\det F^{-1}
        \end{array}
       \]
       and therefore we use (the eigenvalues of) $F^{-1}$ in place of (these of) $S^{-1}$ to define the $size()$ term in the objective function the first two ways illustrated in \S \ref{ssec:size}.
\end{itemize}

\subsection{Avoiding bilinear terms in $x_{ij}$ and $S$}
\label{ssec:xijS}

The last required step in order to obtain a semi-definite formulation of \eqref{eq:ellips_nonconv} is to reformulate the remaining constraints \eqref{eq:SDPconstr_1} as LMI. The Schur Complement Lemma allows to rewrite them as
\[
\left(\begin{array}{cc}
\tilde{S}^{-1} & \tilde{x}_{ij}\\
\tilde{x}_{ij}^\top & 1+\beta_{ij}
\end{array}\right)\succeq 0
 \qquad
\]
However, this requires the inverse of $\tilde{S}$. This can be done in an analogous as in \eqref{eq:auxil_volconstr}, i.e.,  by the optimal solution of the problem
\[
 \displaystyle
 \min \, \Big\{ \tr(G) \,:\,
 \begin{pmatrix} \tilde{S} & \Id\\ \Id & G \end{pmatrix} \succeq 0
 \Big\}
 \; .
\]

\subsection{The final formulation}

By putting everything together, and after choosing one of the approximations for the volume term in the objective described in \S \ref{ssec:size}, we obtain the following semi-definite reformulation for the ESM problem:

\begin{subequations}
\begin{align}
 \min \;
  &
   size(\tilde{S}) + C_1 \sum_{(i,j) \in P} \beta_{ij} +
    C_2 \Big(\sum_{i \in I^+} \nu^+_i + \sum_{j \in I^-} \nu^-_j \Big)
 & \textstyle
   + \, C_3 \tr(G) +  C_4\|\tilde{S}\|_F^2\,/\, 2
   \label{eq:SDPreform_1}\\[0.05cm]
&\begin{pmatrix}
G & \tilde{x}_{ij}\\
\tilde{x}_{ij}^\top & 1+\beta_{ij}
\end{pmatrix}\succeq 0
 & ( \, i \,,\, j \, ) \in P \label{eq:SDPreform_2} \\[0.05cm]
& 1 - \big\langle \, \tilde{x}^+_i (\tilde{x}^+_i)^\top \,,\,
                 \tilde{S}  \, \big\rangle \leq \nu^+_i 
& i \in I^+ \label{eq:SDPreform_3}\\[0.05cm]
& 1 - \big\langle \, \tilde{x}^-_j (\tilde{x}^-_j)^\top \,,\,
                 \tilde{S} \, \big\rangle \leq \nu^-_j
& j \in I^- \label{eq:SDPreform_4}\\[0.05cm]
  & \begin{pmatrix} \tilde{S} & \Id\\ \Id & G \end{pmatrix} \succeq 0
  &  \label{eq:SDPreform_5}\\[0.1cm]
  & \mbox{\eqref{eq:SVMobj2_5}--\eqref{eq:SVMobj2_9}} \nonumber
\end{align}
\label{eq:SDPreform}
\end{subequations}
Although this formulation is convex, it contains a large number of semidefinite constraints; hence, it may not be able to scale to the size required by practical classification problems. Therefore, in the following section we propose a different formulation that is more amenable to (approximate) solution methods.

\section{A more compact nonconvex formulation} 

To make the formulation more practical, we first choose the cheapest of the two ways of measuring the ellipsoid size. Then, we only use the linearization tricks of the previous session, i.e., these of \S \ref{ssec:tildeS} (that give rise to linear constraints) but not these of \S \ref{ssec:xijS} (that give rise to semidefinite constraints). The resulting problem is
\begin{subequations}
\begin{align}
 \min \;
  &\textstyle
    size(\tilde{S}) + C_1 \sum_{(i,j) \in P} \beta_{ij}  +
    C_2 \big (\sum_{i \in I^+} \nu^+_i + \sum_{j \in I^-} \nu^-_j \big)
 & 
   + \, C_3 \|\tilde{S}\|_F^2 \,/\, 2
   \label{eq:finalmodel_1}\\[0.05cm]
& \tilde{S} \succeq 0 & \label{eq:finalmodel_2}\\[0.05cm]
& \tilde{x}_{ij}^\top \,\,\tilde{S}\,\, \tilde{x}_{ij} - 1 \leq \beta_{ij} & (i\,,\,j) \in P\label{eq:finalmodel_3}\\[0.05cm]
&1-(\tilde{x}_i^+)^\top\,\, \tilde{S}\,\, (\tilde{x}_i^+) \leq \nu_i^+& i \in I^+\label{eq:finalmodel_4}\\[0.05cm]
&1-(\tilde{x}_j^-)^\top\,\, \tilde{S} \,\,(\tilde{x}_j^-) \leq \nu_j^-& j \in I^-\label{eq:finalmodel_5}\\[0.05cm]
%
%
  %
  & \mbox{\eqref{eq:SVMobj2_5}--\eqref{eq:SVMobj2_9}} \nonumber
\end{align}
\label{eq:finalmodel}
\end{subequations}

\noindent
The advantage is that all the constraints---save possibly the semidefinite ones required to express $size(\tilde{S})$, depending on the choice discussed in \S \ref{ssec:size}---are linear in $\tilde{S}$ and therefore LMIs. The disadvantage is that, expanding $x$ via \eqref{eq:SVMobj2_5}, the problem contains bilinear terms in $S$ and $\alpha$ and is therefore no longer convex. However, it becomes easier once some appropriately chosen subsets of the variables are fixed, which suggests a block-Gauss-Seidel approach. In particular, when $\alpha = \alpha^*$ is fixed, and therefore $x_{ij}^* = \alpha_{ij}^* x^+_i + ( 1 - \alpha_{ij}^* ) x^-_j$ are fixed as well, the problem becomes the following SDP on the variables $\tilde{S}$, $\beta$, $\nu^+$ and $\nu^-$ 
\begin{subequations}
\begin{align}
\phi^* =  \min \;
  & size(\tilde{S}) + C_1 \sum_{(i,j) \in P} \beta_{ij} +
    C_2 \big( \sum_{i \in I^+} \nu^+_i + \sum_{j \in I^-} \nu^-_j \big)
  & + \, C_3 \|\tilde{S}\|_F^2 \,/\, 2
    \label{eq:ProgramSDP_1}
    \\[0.05cm]
& (\tilde{x}_{ij}^*)^\top \,\,\tilde{S}\,\, \tilde{x}_{ij}^* - 1 \leq \beta_{ij}
& (i\,,\,j) \in P\label{eq:ProgramSDP_3}\\[0.05cm]
& \mbox{\eqref{eq:finalmodel_2} \;\;,\;\;
        \eqref{eq:SVMobj2_6}--\eqref{eq:SVMobj2_8} \;\;,\;\;
        \eqref{eq:finalmodel_4}--\eqref{eq:finalmodel_5}}
\end{align}
\label{eq:ProgramSDP}
\end{subequations}
(where, as usual, $\tilde{x}^*_{ij} = (\,1\,,\, x_{ij}^* \,)$), and therefore can be solved by off-the-shelf methods. On the other hand, when $\tilde{S}=\tilde{S}^*$ is fixed, the problem on the $\alpha$ and $\beta$ variables becomes
\begin{subequations}
\begin{align}
\min \,
  &\textstyle
    C_1 \sum_{(i,j) \in P} \beta_{ij} \label{eq:ProgramQC_1}\\[0.05cm]
& \tilde{x}_{ij}^\top \,\,\tilde{S}^*\,\, \tilde{x}_{ij} - 1 \leq \beta_{ij} & (i\,,\,j) \in P\label{eq:ProgramQC_2}\\[0.05cm]
& \mbox{\eqref{eq:SVMobj2_5} \;\;,\;\; \eqref{eq:SVMobj2_6} \;\;,\;\;
        \eqref{eq:SVMobj2_9} }
\end{align}
\label{eq:ProgramQC}
\end{subequations}
which can be easily solved because each $\beta_{ij}$ only depends on the corresponding $\alpha_{ij}$. The block-Gauss-Seidel algorithm we propose works as follows: we fix the values of the $\alpha_{ij}$ to $1/2$ for every $i$ and $j$. We solve \eqref{eq:ProgramSDP} to obtain a starting separating ellipsoid $\tilde{S}$ and the optimal value of the objective function $\phi^*$. Then, we set the estimate $\bar{\phi} = \phi^*$ to the optimal value of \eqref{eq:ProgramSDP} and we alternatively solve \eqref{eq:ProgramQC} and \eqref{eq:ProgramSDP} a number $r$ of times. If the relative difference between $\bar{\phi}$ and the objective function value $\phi^*$ of the last SDP solved is less than a threshold $\varepsilon$ the algorithm stops, otherwise we update $\bar{\phi}$ with $\phi^*$ and we repeat the process of solving alternatively \eqref{eq:ProgramQC} and \eqref{eq:ProgramSDP} $r$ times. A scheme of the algorithm is given in Algorithm~\ref{alg:GaussSeidelAlg}. Two-set block-Gauss-Seidel approaches where each set is convex are known to converge to a critical point \cite{BlockGS}, but not necessarily to the global optimum. Indeed, it is easy to verify experimentally that the approach can get stuck in local minima of the nonconvex formulation that do not correspond to the global minima of the convex one. However, it is well-known that in learning applications such as the present one there is usually no need to determine provably optimal solutions.

\begin{algorithm}[H]
    \SetKwInOut{Input}{Input}
    \SetKwInOut{Output}{Output}
    \Input{$\varepsilon, r, \alpha=1/2 \times \mathbf{1}_{kl}$}
    $\bar{\phi} = \infty$
    $\tilde{S}, \phi^* \leftarrow$ solve \eqref{eq:ProgramSDP}
    
    \While{$|\phi^*-\bar{\phi}|> \varepsilon$}{
    $\bar{\phi} = \phi^* $

    \For{$i=1, \dots, r$}{
    $\alpha \leftarrow$ solve \eqref{eq:ProgramQC}
    
     $\tilde{S}, \phi^* \leftarrow solve \eqref{eq:ProgramSDP}$ 
    }
    }
    \Output{$\tilde{S}$}
    \caption{Ellipsoidal Separation by block-Gauss-Seidel}
    \label{alg:GaussSeidelAlg}
\end{algorithm}

\section{Lagrangian relaxation for the semi-definite program}
\label{sec:4}

Although \eqref{eq:ProgramSDP} is a smaller SDP than \eqref{eq:SDPreform}, directly using an off-the shelf solver can be prohibitively costly, especially if the number of points in $\mathcal{X}$ is large. Since the folklore is that classification methods do not need very high accuracy solutions, and the problem is solved iteratively anyway, we implemented a proximal bundle method \cite{Fr20} in order to obtain an approximate solution in a reasonable amount of time. To this aim, we first relax all the constraints except the non-negativity ones and introduce the Lagrangian multipliers $\lambda$, $\mu$, $\xi$:
\begin{subequations}
\begin{align}
\min \;
    &\textstyle
    size(\tilde{S}) + 
     \sum_{(i,j)\in P}\big( ( C_1 - \lambda_{ij})\beta_{ij} + \lambda_{ij}\tilde{x}_{ij}\tilde{S} \tilde{x}_{ij} - \lambda_{ij}\big) \nonumber \\
    &\textstyle
     +\, \sum_{i \in I^+}\big( ( C_2 - \mu_i)\nu_i^+ -\mu_i (\tilde{x}^+_i)^\top \tilde{S} \tilde{x}_i^+ + \mu_i\big) \nonumber \\
    &\textstyle
     + \, \sum_{j \in I^-}\big( ( C_2 - \xi_j)\nu_j^- - \xi_j(\tilde{x}^-_j)^\top \tilde{S} \tilde{x}_j^- + \xi_j\big)
     + C_3 \|\tilde{S}\|_F^2 \,/\, 2
     \label{eq:SDPrelaxed_1}\\[0.05cm]
 & \tilde{S} \succeq 0 \nonumber \\[0.05cm]
 & 0 \leq \beta_{ij} \leq U_{ij}
 & ( \, i \,,\, j \, ) \in P \label{eq:SDP_lagrangian_4} \\[0.05cm]
 & 0 \leq \nu^+_i \leq U^+_i
 & i \in I^+ \label{eq:SDP_lagrangian_5} \\[0.05cm]
 & 0 \leq \nu^-_j \leq U^-_j
 & j \in I^- \label{eq:SDP_lagrangian_6}
%
%
%
\end{align}
\label{eq:SDPrelaxed}
\end{subequations}
The bounds $U$, $U^+$ and $U^-$ are added in order to ensure that the relaxation is bounded, and they can be easily computed by considering some worst-case scenarios:
\begin{itemize}
\item $\beta$: for every $i$ and $j$, the quantity $\beta_{ij}$ determines how far the point $\tilde{x}_{ij}$ is inside the ellipsoid $\tilde{S}$. Therefore, as worst-case scenario for the $\beta$'s we consider the smallest sphere containing as few points as possible. To do this, we consider the point $z \in\{\tilde{x}_{ij} \mid i\le n^+ \,,\, j \le n^-\}$ closest to the origin and the sphere centered in the origin having radius equal to $r = \|z\|_2$. For every $i$ and $j$ we can then set $U_{ij} = \|\tilde{x}_{ij}\|^2_2 - r^2$.
\item $\nu^+$ and $\nu^-$: For every $i \in I^+$ and for every $j\in I^-$ the quantities $\nu^+_i$ and $\nu^-_j$ determine, respectively, how far the point $\tilde{x}_i^+$ and $\tilde{x}_j^-$ are outside of the ellipsoid $\tilde{S}$. Therefore, as a worst-case scenario we can consider the smallest sphere containing all the points considered, having radius $R = \max_{i,j}\{\|\tilde{x}_i^+\|_2, \|\tilde{x}_j^-\|_2\}$ and centered in the origin. As an upper bound for the variable $\nu^+_i$ we can then set $U^+_i = R^2 - \|\tilde{x}_i^+\|_2^2$, and a completely analogous formula works for $U^-_j$.
\end{itemize}
Problem \eqref{eq:SDPrelaxed} can be solved very efficiently since each variable $\beta_{ij}$, $\nu^+_i$ and $\nu^-_j$ is independent on all others and can be optimised separately with a closed formula. What remains is the problem in $\tilde{S}$-space:
\begin{equation}
\label{eq:sdpsep}
\begin{array}{lll}
\min \; & size(\tilde{S}) + 
          \sum_{(i,j) \in P} \tilde{x}_{ij}^\top \tilde{S} \tilde{x}_{ij}
          \\[0.05cm]
        & - \sum_{i \in I^+} (\tilde{x}^+_i)^\top \tilde{S} \tilde{x}^+_i
          - \sum_{j \in I^-} (\tilde{x}^-_j)^\top \tilde{S} \tilde{x}^-_j
   + C_3 \|\tilde{S}\|_F^2 \,/\, 2 \\[0.05cm]
s.t. & \tilde{S} \succeq 0
\end{array}
\end{equation}
The exact form and complexity of \eqref{eq:sdpsep} now crucially depends on the choice for the $size(\tilde{S})$ term. The box-size formulation, and a fortiori the one using $\log\det$, result in different (relatively) small-scale SDPs that can be efficiently solved with an off-the-shelf solver such as {\tt{Mosek}} for moderate number of features. However, as the size of $\tilde{S}$ grows, such SDPs can still be excessively costly to solve several times within the Lagrangian approaches, especially taking into account the need of repeating the fitting process for several conbinations of the hyperparameters. However, the simplest choice $size(\tilde{S}) = -\tr(\tilde{S})$ here shows a decisive advantage, since rearranging terms leads to a problem of the form
\[
 \min \big\{ \, \big\langle \, C \,,\, \tilde{S} \, \big\rangle +
                C_3 \|\tilde{S}\|_F^2 \,/\, 2 \;:\, \tilde{S} \succeq 0
             \, \big\}
\]
for some appropriate matrix $C$ depending on the values of the Lagrangian multipliers, that can then be easily recast as
\begin{equation}
    \label{eq:sdpapprox}
 \min \big\{ \, \| \, C \,/\, C_3 + \tilde{S} \, \|_F^2
             \;:\, \tilde{S} \succeq 0
             \, \big\}
 \;,
\end{equation}
i.e., as the projection of the given matrix onto the cone of semidefinite ones. This is well-known to be an easy problem that can be solved by means of the eigenvalue decomposition $C \,/\, C_3 = H \Gamma H^T$, with $\Gamma$ the diagonal matrix of eigenvalues, and then zeroing out any negative eigenvalue \cite{SDPProj}. As already noted, this may then produce degenerate ellipsoids, but it has a roughly cubic complexity on the size of $\tilde{S}$ and therefore allows the approach to scale to way larger problems than any of the other alternatives. This is why this approach has been chosen for our computational results. 

\smallskip
\noindent
The optimal solution to \eqref{eq:sdpsep} provides the function value and the subgradient, which can be used in standard ways to drive the search for the optimal Lagrangian multipliers. In practice we will not seek any accurate solution but be content with the multipliers that we obtain after a prescribed number of iterations.

\section{Numerical results}

We have implemented the ESM classifier using Python as programming language, in the context of the {\tt scikit-learn} machine learning framework \cite{scikit-learn} in order to exploit its prediction pipeline (i.e., grid-search for hyper-parameters, random train/test splitting, etc.). The small-scale quadratic master problem of the Bundle approach, has been solved with {\tt Gurobi}, while the quadratic problem \eqref{eq:ProgramQC} has been solved with {\tt Mosek}. For the three hyper-parameters of our model, $C_1$, $C_2$ and $C_3$, we performed a grid-search first inside the interval $[ \, \text{\tt 1e-12} \,,\, \text{\tt 1e+4} \, ]$ (by varying the order of magnitude alone) and then on a smaller grid according to the best values found in the first execution. We have let the Bundle parameter \textsf{max\_bundle\_it} vary in $\{ \,20 \,,\, 40 \, \}$, and number of iterations of the Gauss-Seidel approach $r$ in $\{ \, 5 \,,\, 8 \, \}$, while the other parameter $\varepsilon$ was set to {\tt 1e-4}, since preliminary experiments showed that this values corresponds to a good tradeoff between quality of the solution and execution time. As described in Section \ref{sec:4}, using the formulation \eqref{eq:sdpapprox} for the SDP problem inside the Gauss-Seidel Algorithm might yield a degenerate ellipsoid depending on the values of the hyper-parameters considered, which is particularly troubling if the topmost diagonal element ($s$ in \S \ref{ssec:tildeS}) is zero, since this does not allow to scale $\tilde{S}$. If the generated ellipsoid is degenerate we therefore add a small multiple of the identity to $\tilde{S}$. In all our experiments, combination of hyper-parameters corresponding to a degenerate ellipsoid have always been excluded by the grid search, showing that formulation \eqref{eq:sdpapprox} can indeed be used.

\smallskip
\noindent
As for the quality assessment of the classifier, we have considered the following scores:
\begin{itemize}
\item the accuracy with rejection (AR) given by the formula $(\, |X_{wc}| - \xi|X_{rej}| \, ) \,/\, |X|$, where $X_{wc}$ is the set of well-classified points, $X_{rej}$ is the set of rejected points and $\xi \in [0,0.5)$ is the rejection cost parameter;
\item the accuracy (A) on the classified points, i.e., not considering the rejected points;
\item the percentage of well classified points (WC);
\item the percentage of misclassified points (MC);
\item the percentage of rejected points (R).
\end{itemize}
Since the most interesting characteristic of our classifier is arguably the shape of its rejection zone, we have used the accuracy with rejection as score to determine the best hyper-parameters. We have compared the ESM approach with the Support Vector Machine (SVM) binary classifier with the three general-purpose kernel: Gaussian (GSVM), Polynomial (PSVM) and Sigmoid (SSVM). By default SVM doesn't have a rejection zone, therefore we have considered all the points inside its margin as rejected. In order to train SVM we have considered the following hyper-parameters: for PSVM the degree varied in $\{ \, 2 \,,\, \dots \,,\, 8 \, \}$, for GSVM and SSVM the $\gamma$ parameter varied in $[ \, \text{\tt 1e-4} \,,\, \text{\tt 1e+4} \, ]$, and for SSVM the parameter coef$_0$ varied in $[ \, \text{\tt 1e-2} \,,\, \text{\tt 1e+1} \, ]$. Moreover, for all the kernels the parameter $C$ varied in the interval $[ \, \text{\tt 1e-4} \,,\, \text{\tt 1e+4} \, ]$. For all approaches, we have performed a 10-fold cross validation splitting each dataset into training (50\% for the datasets with less than 1000 samples and 30\% for the others) and testing (50\% for the datasets with less than 1000 samples and 70\% for the others) in order to obtain the best hyper-parameters. All the numerical experiments have been executed on a Mac Studio with a 32-core M3 Ultra CPU and 512 GB of RAM.

\begin{figure}[h!]
\begin{center}
\includegraphics[scale=0.28]{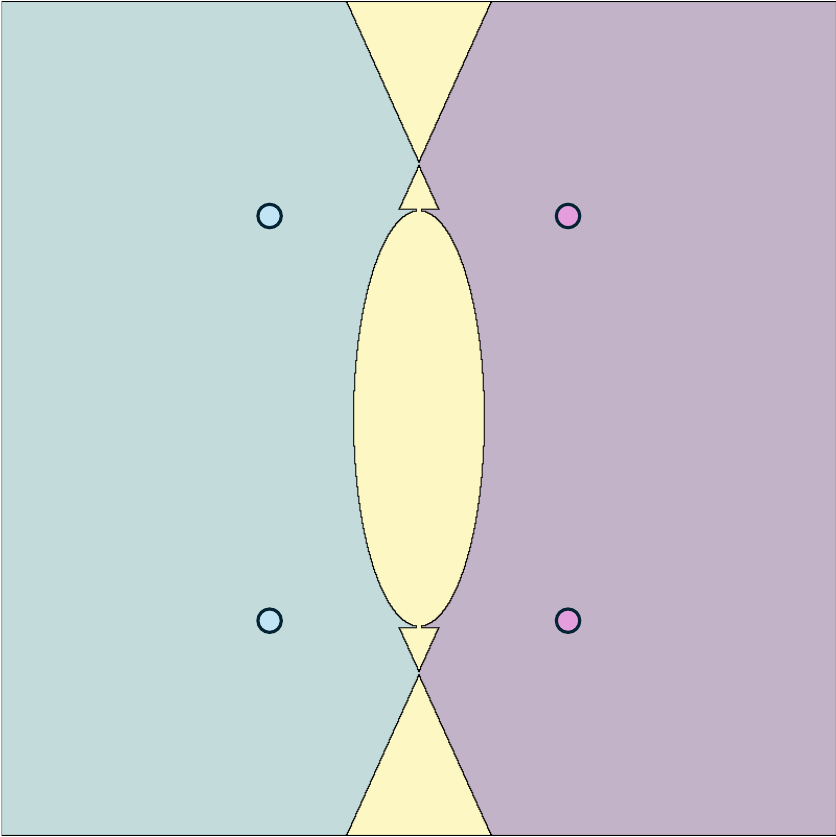}\;\;%
\includegraphics[scale=0.28]{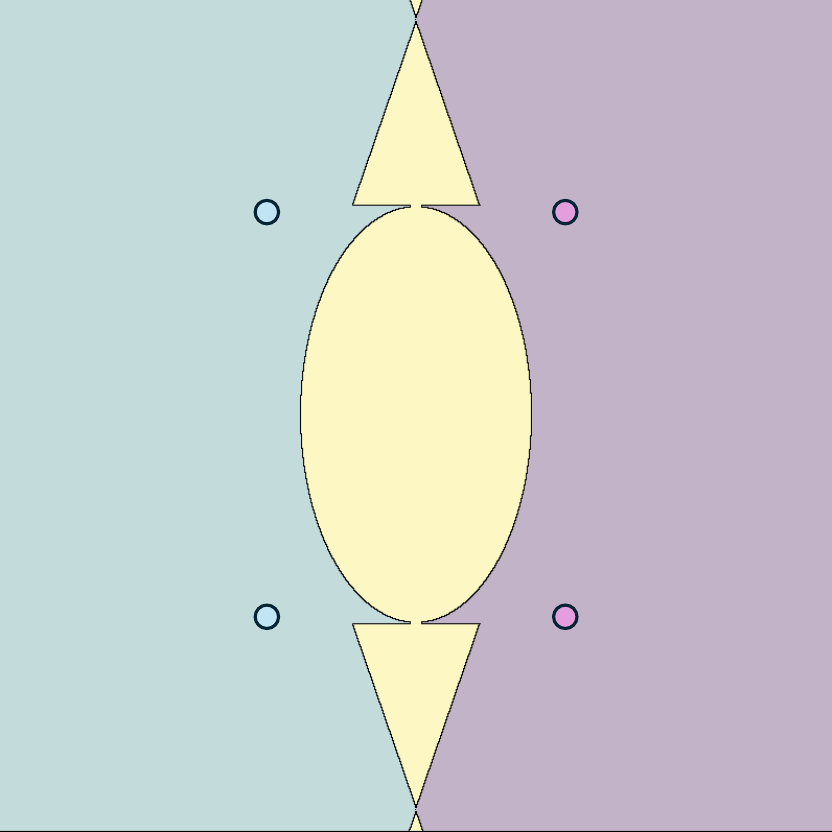}\;\;%
\includegraphics[scale=0.28]{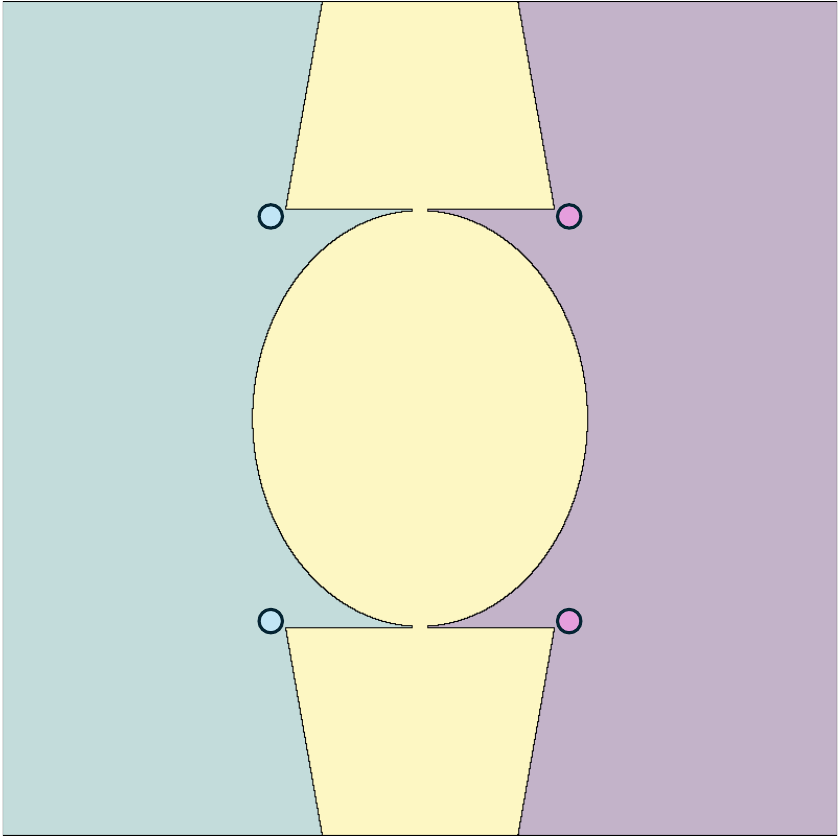}%
\end{center}
\caption{Decision borders for the ESM classifier with $C_3$ equals to 0.01, 0.1 and 0.5 (from left to right)}
\label{fig:ESMborders}
\end{figure}

\begin{figure}[h!]
\begin{center}
\includegraphics[scale=0.28]{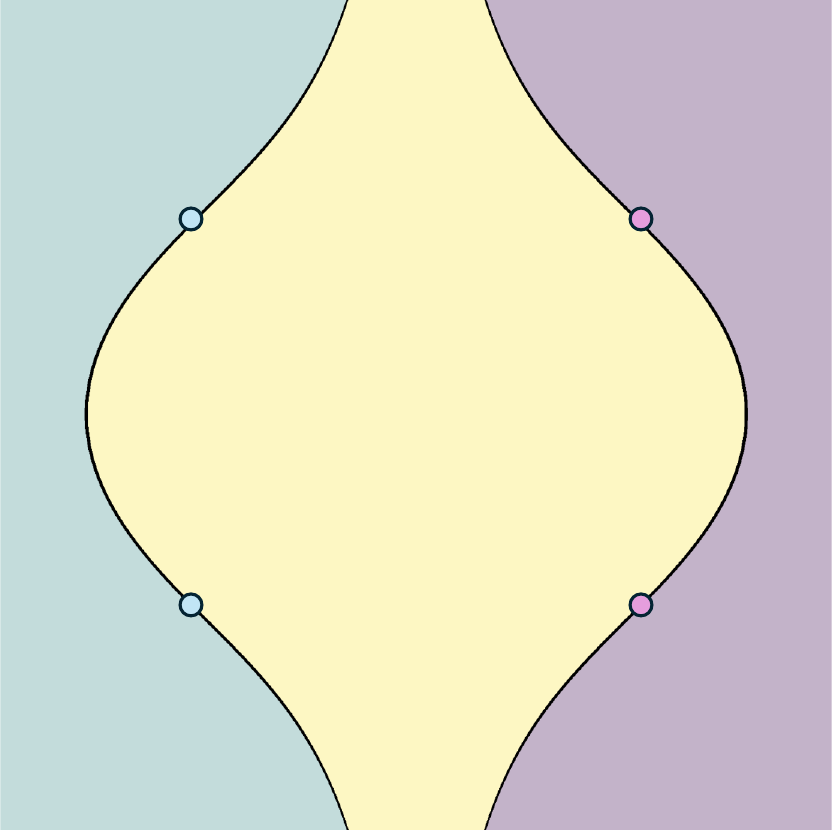}\;\;%
\includegraphics[scale=0.28]{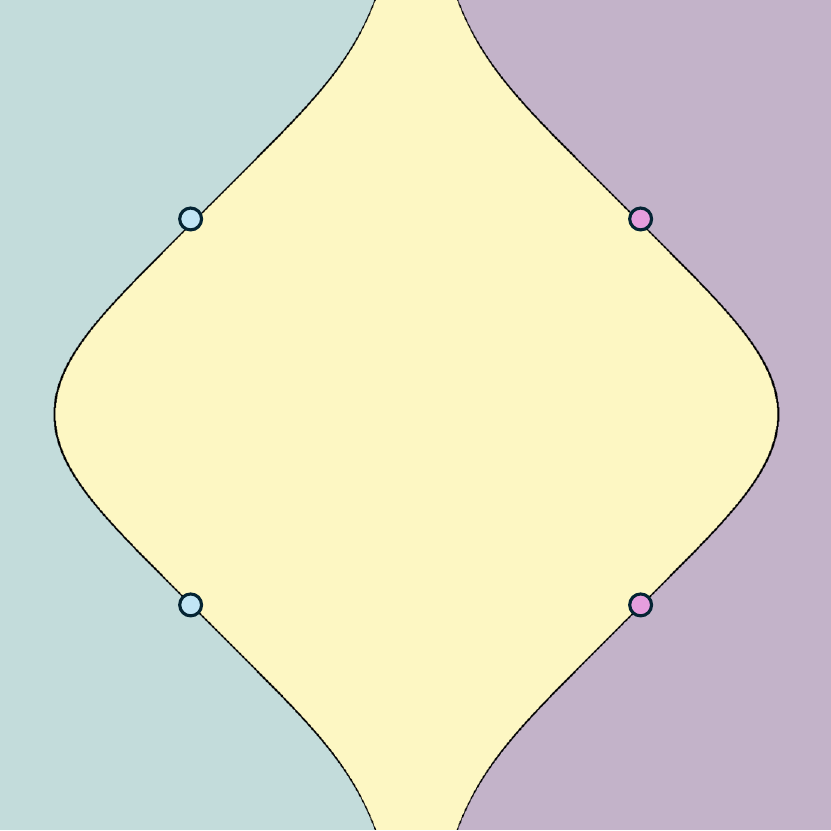}\;\;%
\includegraphics[scale=0.28]{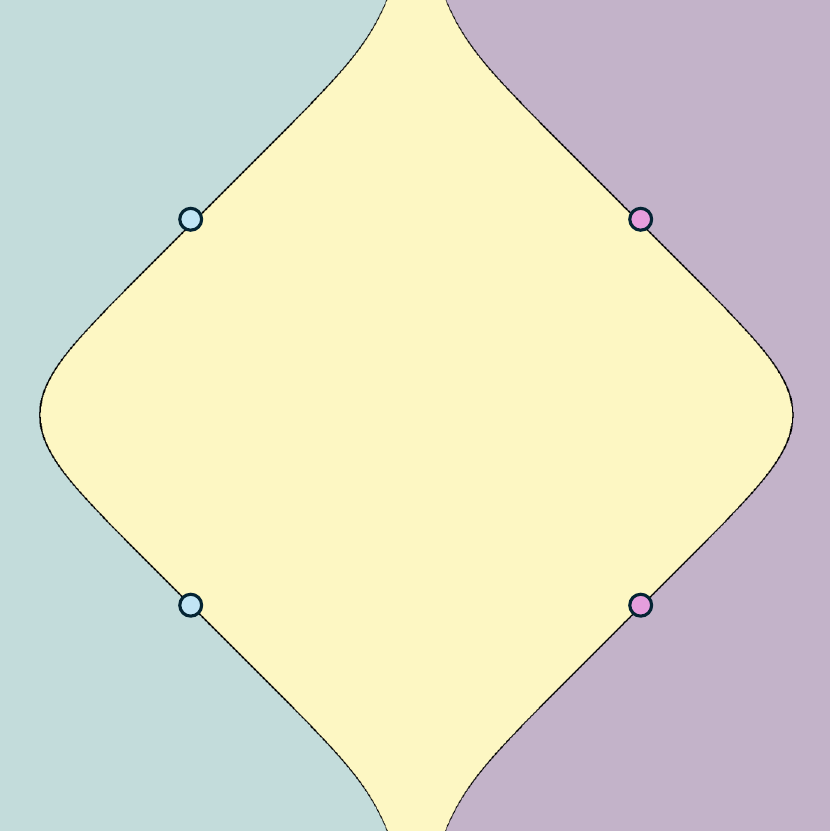}%
\end{center}
\caption{Decision borders for SVM with polynomial kernel with its hyper-parameter degree equals to 3, 5 and 7 (from left to right)}
\label{fig:polyborders}
\end{figure}

\begin{figure}[h!]
\begin{center}
\includegraphics[scale=0.28]{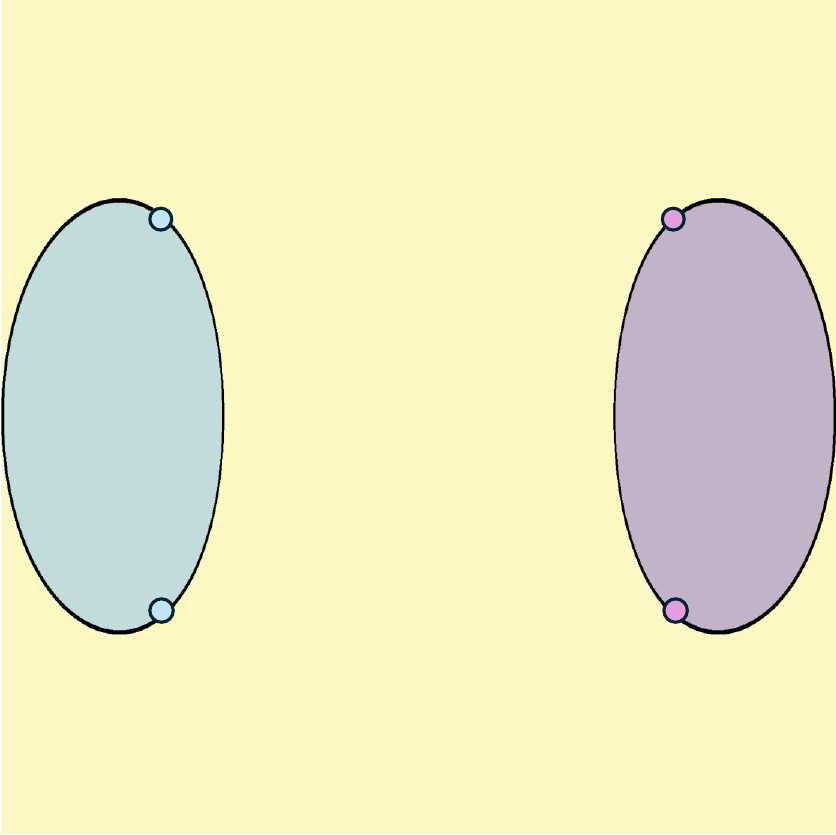}\;\;%
\includegraphics[scale=0.28]{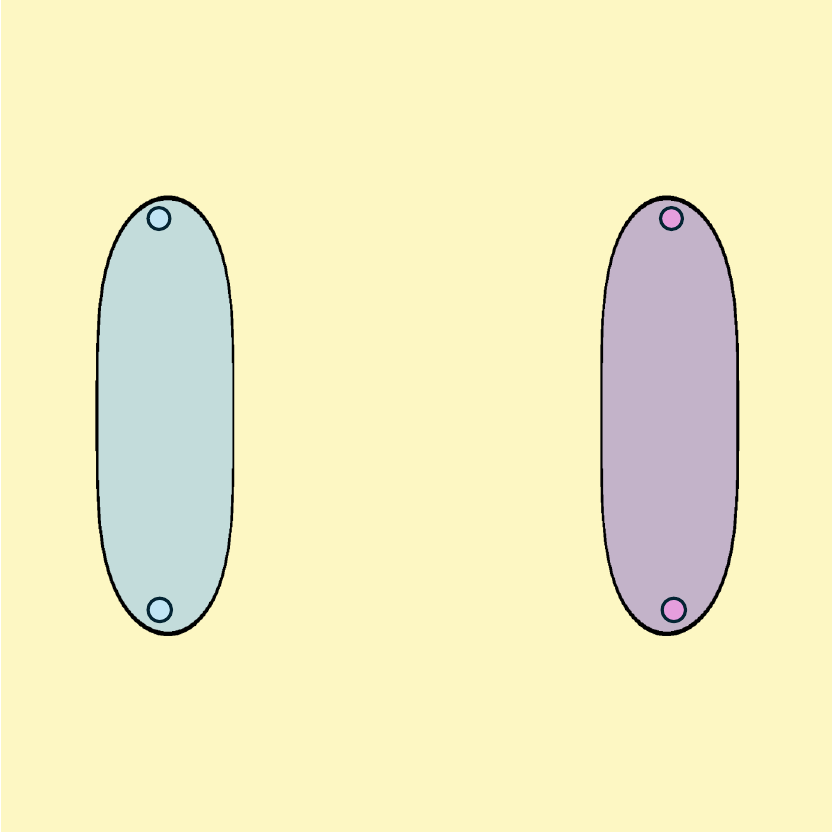}\;\;%
\includegraphics[scale=0.28]{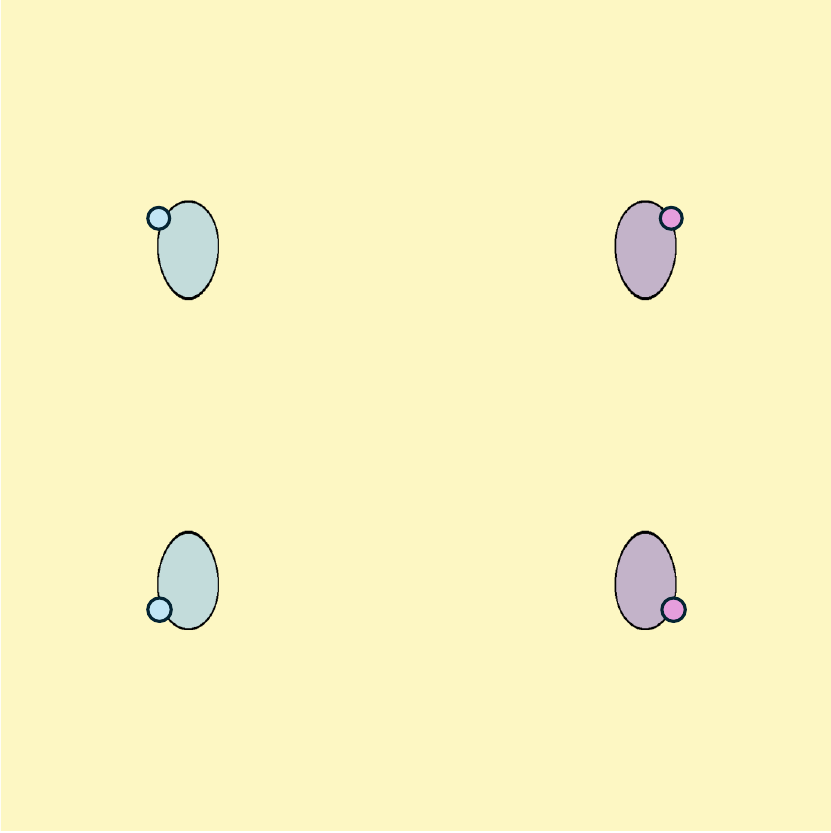}
\end{center}
\caption{Decision borders for SVM with Gaussian kernel with its hyper-parameter $\gamma$ equals to 0.3, 0.5 and 0.7 (from left to right)}
\label{fig:gaussborders}
\end{figure}

\begin{figure}[h!]
\begin{center}
\includegraphics[scale=0.28]{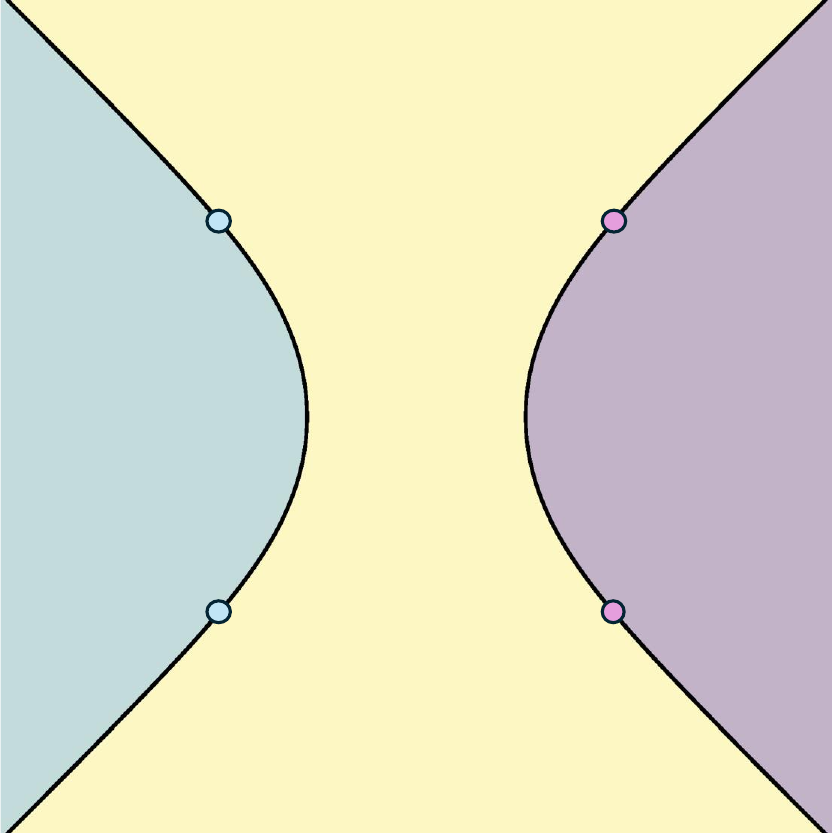}\;\;%
\includegraphics[scale=0.28]{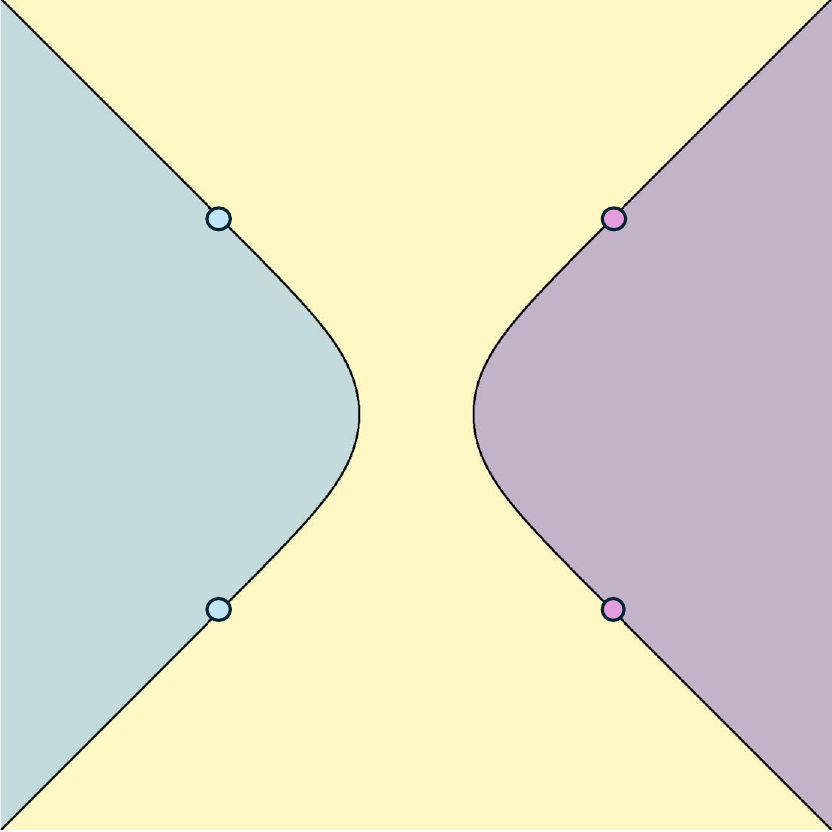}\;\;%
\includegraphics[scale=0.28]{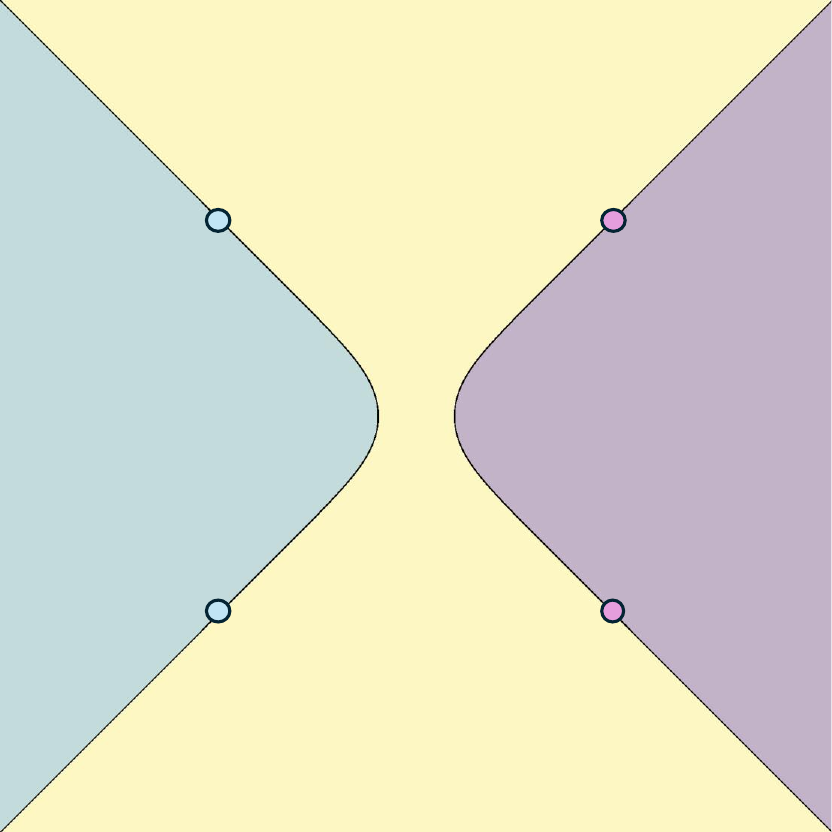}
\end{center}
\caption{Decision borders for SVM with Sigmoidal kernel with its hyper-parameter $\gamma$ equals to 1, 2 and 3 (from left to right)}
\label{fig:sigmoidborders}
\end{figure}

\smallskip
\noindent
Interestingly, the shape of the rejection zones can be expected to vary significantly among the different approaches. This is shown in Figures \ref{fig:ESMborders}--\ref{fig:sigmoidborders} for a trivial instance by plotting the decision boundaries for ESM and SVM with the kernels we have considered. The pictures confirm that ESM has a quite different decision strategy with respect to SVM, regardless of the chosen kernel. 

\smallskip
\noindent
The dataset we have considered have been taken from the libSVM repository \cite{LIBSVM} and from the UCI Machine Learning repository \cite{UCI}. The names and the dimensions of the datasets are reported in Table \ref{tab:realdat}.

\begin{table}[h!]
\begin{tabular}{| c | c | c |}
\hline
Name & n. points & n. features\\
\hline
    a1a &1605 &119 \\
    \hline
 Australian&690 &14 \\
 \hline
 Breast &683 &10 \\
 \hline
 Breast Wisconsin &569 & 29 \\
 \hline
  Diabetes& 768&8 \\
  \hline
  Divorce&170 & 54\\
  \hline
  Fertility &100 &9\\
  \hline
  Flowmeters&87 &36 \\
  \hline
  Gallstone& 319 & 38\\
 \hline
 Heart &270 &13 \\
 \hline
 Ionosphere &351 &34 \\
 \hline
  Liver & 345&5 \\
  \hline
  Mesothelioma&324 & 34\\
  \hline
  Sonar&208 & 60 \\
 \hline
 svmguide3 &1243 &22 \\
 \hline
 w1a &2477 &300 \\
 \hline
\end{tabular}
\caption{Dataset information}
\label{tab:realdat}
\end{table}

\begin{figure}[h!]
\begin{center}
\begin{subfigure}{0.3\textwidth}
\includegraphics[scale=0.25]{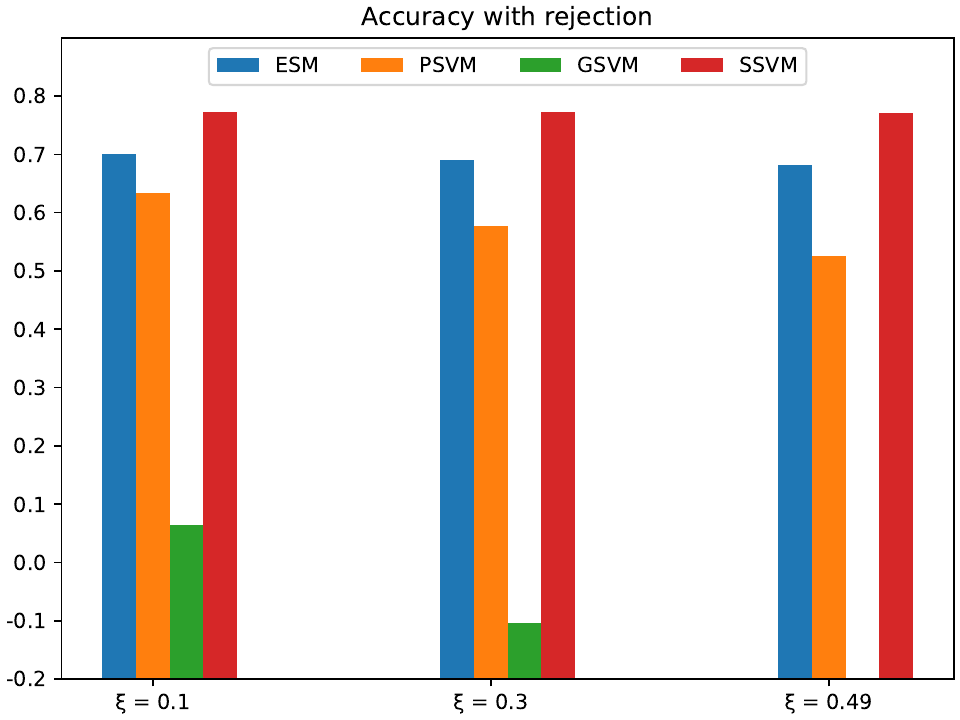}
\end{subfigure}
\hspace{1cm}
\begin{subfigure}{0.3\textwidth}
\includegraphics[scale=0.25]{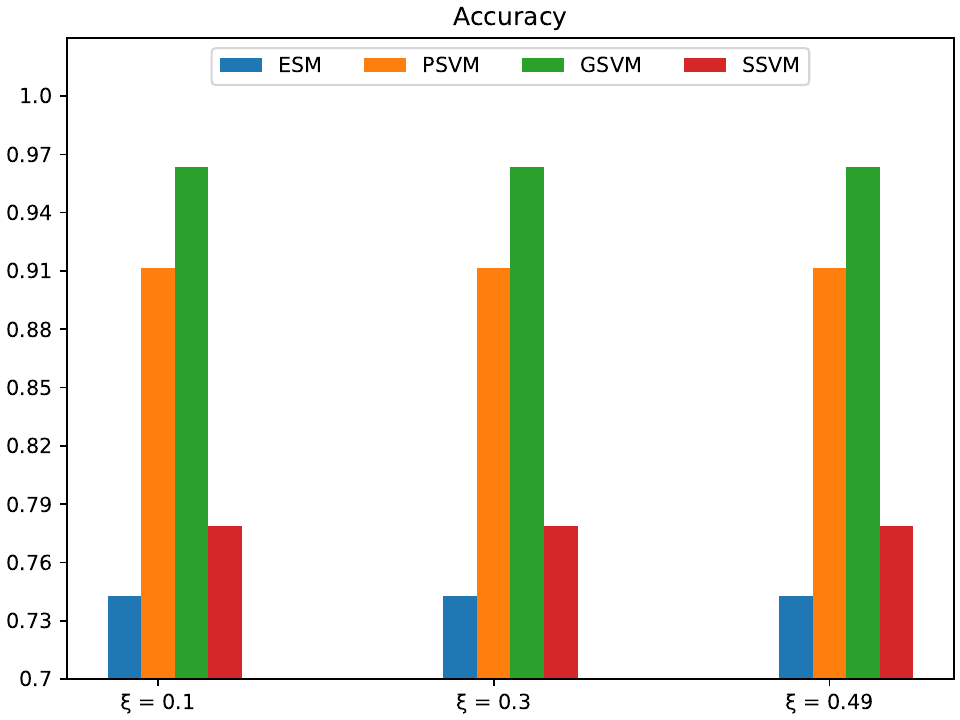}
\end{subfigure}
\vspace{0.3cm}
\begin{subfigure}{0.3\textwidth}
\includegraphics[scale=0.25]{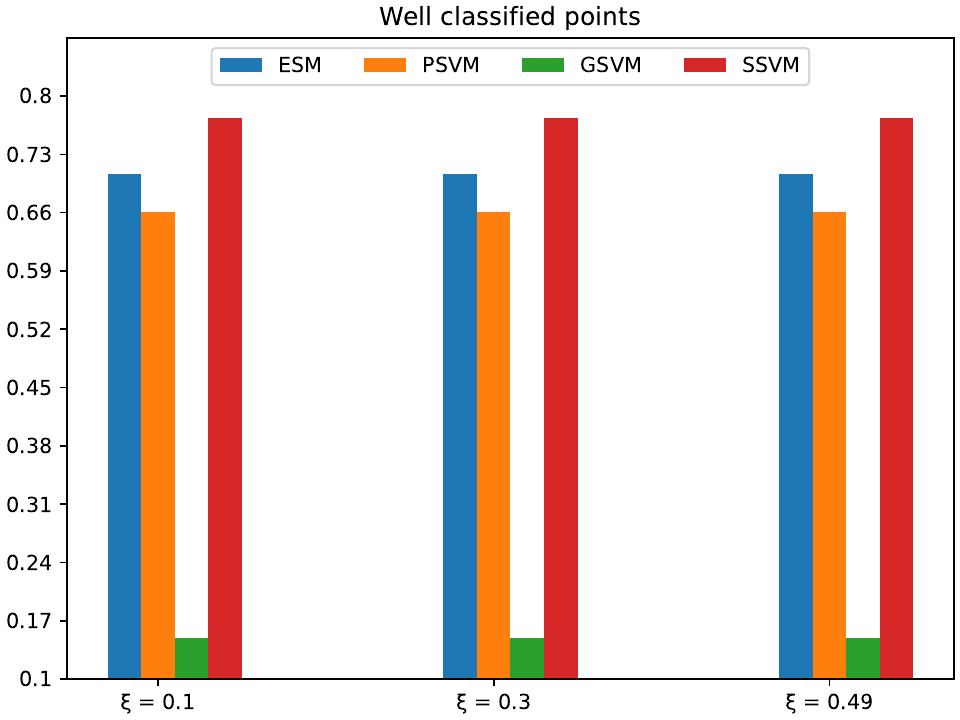}
\end{subfigure}
\hspace{0.1cm}
\begin{subfigure}{0.3\textwidth}
\includegraphics[scale=0.25]{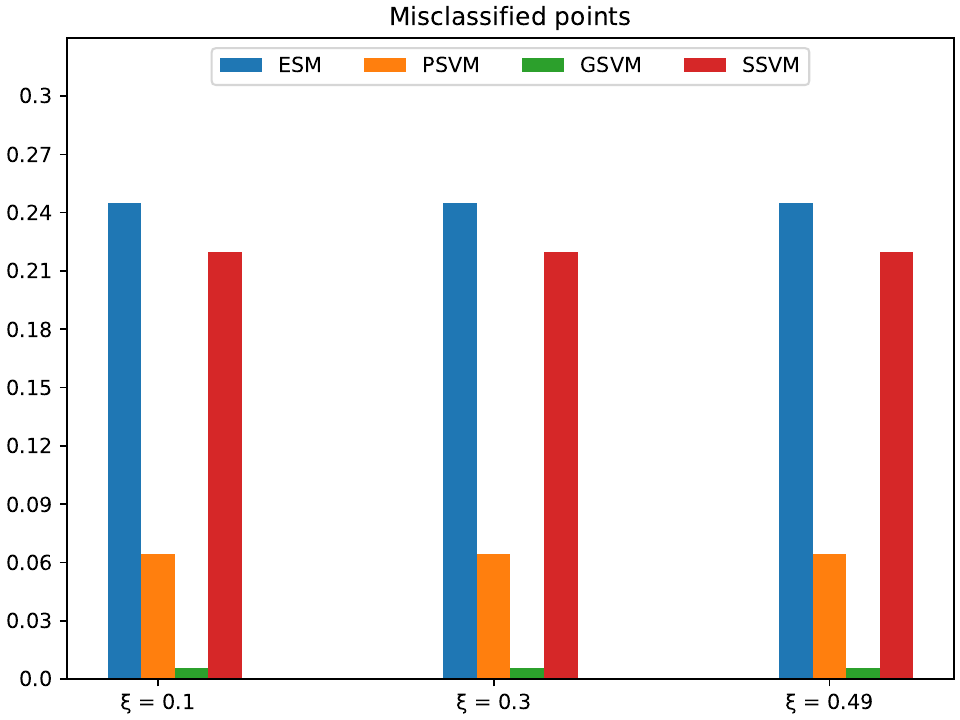}
\end{subfigure}
\hspace{0.1cm}
\begin{subfigure}{0.3\textwidth}
\includegraphics[scale=0.25]{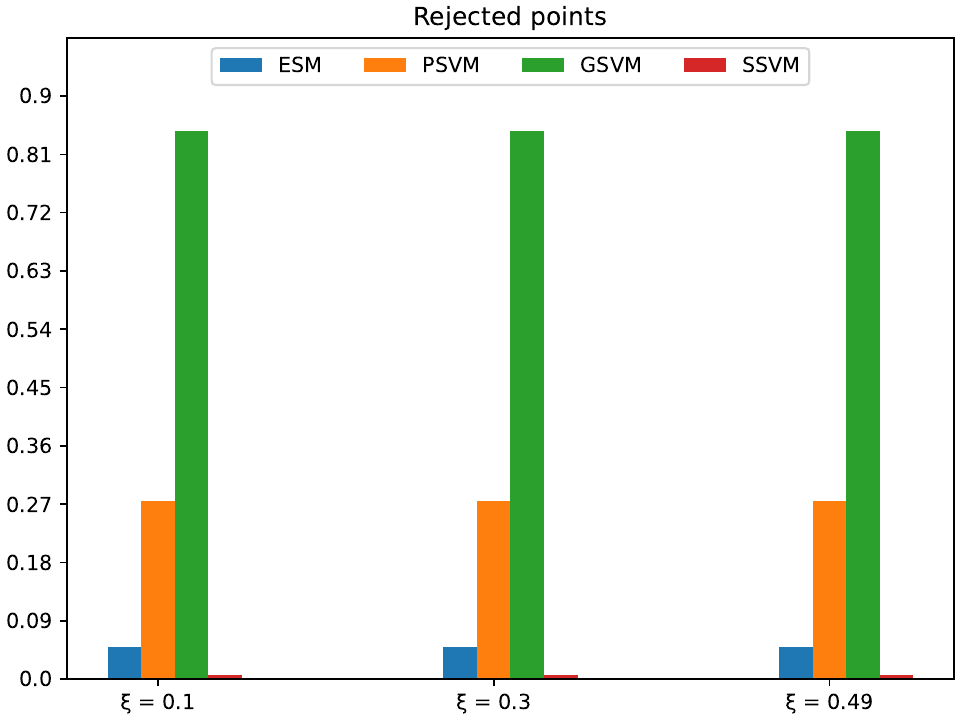}
\end{subfigure}
\end{center}
\caption{Numerical results for the real dataset a1a}
\label{fig:a1a}
\end{figure}

\begin{figure}[h!]
\begin{center}
\begin{subfigure}{0.3\textwidth}
\includegraphics[scale=0.25]{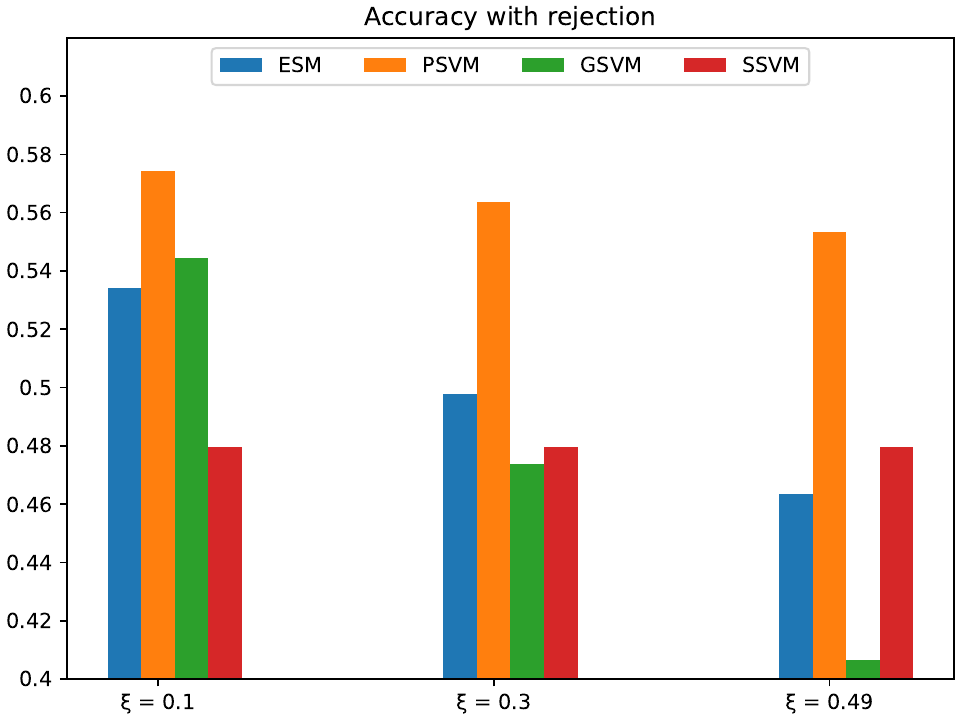}
\end{subfigure}
\hspace{1cm}
\begin{subfigure}{0.3\textwidth}
\includegraphics[scale=0.25]{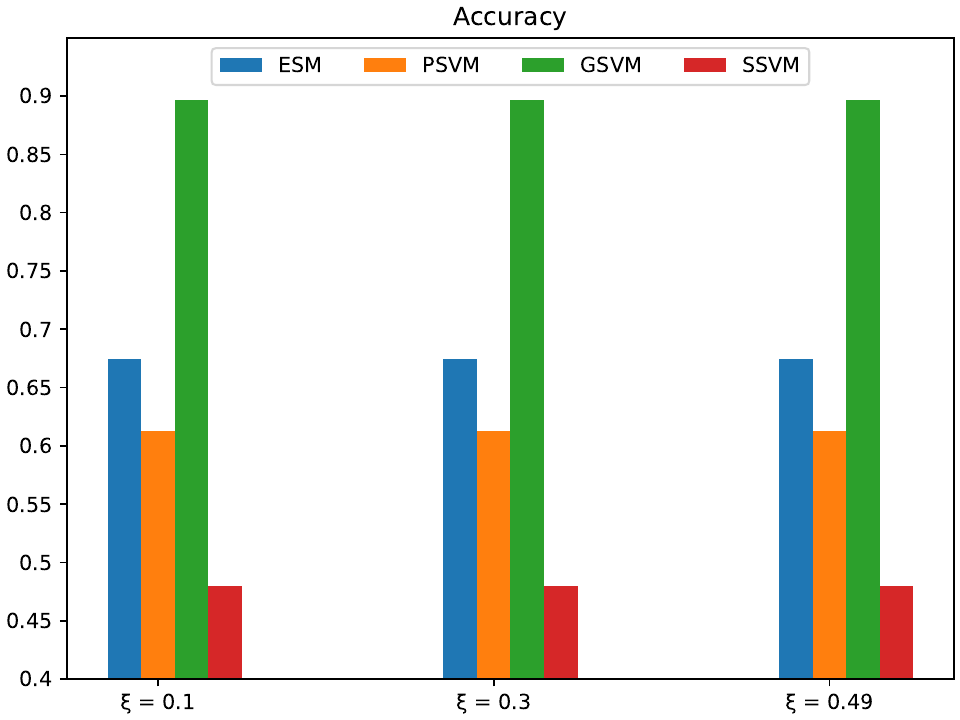}
\end{subfigure}
\vspace{0.3cm}
\begin{subfigure}{0.3\textwidth}
\includegraphics[scale=0.25]{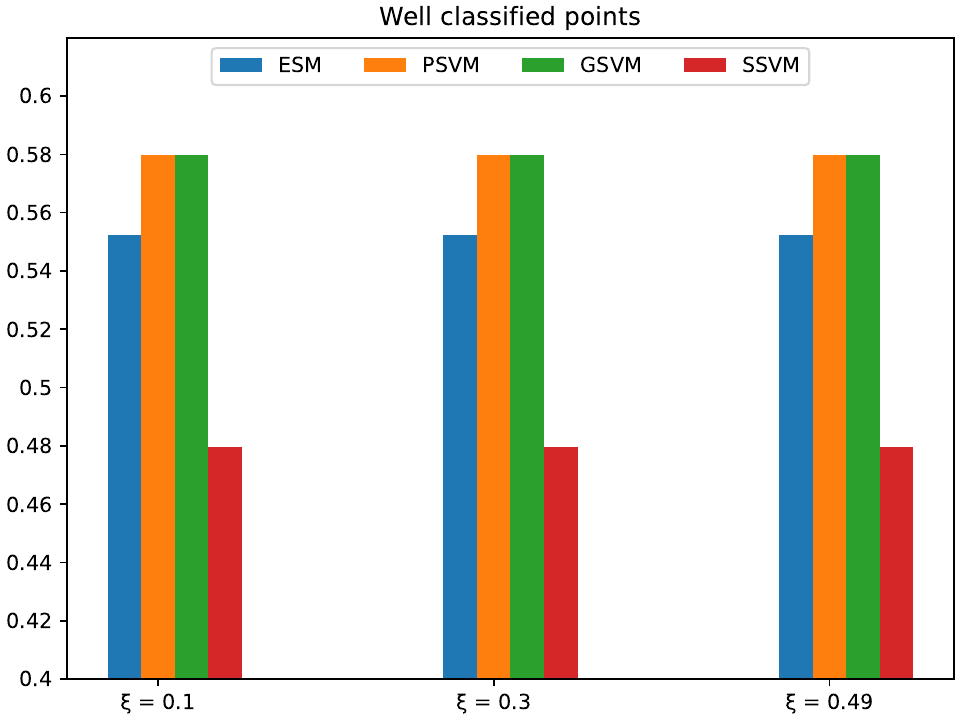}
\end{subfigure}
\hspace{0.1cm}
\begin{subfigure}{0.3\textwidth}
\includegraphics[scale=0.25]{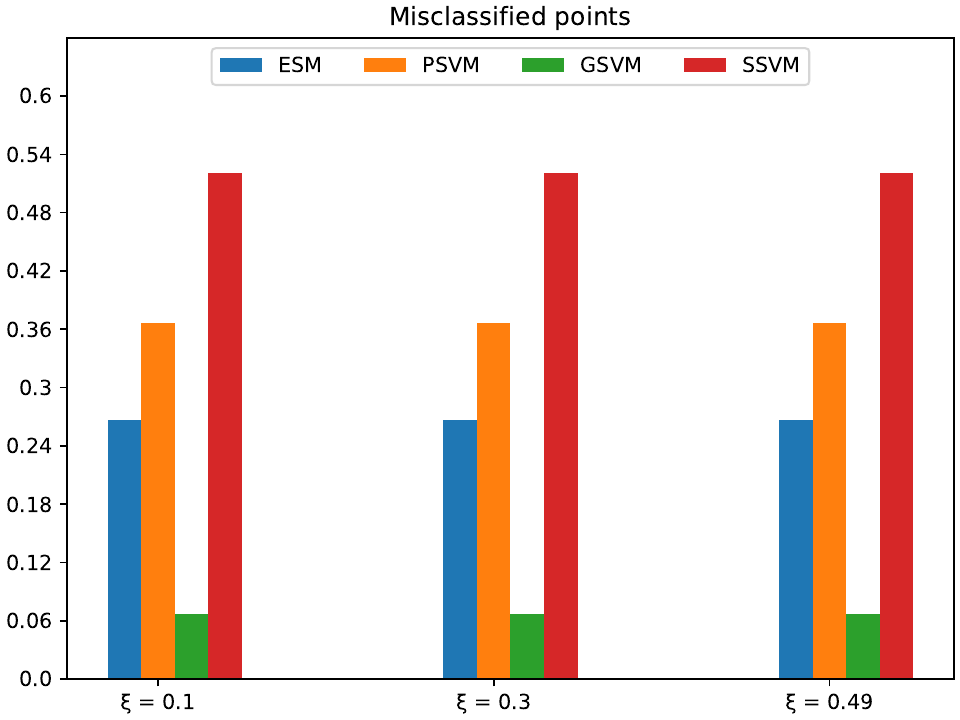}
\end{subfigure}
\hspace{0.1cm}
\begin{subfigure}{0.3\textwidth}
\includegraphics[scale=0.25]{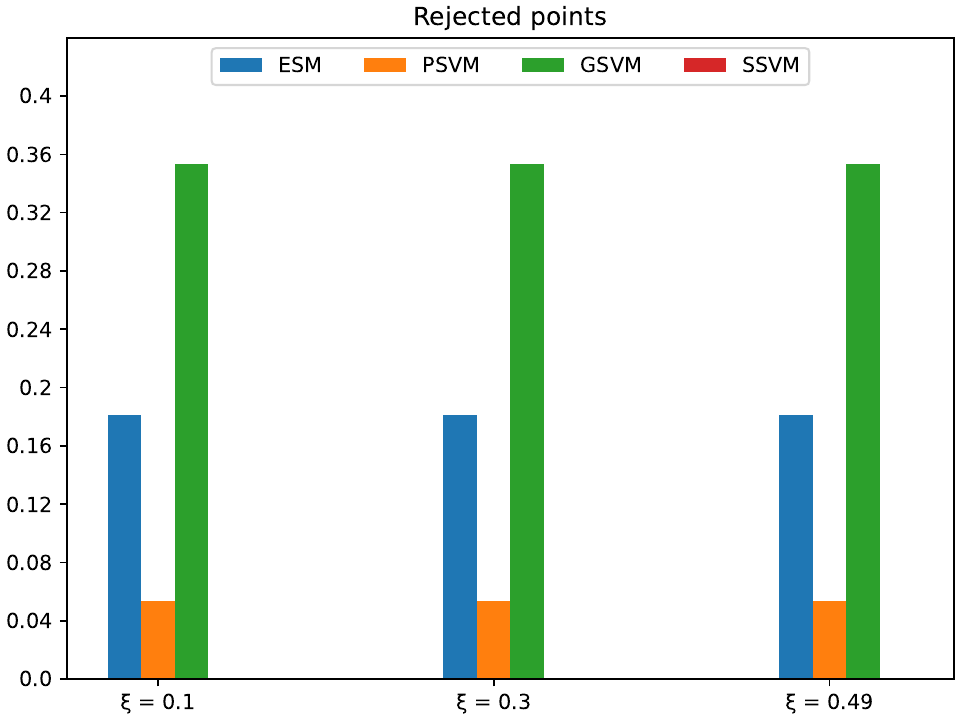}
\end{subfigure}
\end{center}
\caption{Numerical results for the real dataset Australian}
\label{fig:austr}
\end{figure}

\begin{figure}[h!]
\begin{center}
\begin{subfigure}{0.3\textwidth}
\includegraphics[scale=0.25]{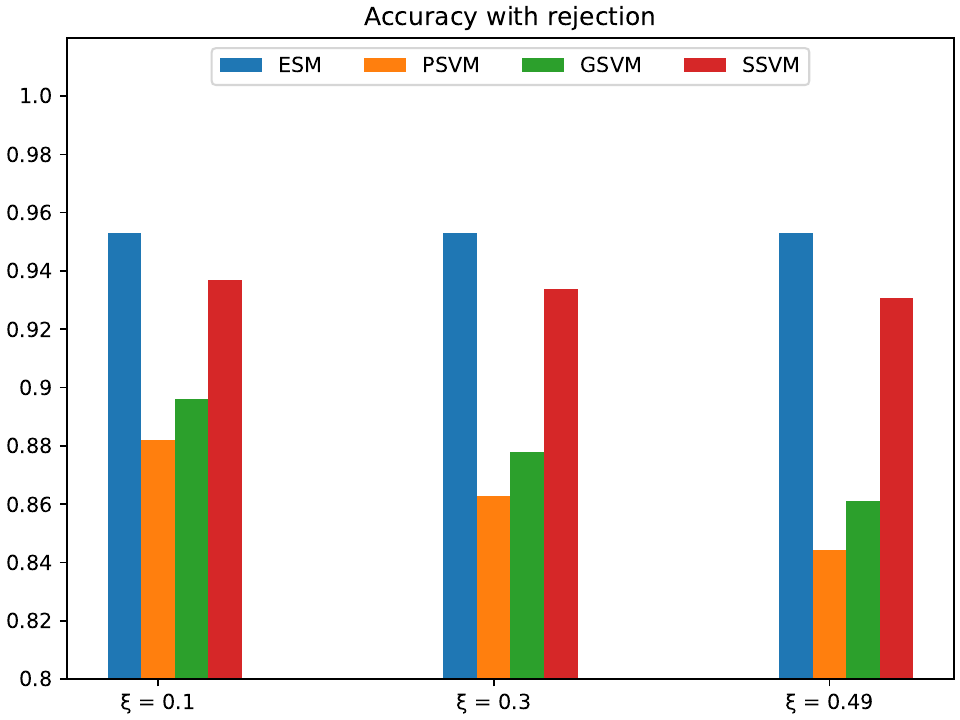}
\end{subfigure}
\hspace{1cm}
\begin{subfigure}{0.3\textwidth}
\includegraphics[scale=0.25]{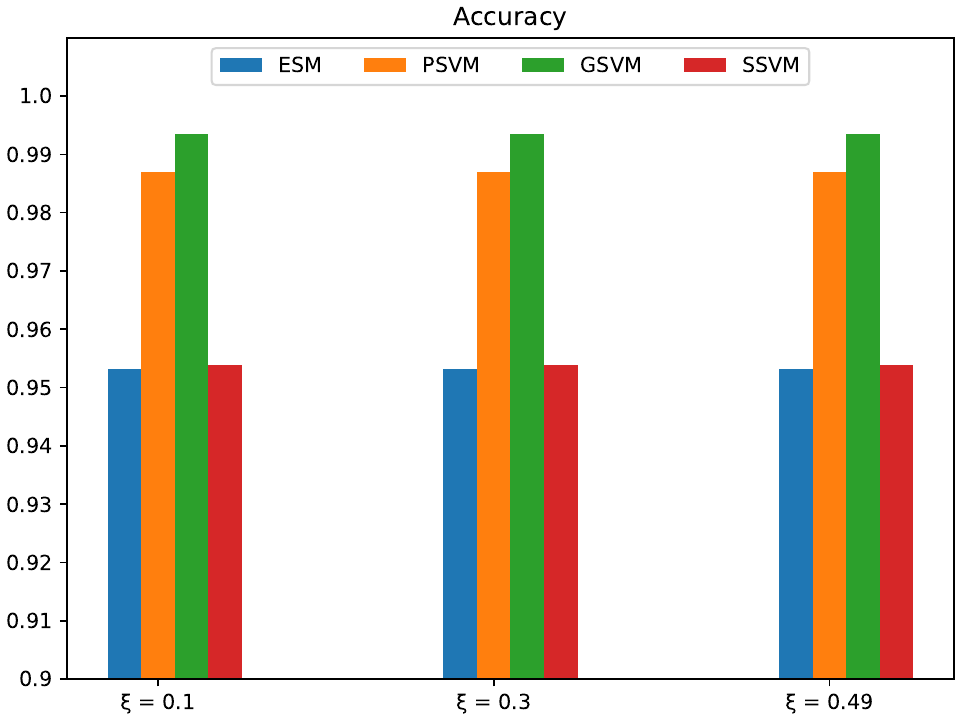}
\end{subfigure}
\vspace{0.3cm}
\begin{subfigure}{0.3\textwidth}
\includegraphics[scale=0.25]{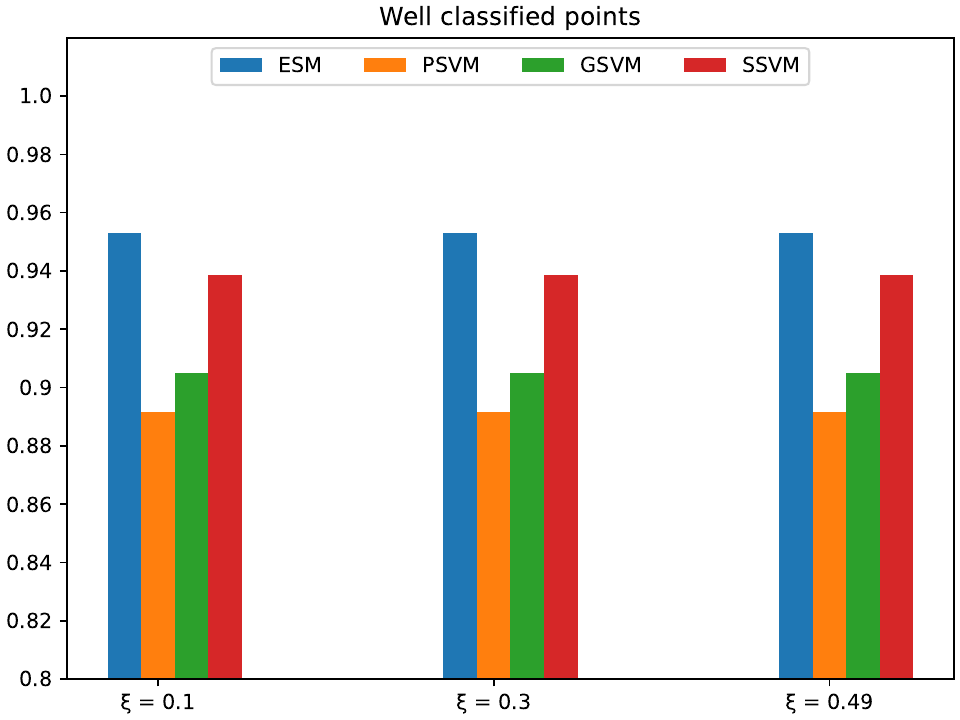}
\end{subfigure}
\hspace{0.1cm}
\begin{subfigure}{0.3\textwidth}
\includegraphics[scale=0.25]{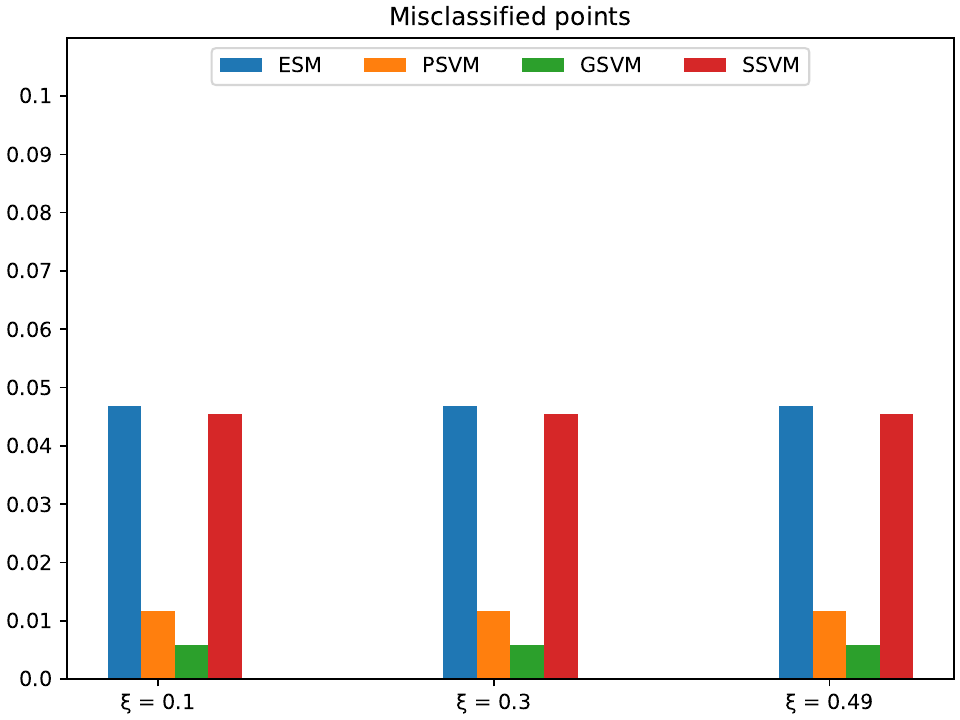}
\end{subfigure}
\hspace{0.1cm}
\begin{subfigure}{0.3\textwidth}
\includegraphics[scale=0.25]{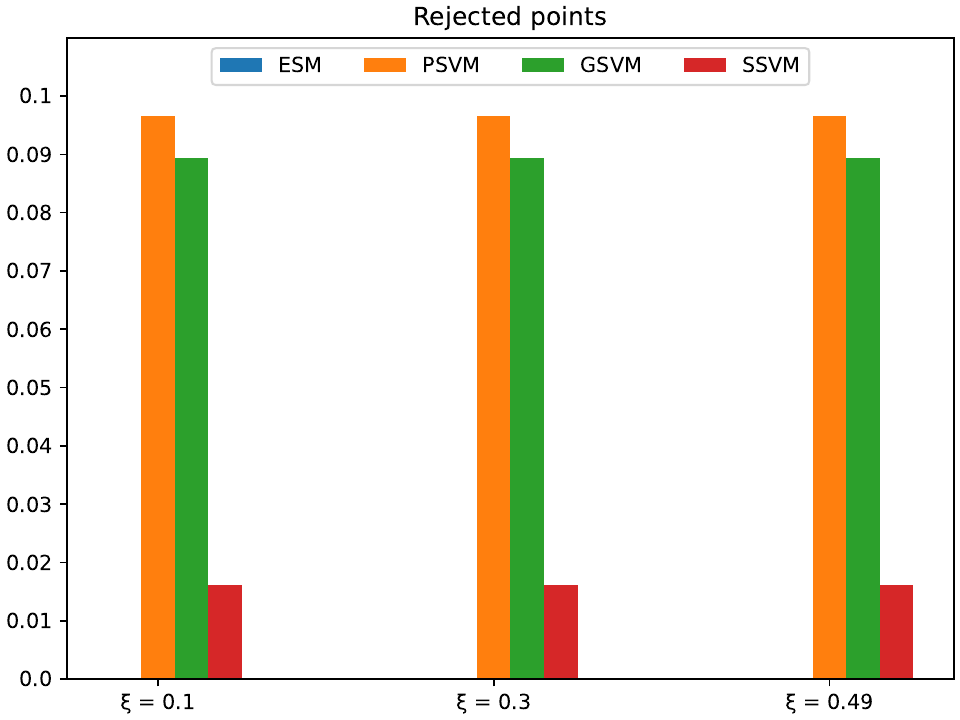}
\end{subfigure}
\end{center}
\caption{Numerical results for the real dataset Breast}
\label{fig:breast}
\end{figure}

\begin{figure}[h!]
\begin{center}
\begin{subfigure}{0.3\textwidth}
\includegraphics[scale=0.25]{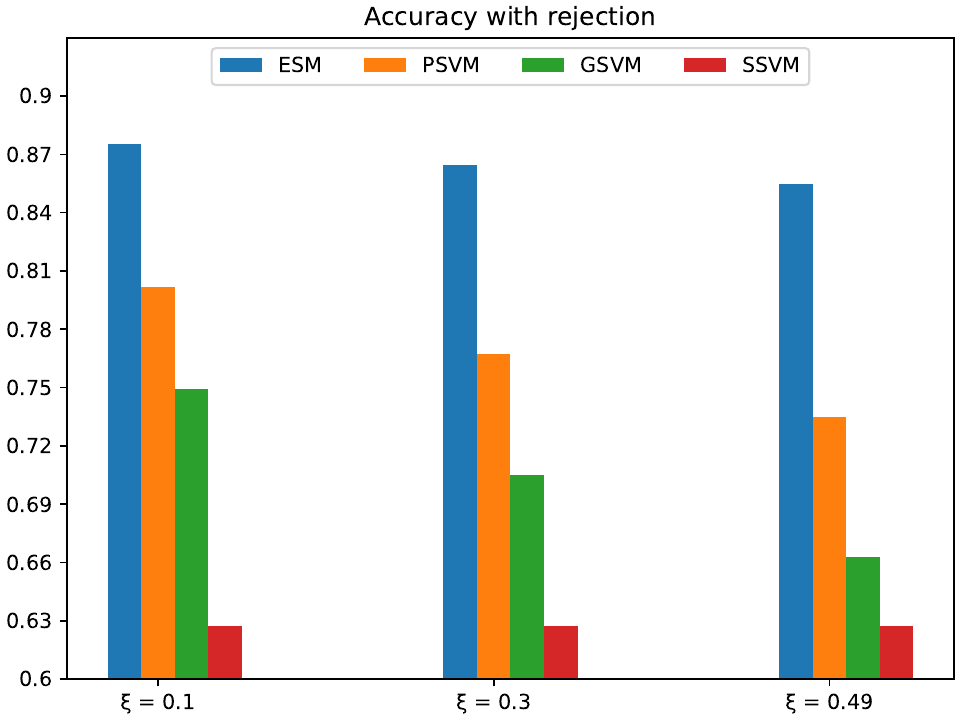}
\end{subfigure}
\hspace{1cm}
\begin{subfigure}{0.3\textwidth}
\includegraphics[scale=0.25]{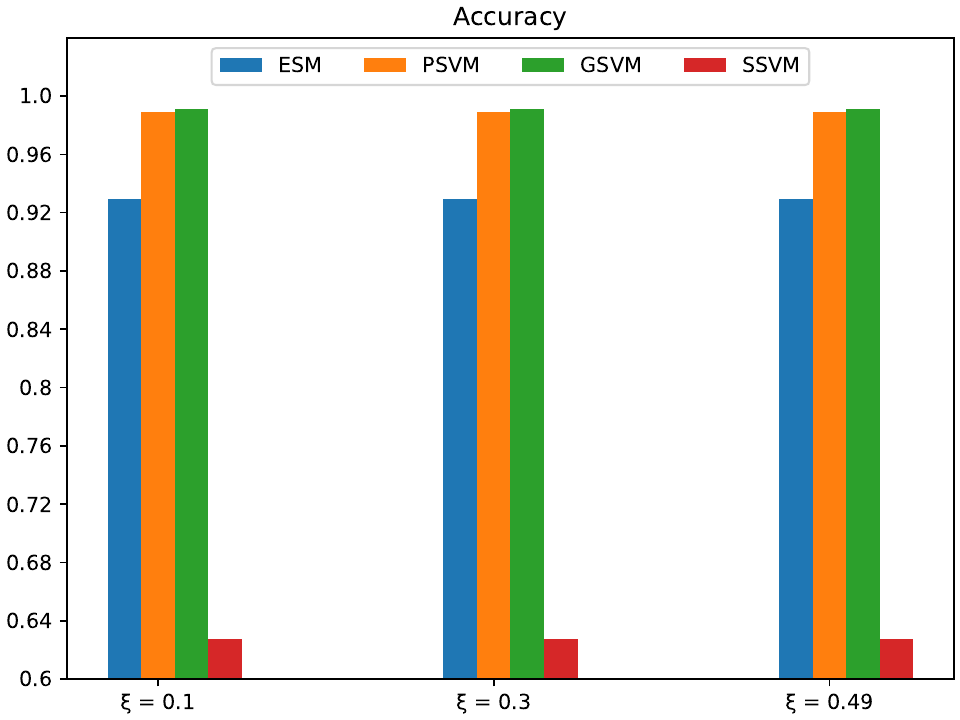}
\end{subfigure}
\vspace{0.3cm}
\begin{subfigure}{0.3\textwidth}
\includegraphics[scale=0.25]{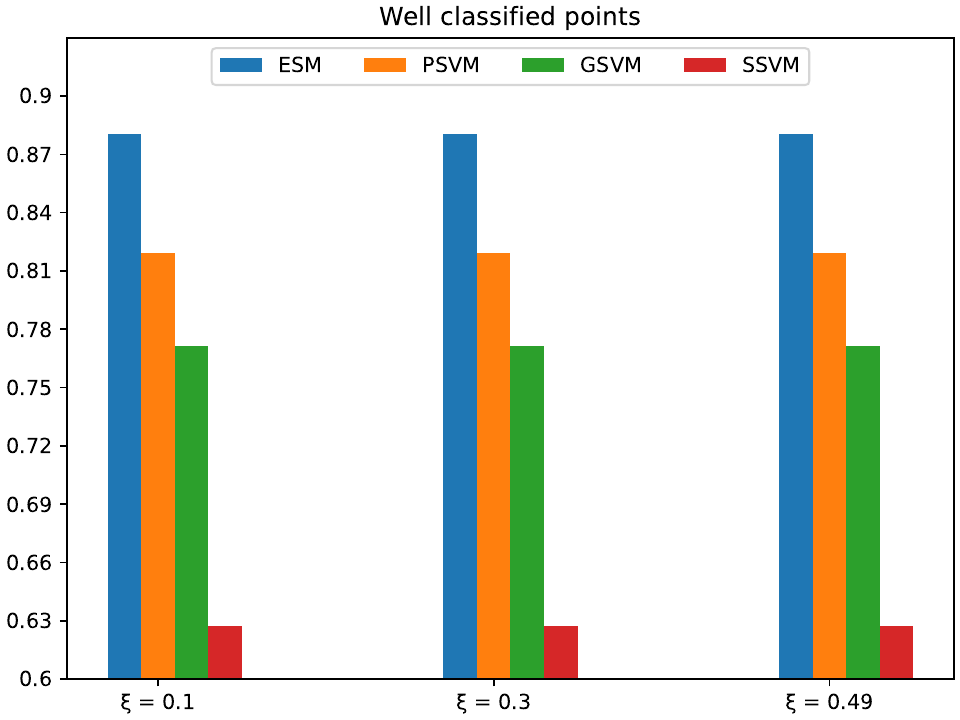}
\end{subfigure}
\hspace{0.1cm}
\begin{subfigure}{0.3\textwidth}
\includegraphics[scale=0.25]{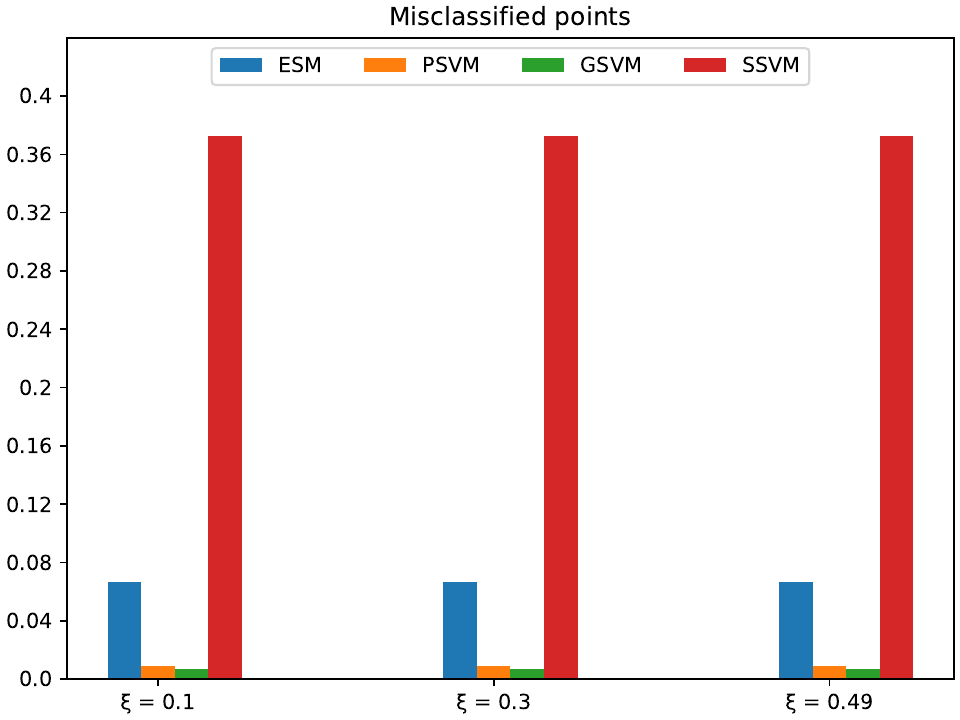}
\end{subfigure}
\hspace{0.1cm}
\begin{subfigure}{0.3\textwidth}
\includegraphics[scale=0.25]{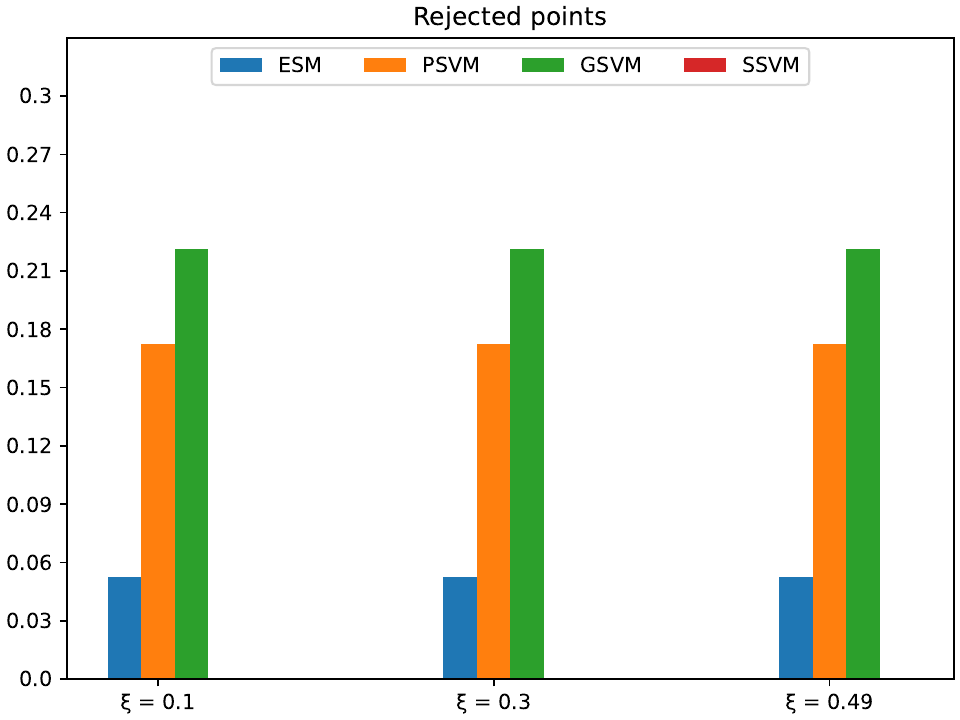}
\end{subfigure}
\end{center}
\caption{Numerical results for the real dataset Breast Wisconsin}
\label{fig:breastwisc}
\end{figure}

\begin{figure}[h!]
\begin{center}
\begin{subfigure}{0.3\textwidth}
\includegraphics[scale=0.25]{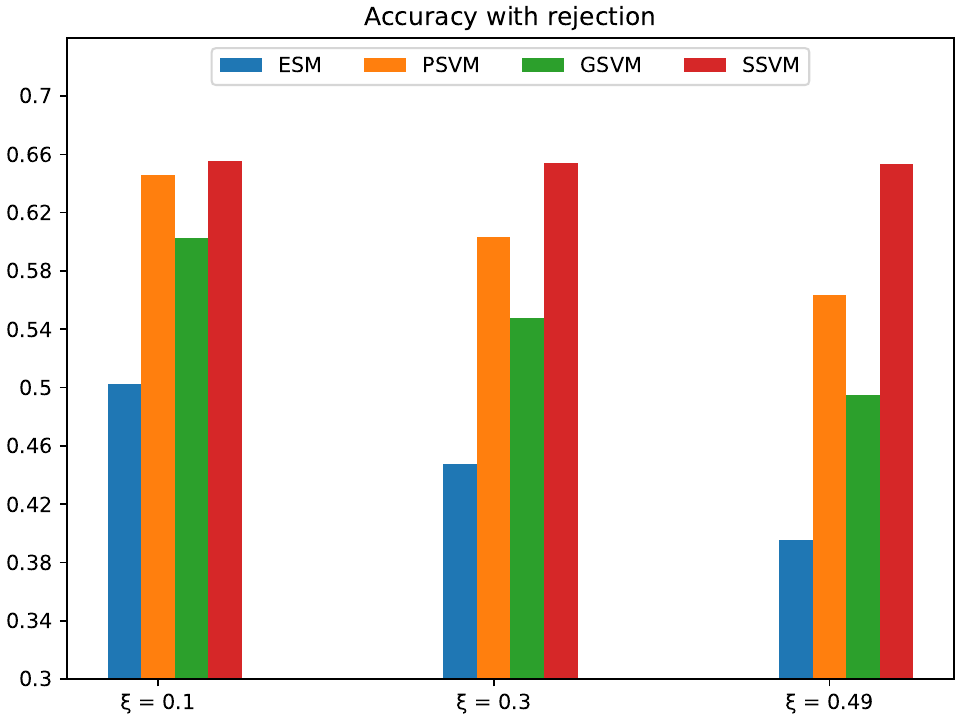}
\end{subfigure}
\hspace{1cm}
\begin{subfigure}{0.3\textwidth}
\includegraphics[scale=0.25]{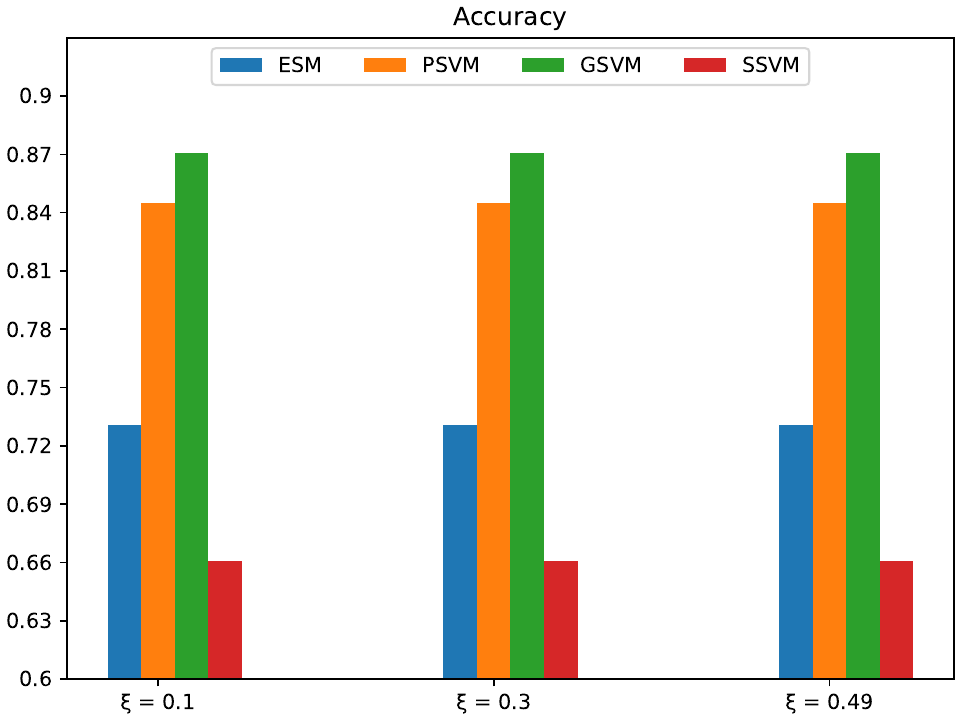}
\end{subfigure}
\vspace{0.3cm}
\begin{subfigure}{0.3\textwidth}
\includegraphics[scale=0.25]{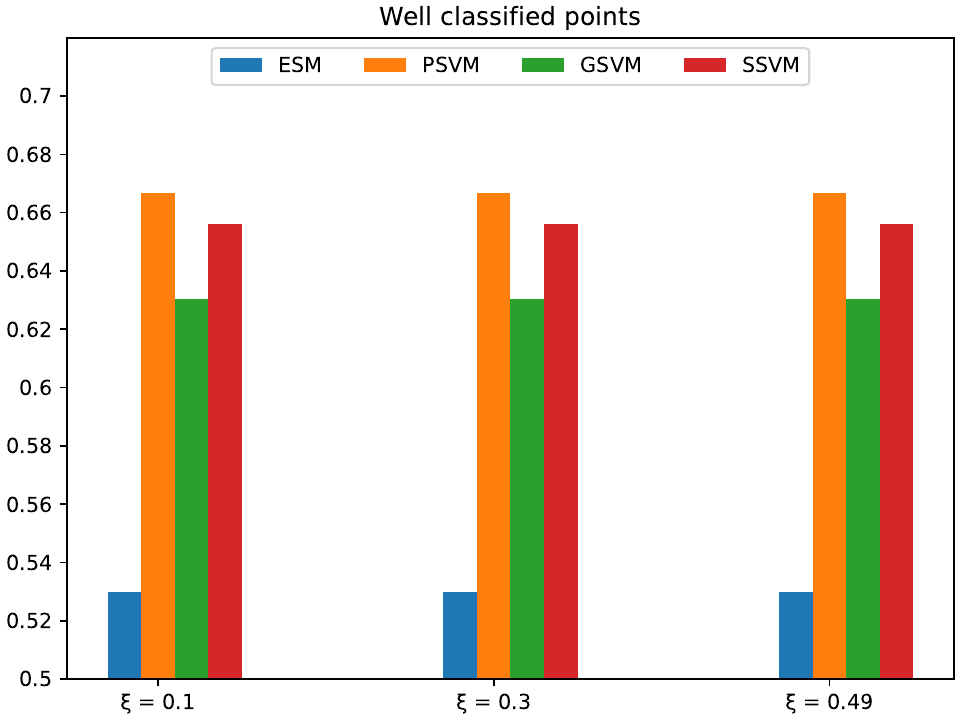}
\end{subfigure}
\hspace{0.1cm}
\begin{subfigure}{0.3\textwidth}
\includegraphics[scale=0.25]{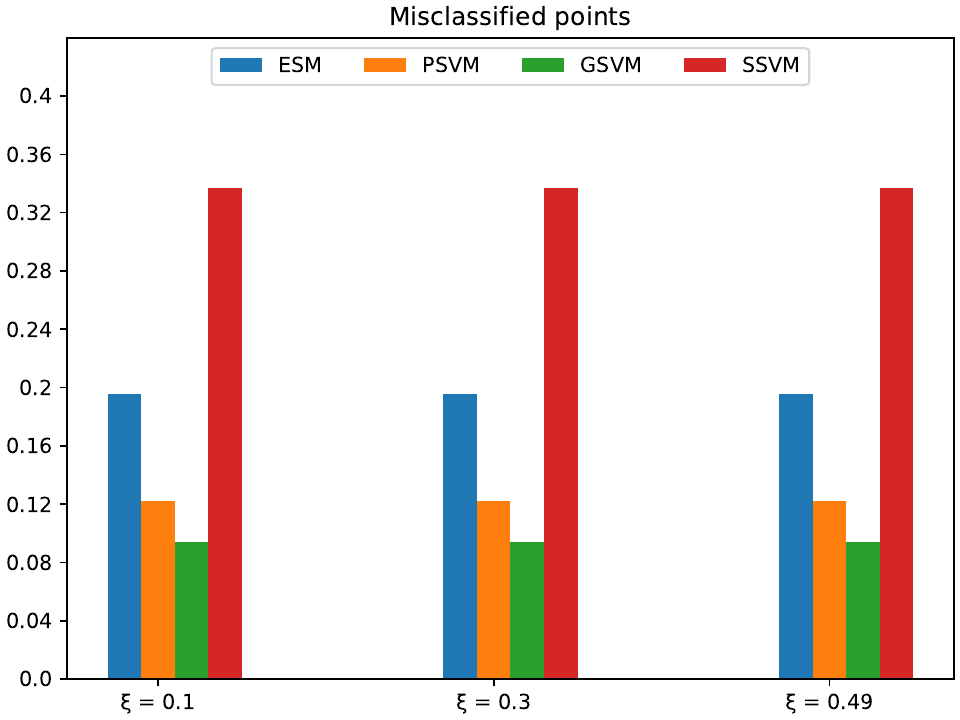}
\end{subfigure}
\hspace{0.1cm}
\begin{subfigure}{0.3\textwidth}
\includegraphics[scale=0.25]{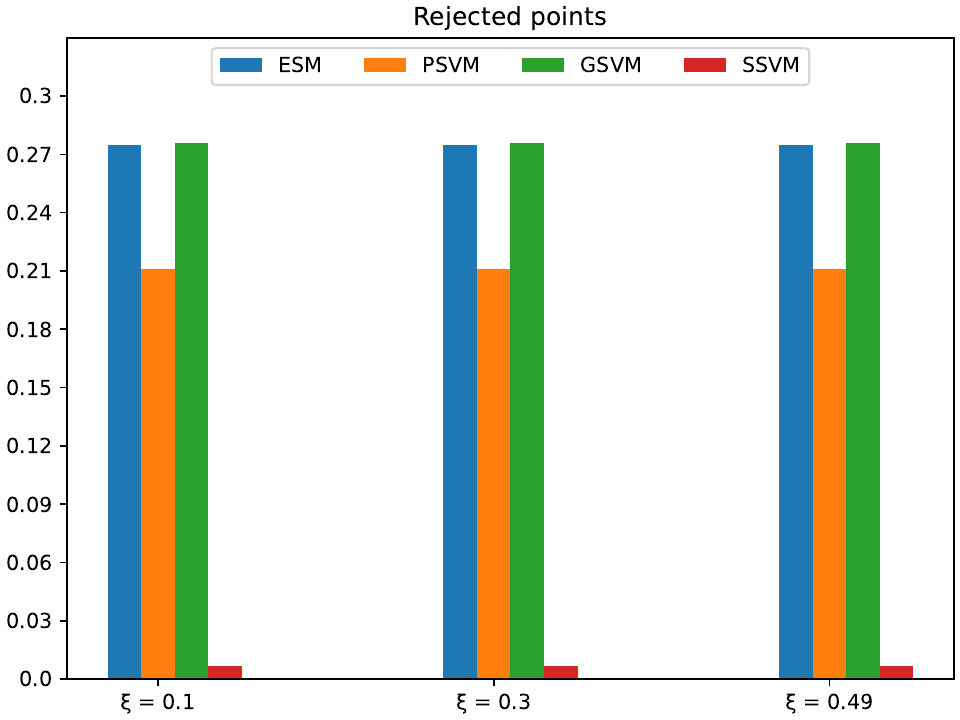}
\end{subfigure}
\end{center}
\caption{Numerical results for the real dataset Diabetes}
\label{fig:diab}
\end{figure}

\begin{figure}[h!]
\begin{center}
\begin{subfigure}{0.3\textwidth}
\includegraphics[scale=0.25]{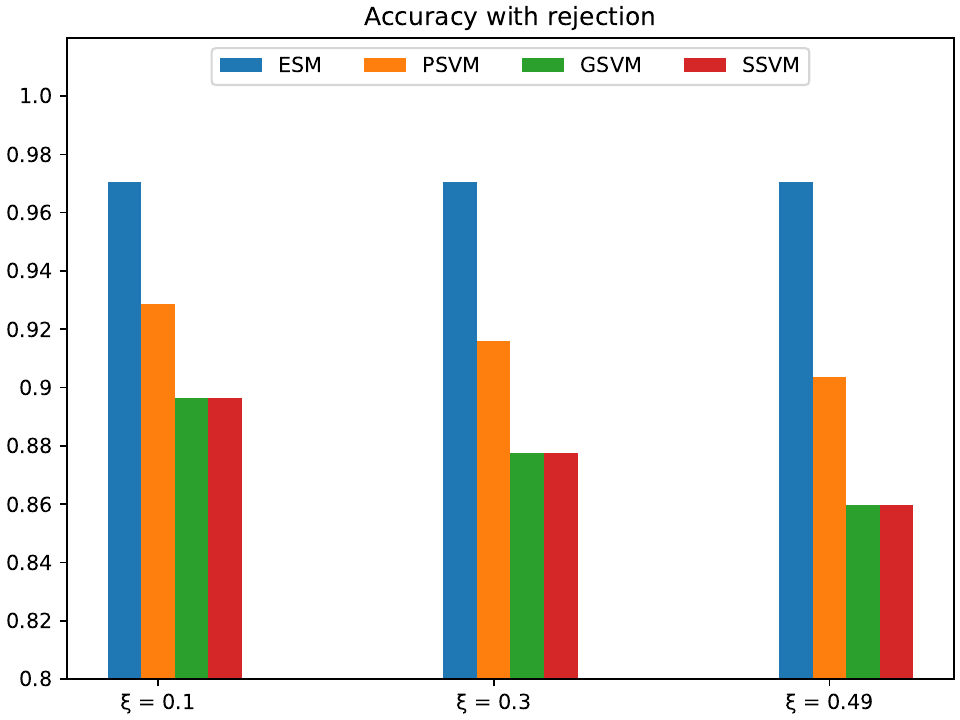}
\end{subfigure}
\hspace{1cm}
\begin{subfigure}{0.3\textwidth}
\includegraphics[scale=0.25]{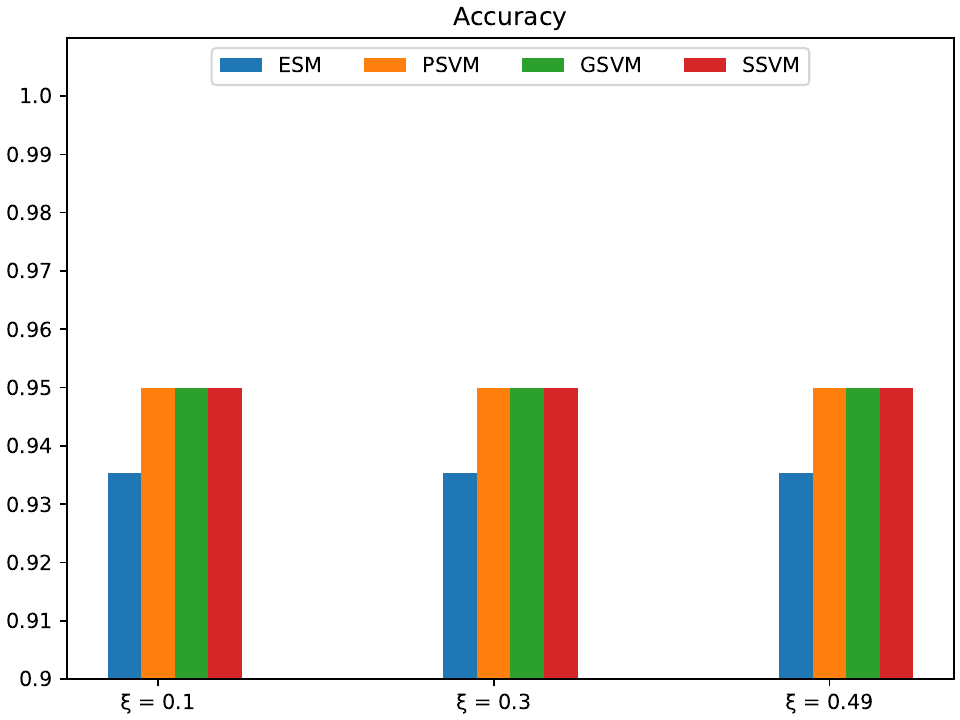}
\end{subfigure}
\vspace{0.3cm}
\begin{subfigure}{0.3\textwidth}
\includegraphics[scale=0.25]{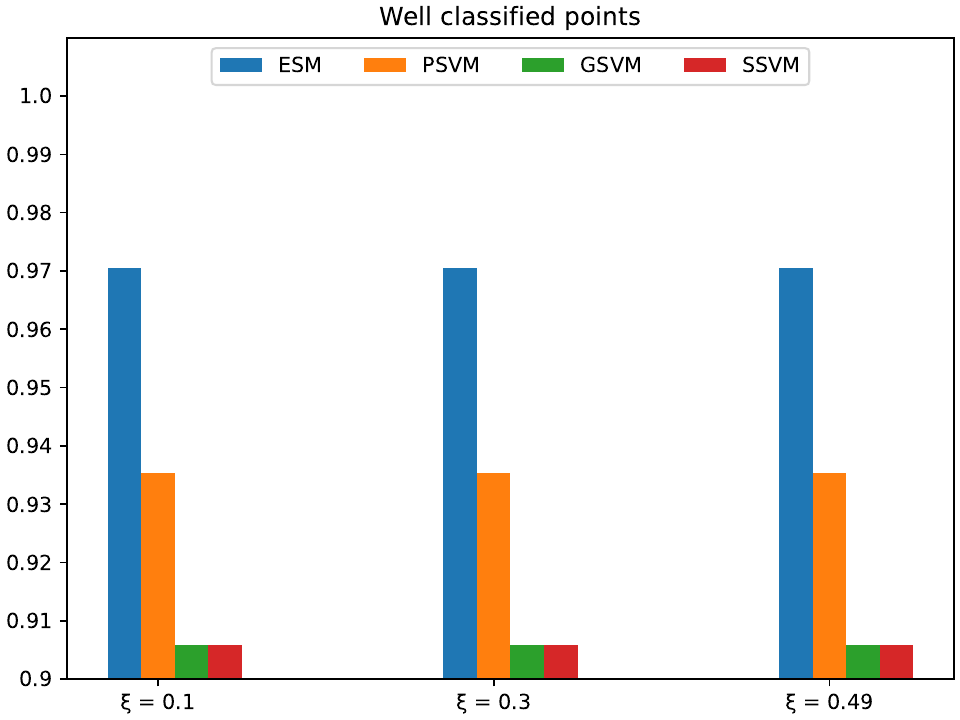}
\end{subfigure}
\hspace{0.1cm}
\begin{subfigure}{0.3\textwidth}
\includegraphics[scale=0.25]{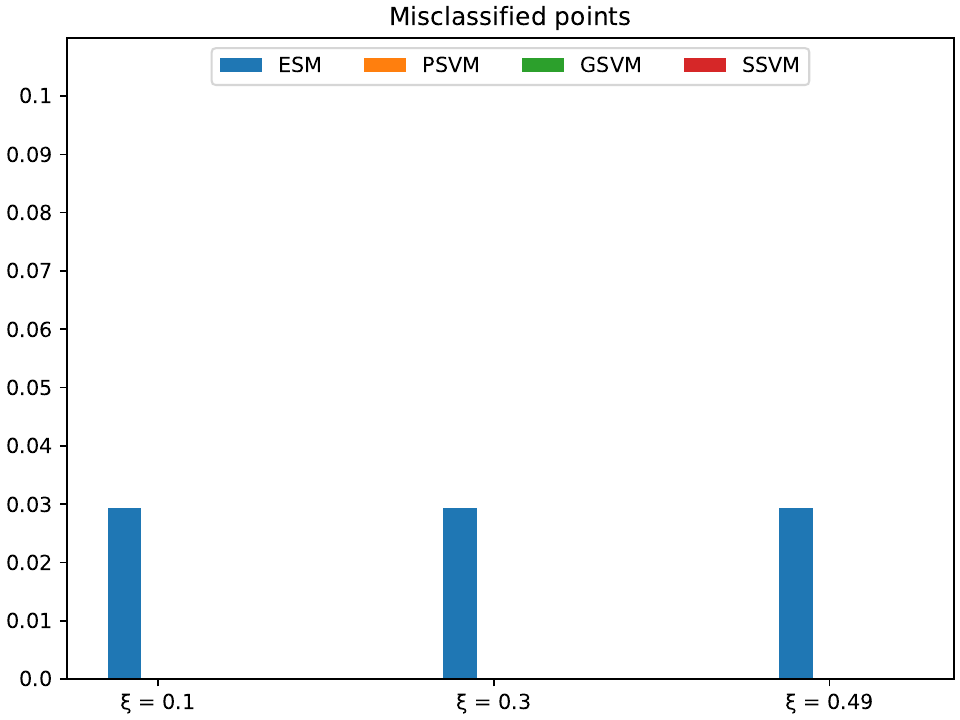}
\end{subfigure}
\hspace{0.1cm}
\begin{subfigure}{0.3\textwidth}
\includegraphics[scale=0.25]{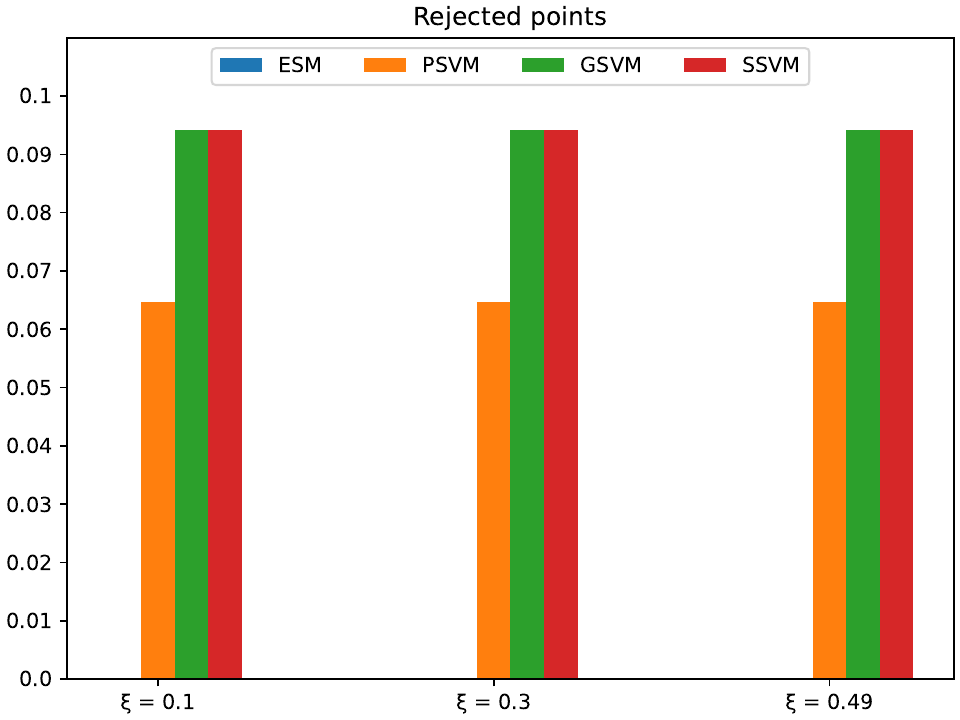}
\end{subfigure}
\end{center}
\caption{Numerical results for the real dataset Divorce}
\label{fig:divor}
\end{figure}

\begin{figure}[h!]
\begin{center}
\begin{subfigure}{0.3\textwidth}
\includegraphics[scale=0.25]{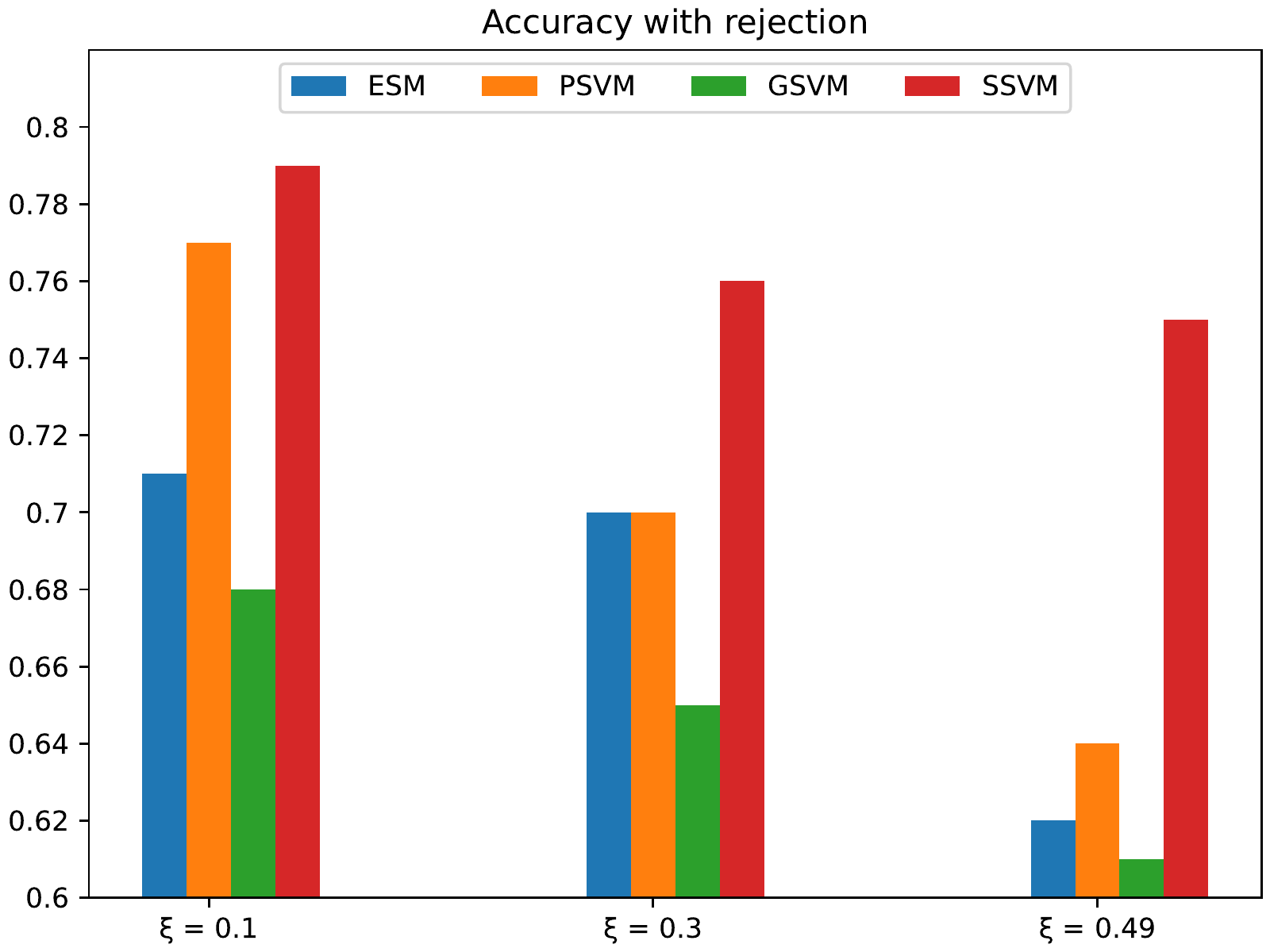}
\end{subfigure}
\hspace{1cm}
\begin{subfigure}{0.3\textwidth}
\includegraphics[scale=0.25]{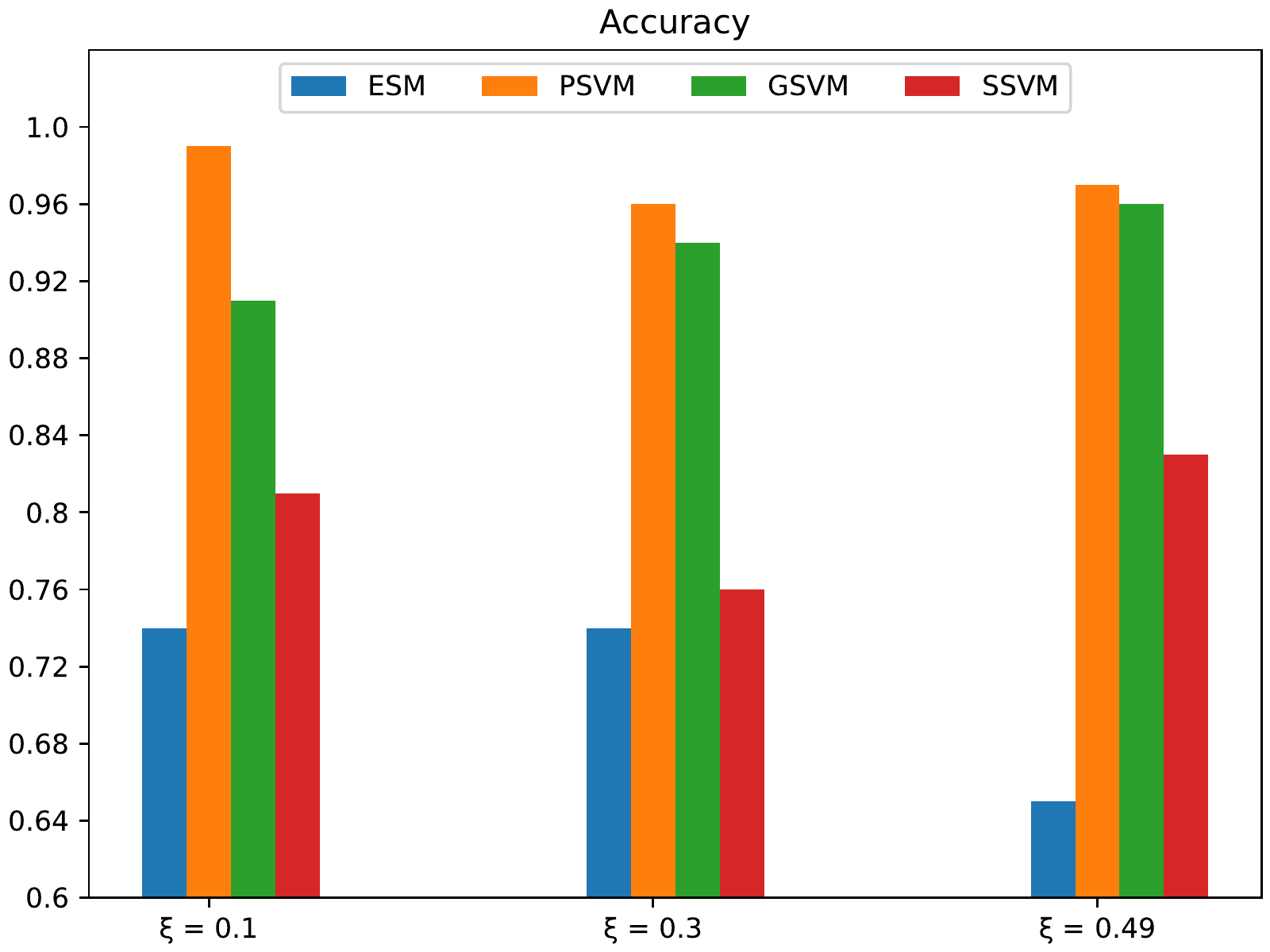}
\end{subfigure}
\vspace{0.3cm}
\begin{subfigure}{0.3\textwidth}
\includegraphics[scale=0.25]{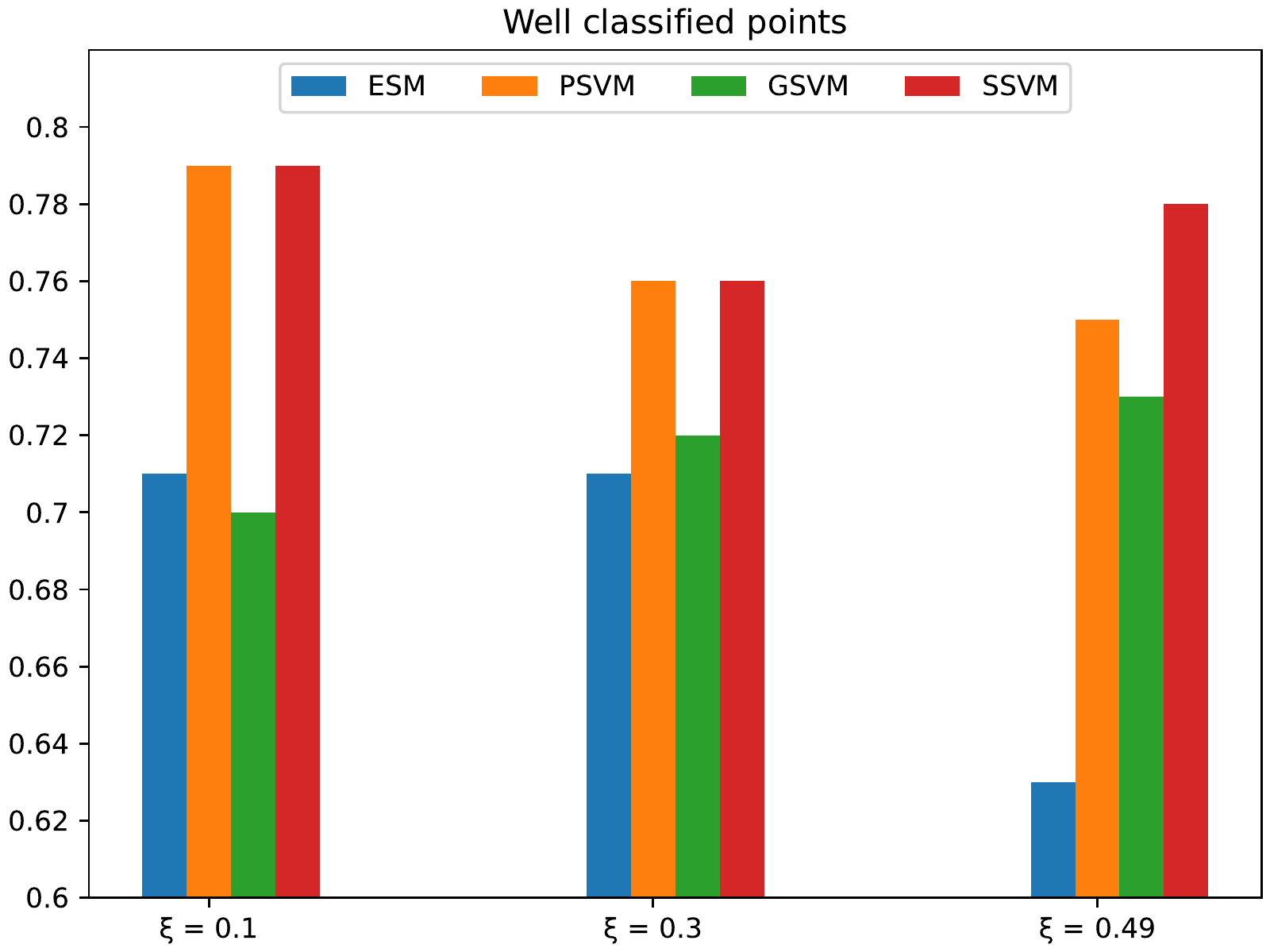}
\end{subfigure}
\hspace{0.1cm}
\begin{subfigure}{0.3\textwidth}
\includegraphics[scale=0.25]{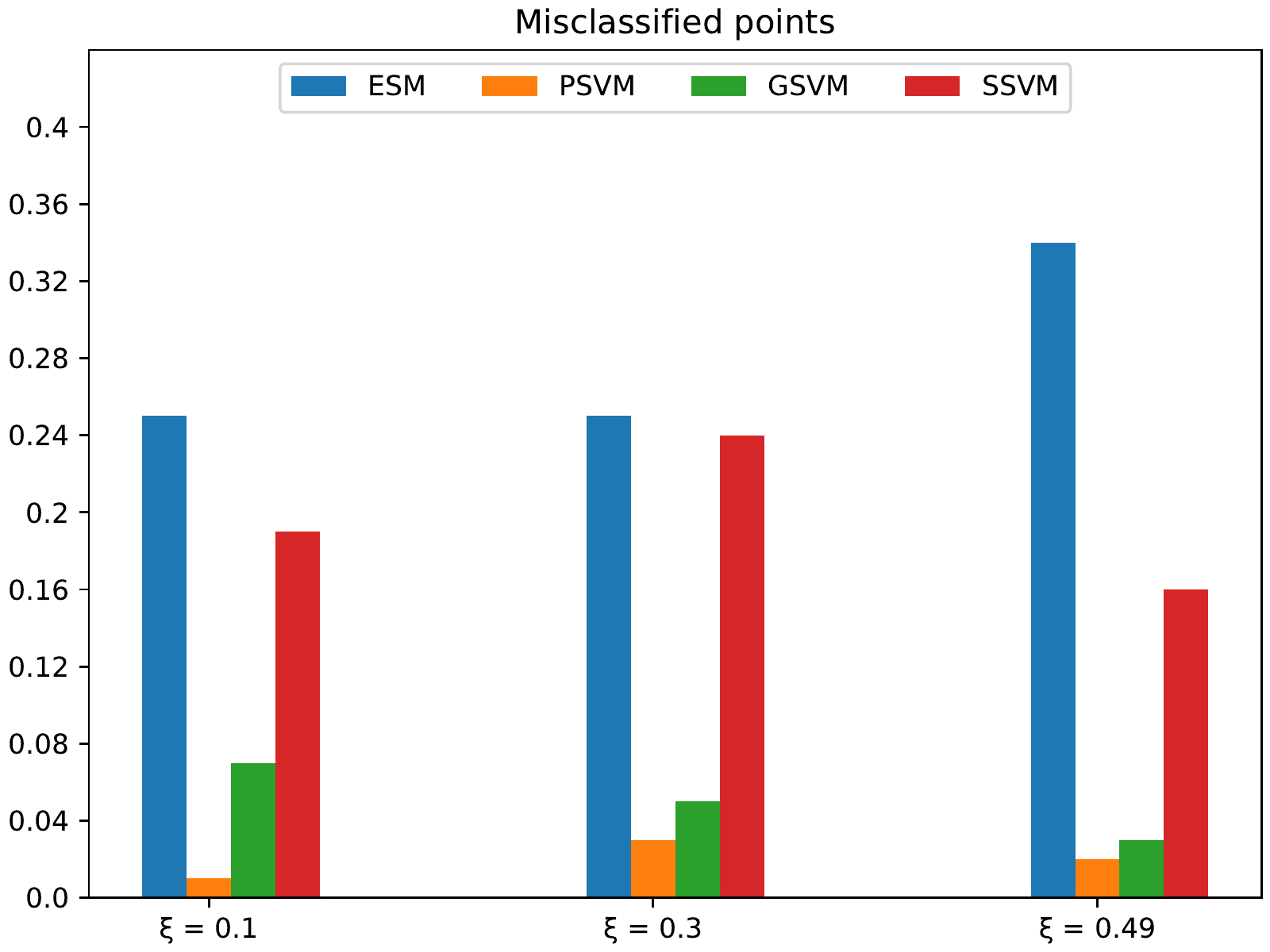}
\end{subfigure}
\hspace{0.1cm}
\begin{subfigure}{0.3\textwidth}
\includegraphics[scale=0.25]{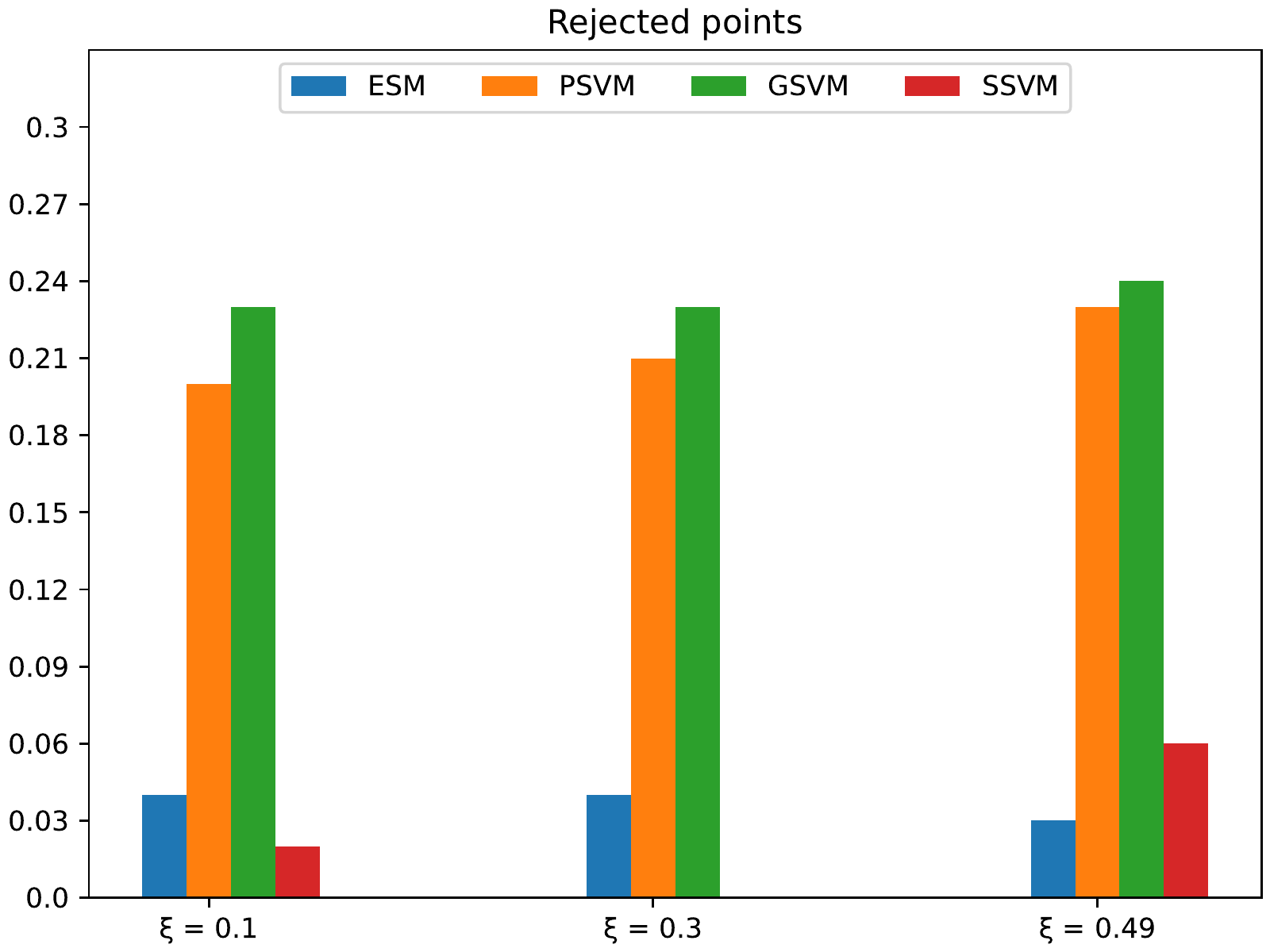}
\end{subfigure}
\end{center}
\caption{Numerical results for the real dataset Fertility}
\label{fig:fertility}
\end{figure}

\begin{figure}[h!]
\begin{center}
\begin{subfigure}{0.3\textwidth}
\includegraphics[scale=0.25]{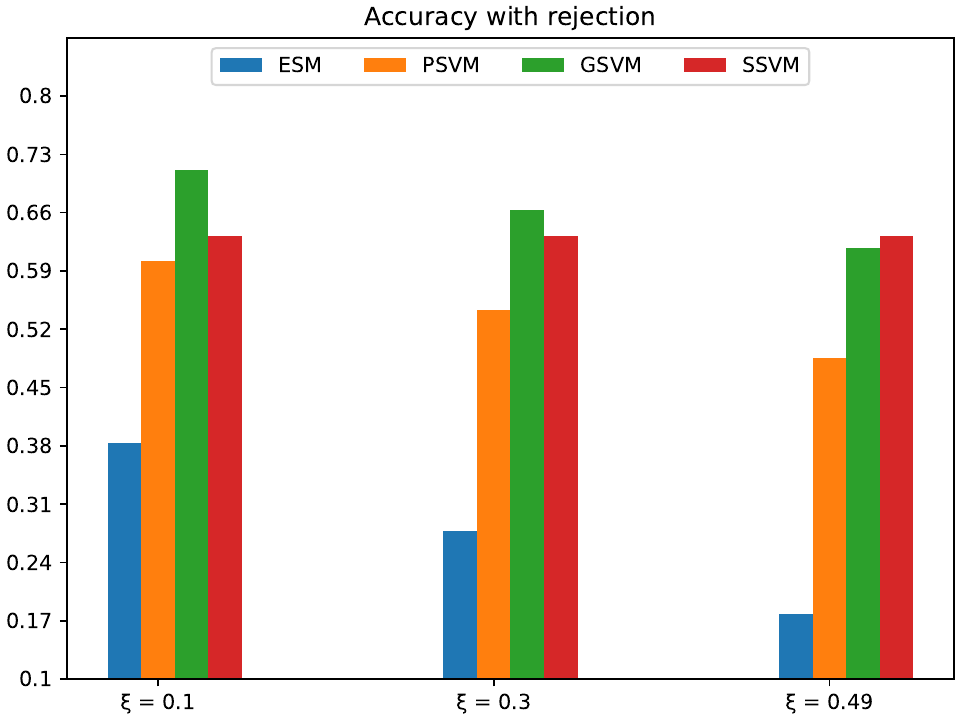}
\end{subfigure}
\hspace{1cm}
\begin{subfigure}{0.3\textwidth}
\includegraphics[scale=0.25]{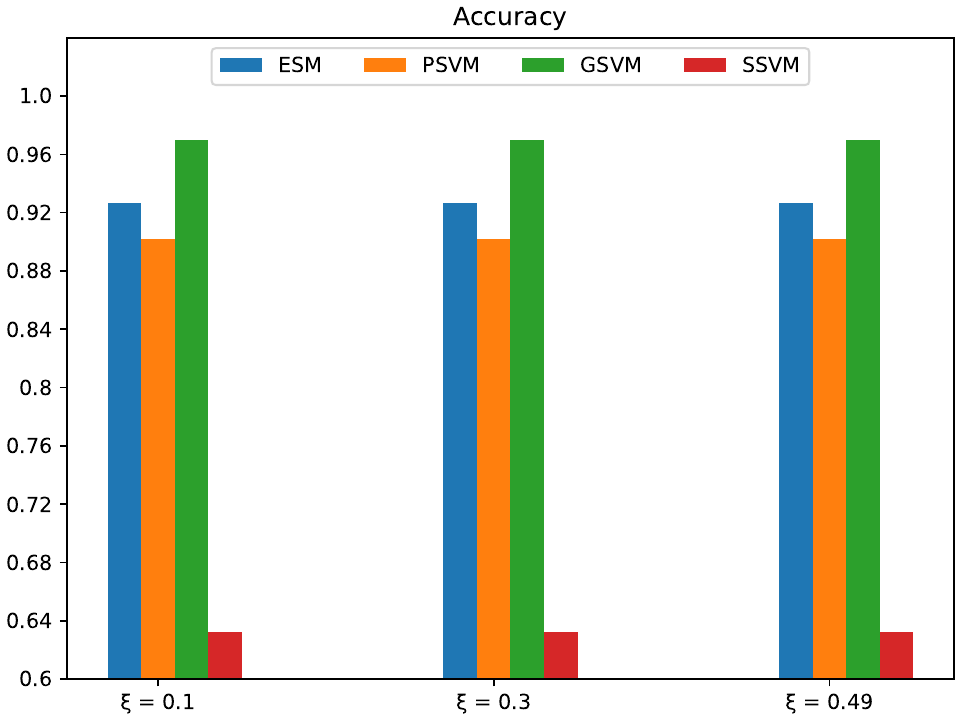}
\end{subfigure}
\vspace{0.3cm}
\begin{subfigure}{0.3\textwidth}
\includegraphics[scale=0.25]{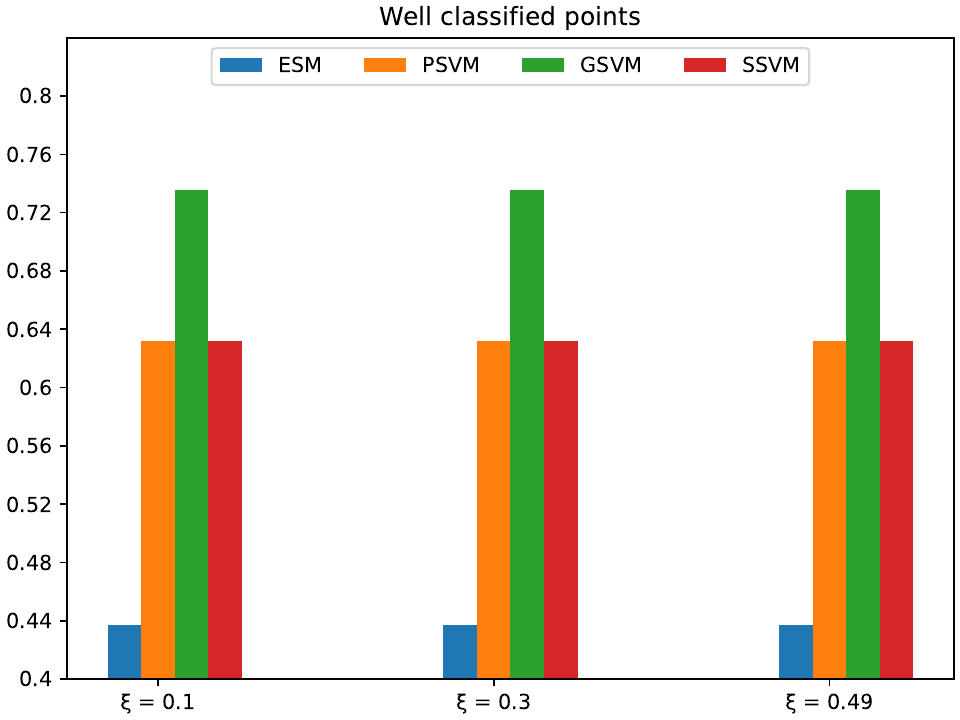}
\end{subfigure}
\hspace{0.1cm}
\begin{subfigure}{0.3\textwidth}
\includegraphics[scale=0.25]{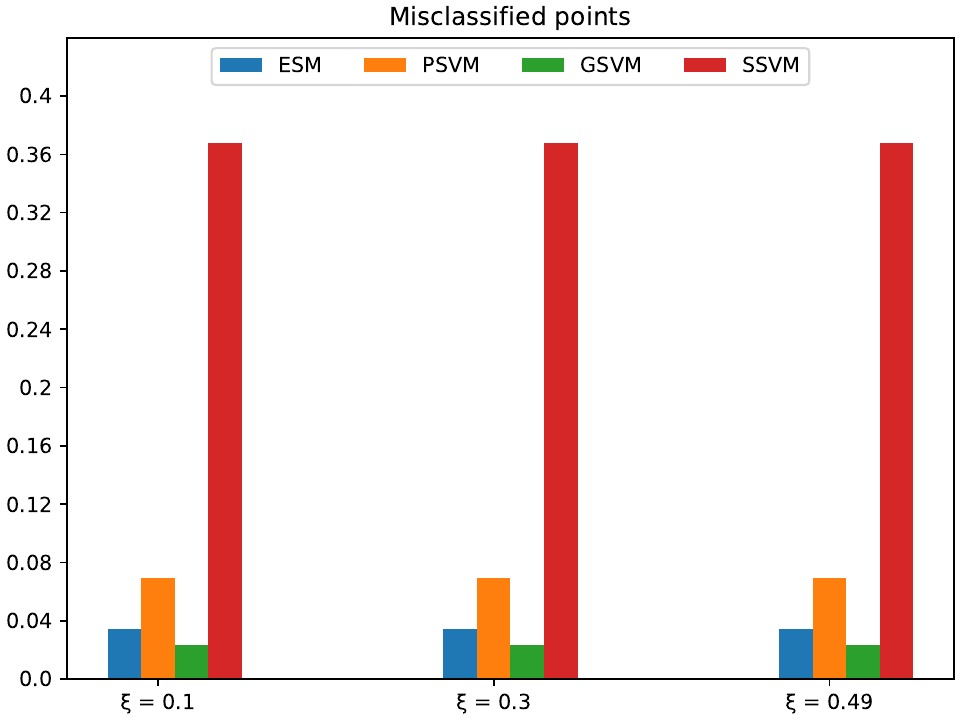}
\end{subfigure}
\hspace{0.1cm}
\begin{subfigure}{0.3\textwidth}
\includegraphics[scale=0.25]{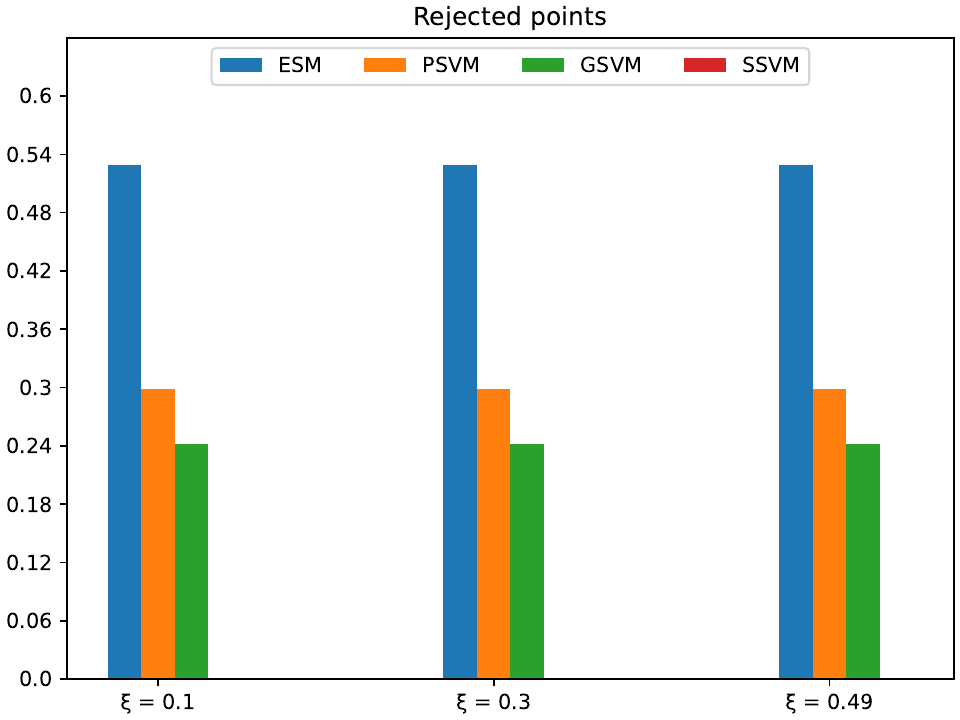}
\end{subfigure}
\end{center}
\caption{Numerical results for the real dataset Flowmeters}
\label{fig:flowm}
\end{figure}

\begin{figure}[h!]
\begin{center}
\begin{subfigure}{0.3\textwidth}
\includegraphics[scale=0.25]{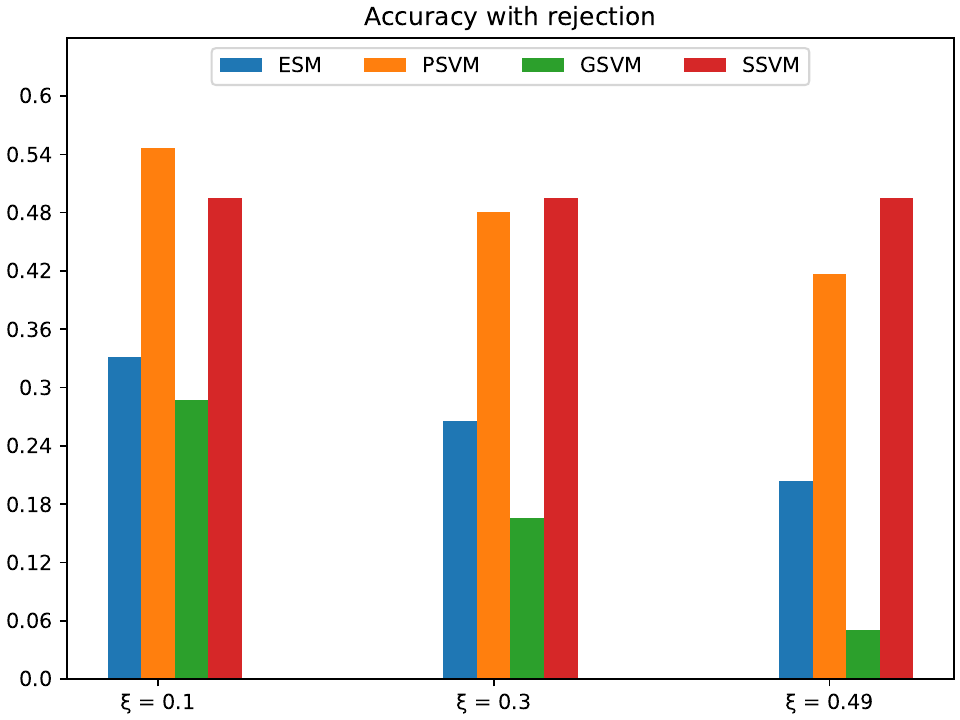}
\end{subfigure}
\hspace{1cm}
\begin{subfigure}{0.3\textwidth}
\includegraphics[scale=0.25]{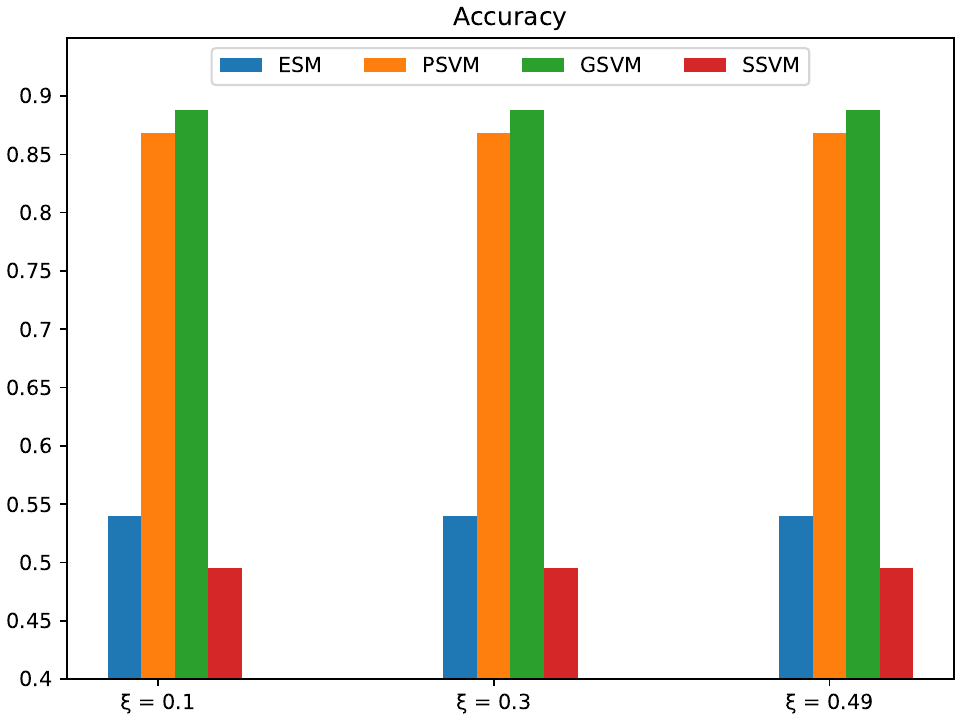}
\end{subfigure}
\vspace{0.3cm}
\begin{subfigure}{0.3\textwidth}
\includegraphics[scale=0.25]{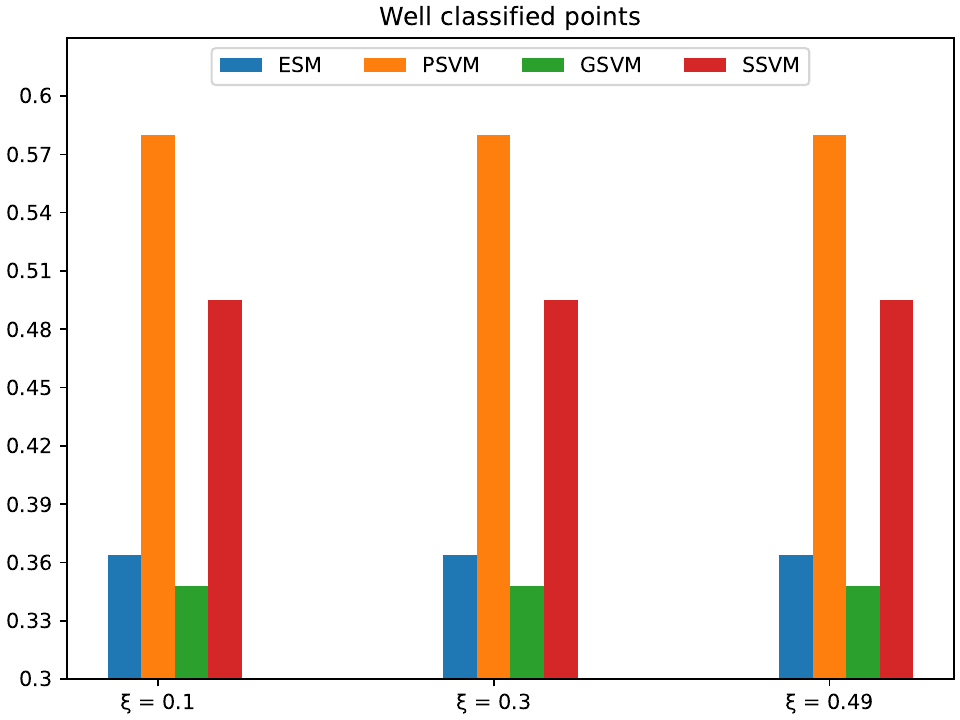}
\end{subfigure}
\hspace{0.1cm}
\begin{subfigure}{0.3\textwidth}
\includegraphics[scale=0.25]{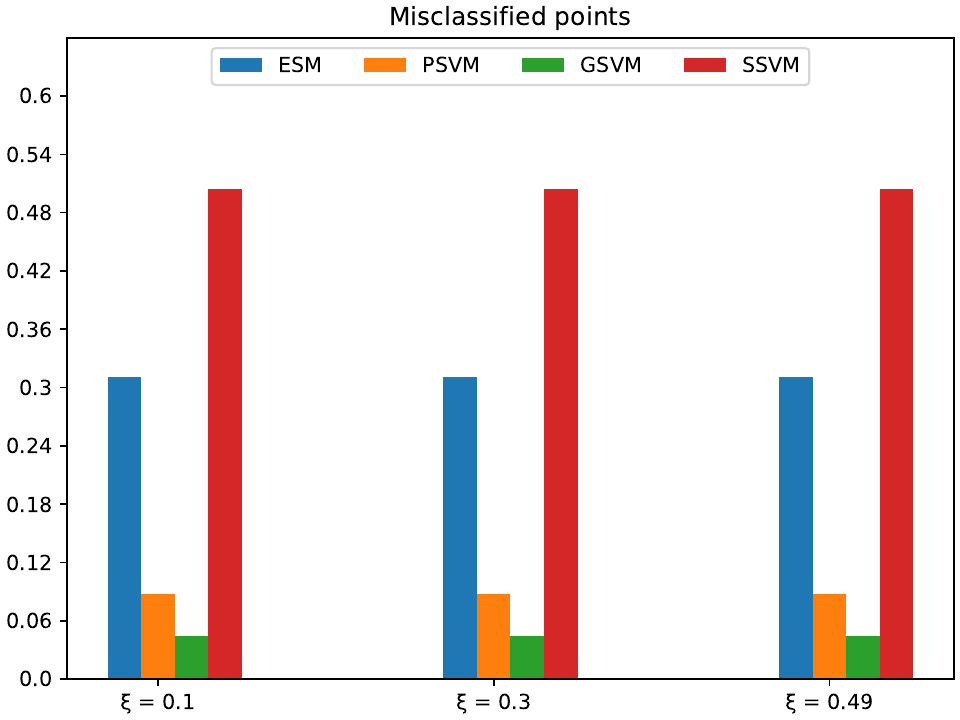}
\end{subfigure}
\hspace{0.1cm}
\begin{subfigure}{0.3\textwidth}
\includegraphics[scale=0.25]{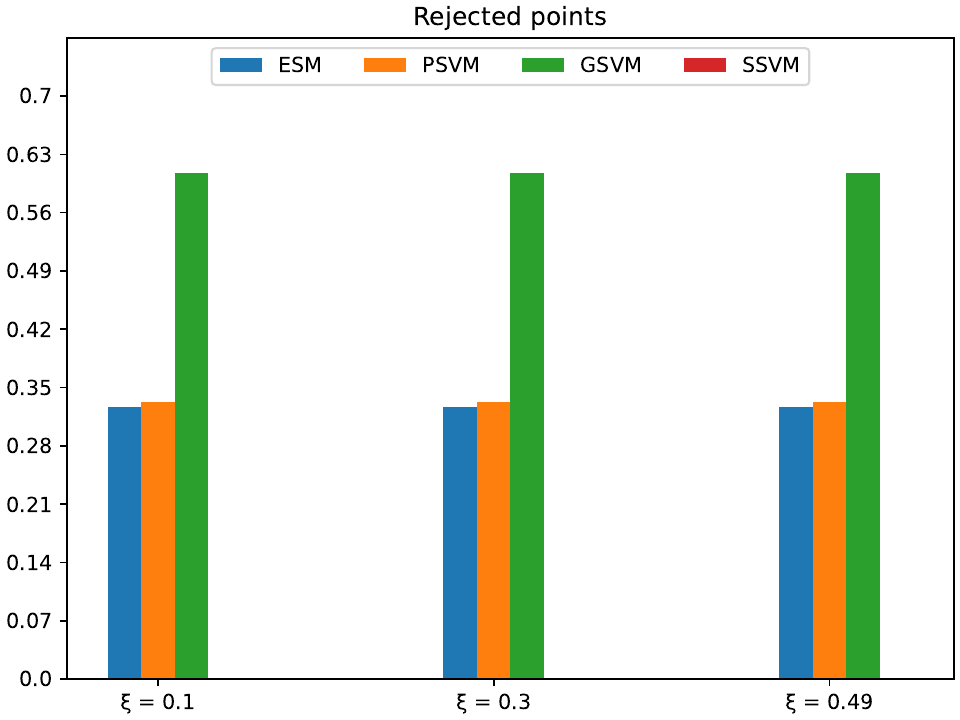}
\end{subfigure}
\end{center}
\caption{Numerical results for the real dataset Gallstone}
\label{fig:gallst}
\end{figure}

\begin{figure}[h!]
\begin{center}
\begin{subfigure}{0.3\textwidth}
\includegraphics[scale=0.25]{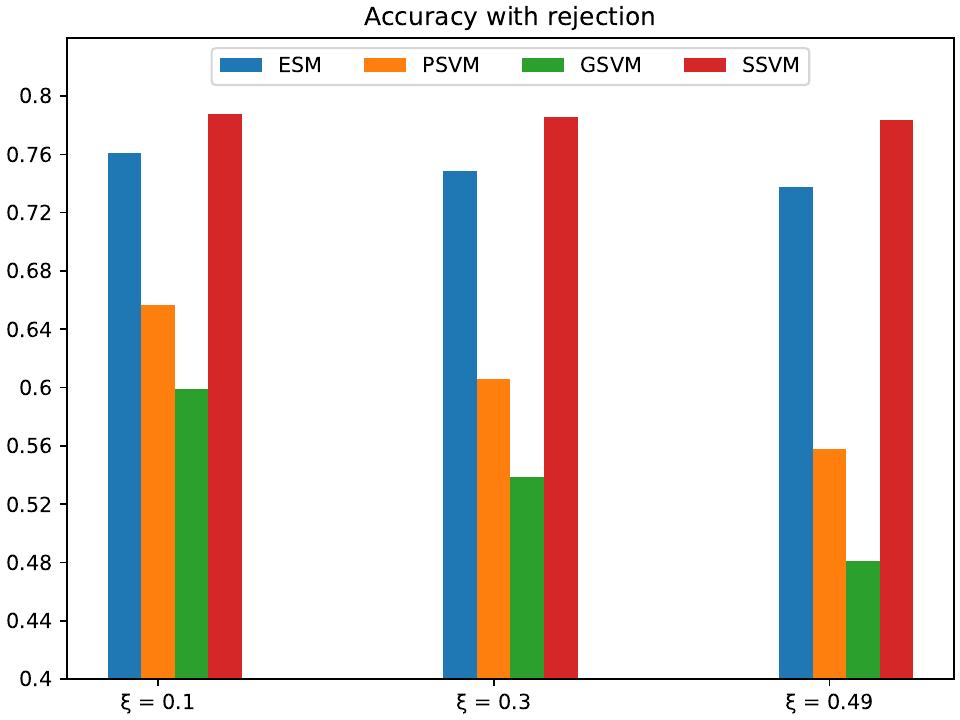}
\end{subfigure}
\hspace{1cm}
\begin{subfigure}{0.3\textwidth}
\includegraphics[scale=0.25]{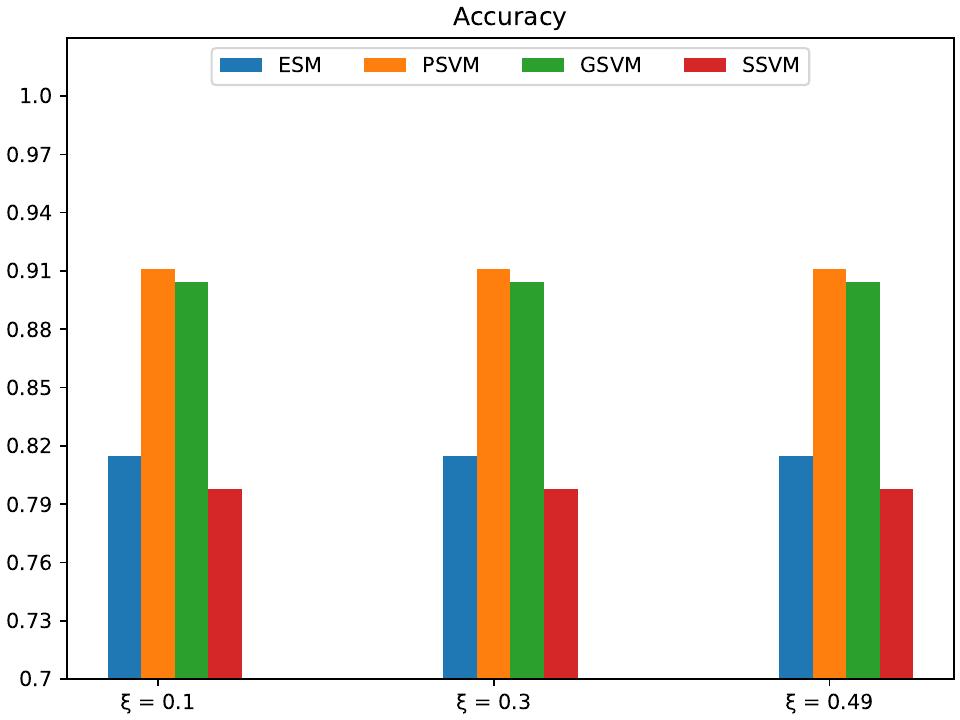}
\end{subfigure}
\vspace{0.3cm}
\begin{subfigure}{0.3\textwidth}
\includegraphics[scale=0.25]{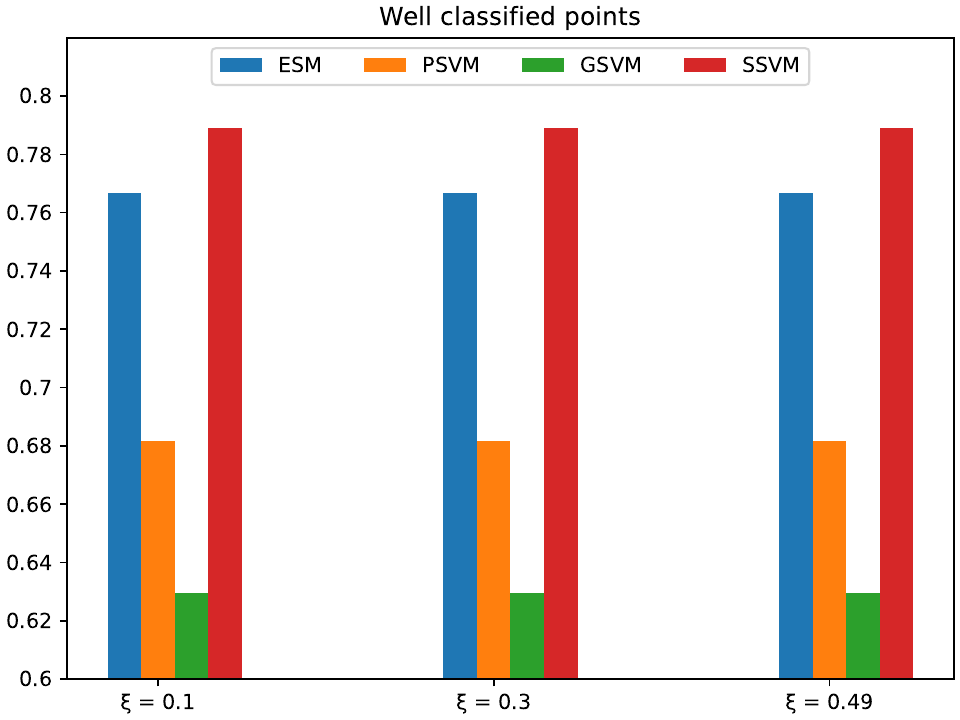}
\end{subfigure}
\hspace{0.1cm}
\begin{subfigure}{0.3\textwidth}
\includegraphics[scale=0.25]{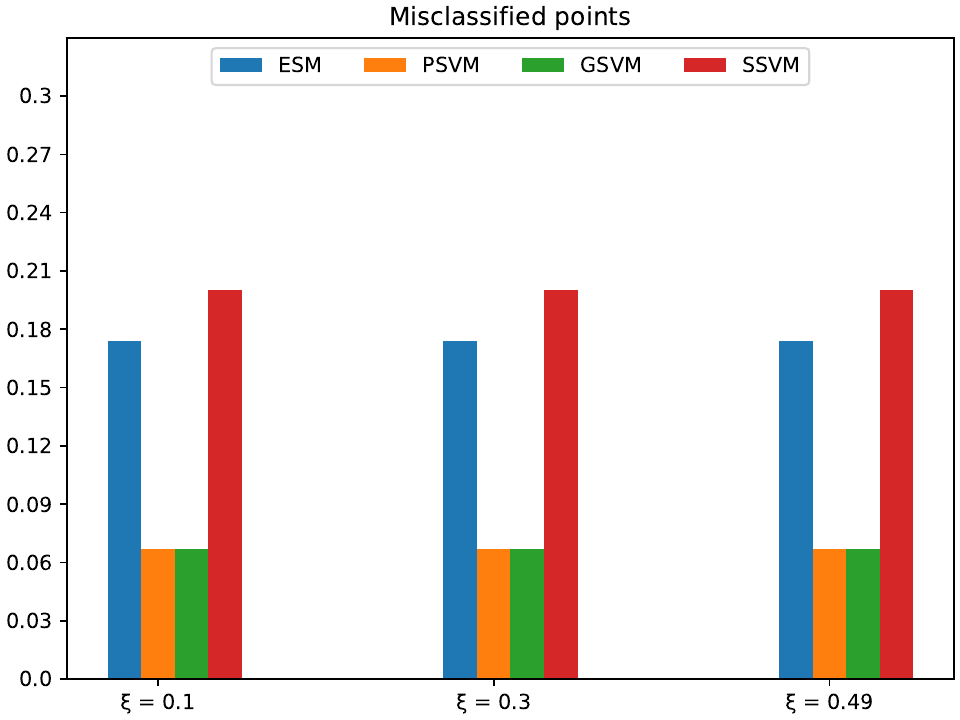}
\end{subfigure}
\hspace{0.1cm}
\begin{subfigure}{0.3\textwidth}
\includegraphics[scale=0.25]{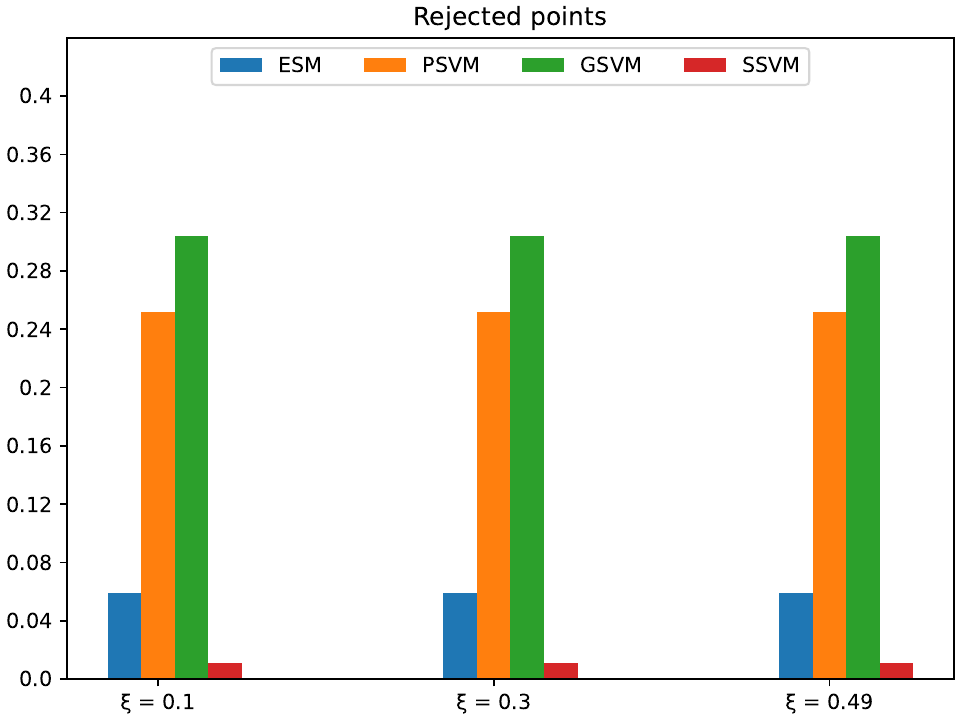}
\end{subfigure}
\end{center}
\caption{Numerical results for the real dataset Heart}
\label{fig:heart}
\end{figure}

\begin{figure}[h!]
\begin{center}
\begin{subfigure}{0.3\textwidth}
\includegraphics[scale=0.25]{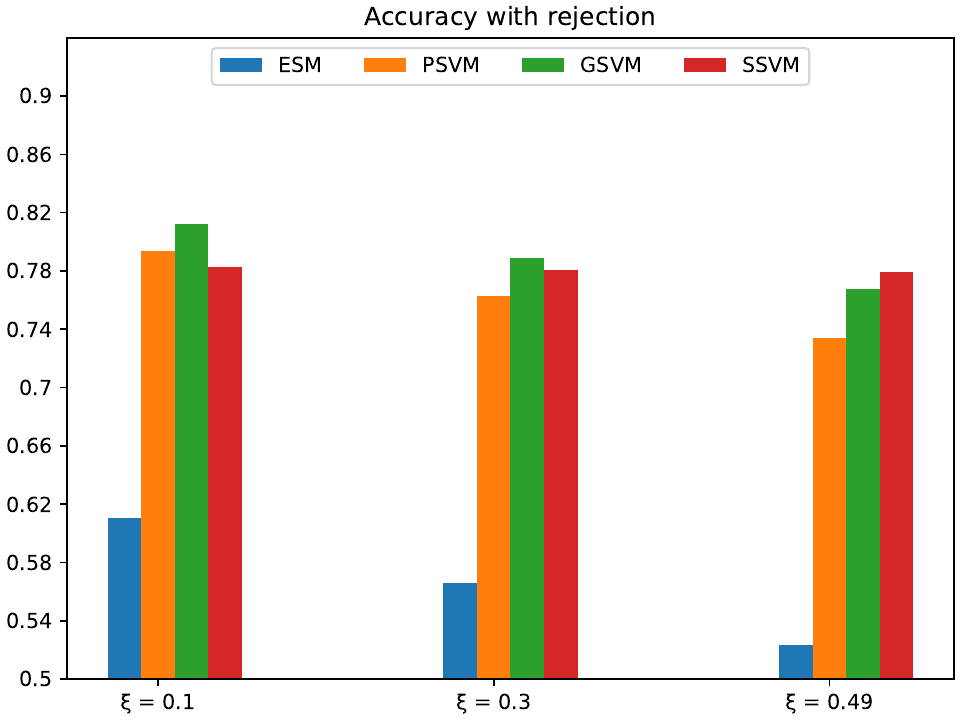}
\end{subfigure}
\hspace{1cm}
\begin{subfigure}{0.3\textwidth}
\includegraphics[scale=0.25]{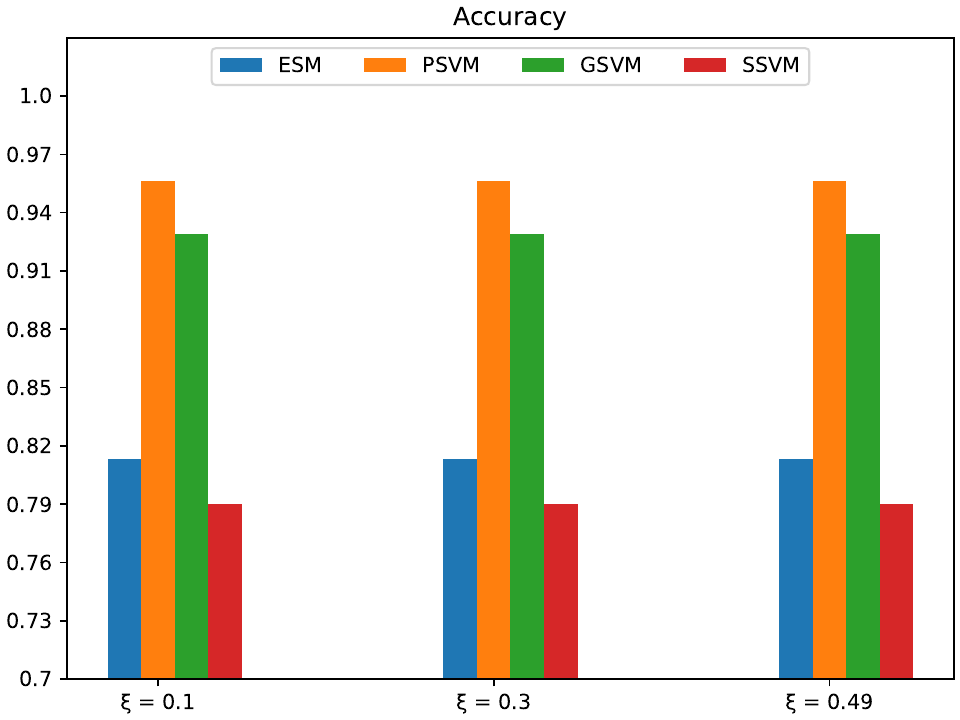}
\end{subfigure}
\vspace{0.3cm}
\begin{subfigure}{0.3\textwidth}
\includegraphics[scale=0.25]{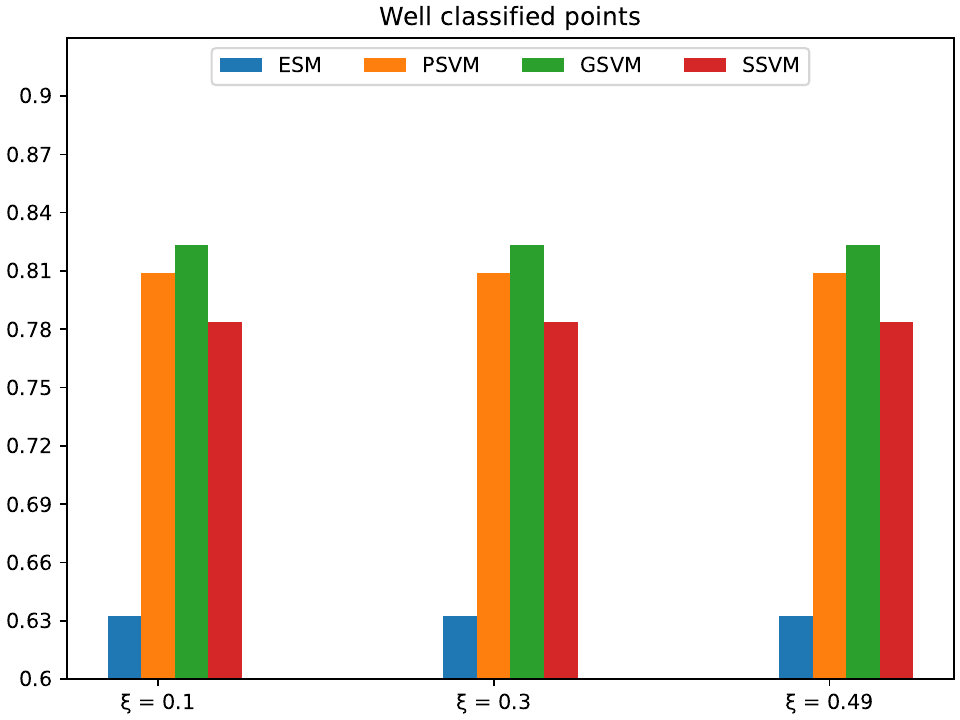}
\end{subfigure}
\hspace{0.1cm}
\begin{subfigure}{0.3\textwidth}
\includegraphics[scale=0.25]{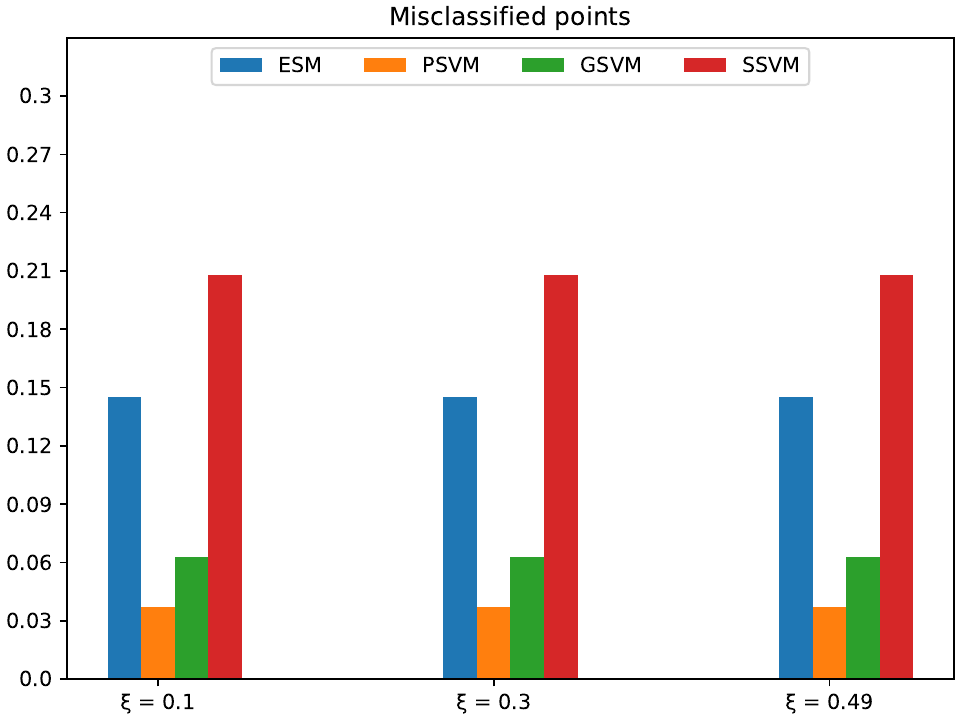}
\end{subfigure}
\hspace{0.1cm}
\begin{subfigure}{0.3\textwidth}
\includegraphics[scale=0.25]{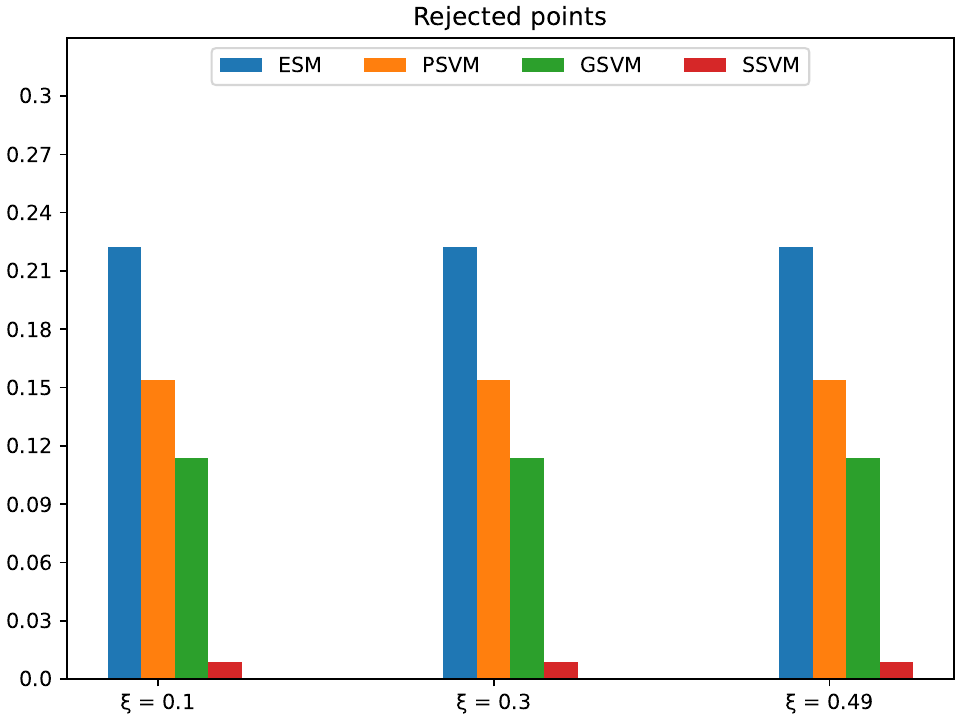}
\end{subfigure}
\end{center}
\caption{Numerical results for the real dataset Ionosphere}
\label{fig:iono}
\end{figure}

\begin{figure}[h!]
\begin{center}
\begin{subfigure}{0.3\textwidth}
\includegraphics[scale=0.25]{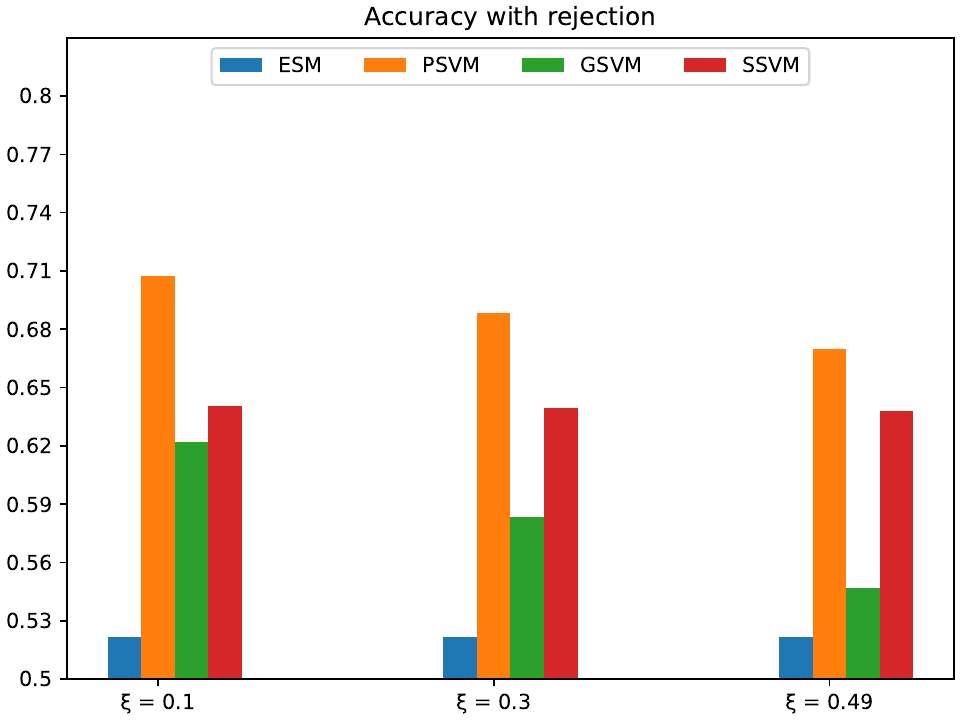}
\end{subfigure}
\hspace{1cm}
\begin{subfigure}{0.3\textwidth}
\includegraphics[scale=0.25]{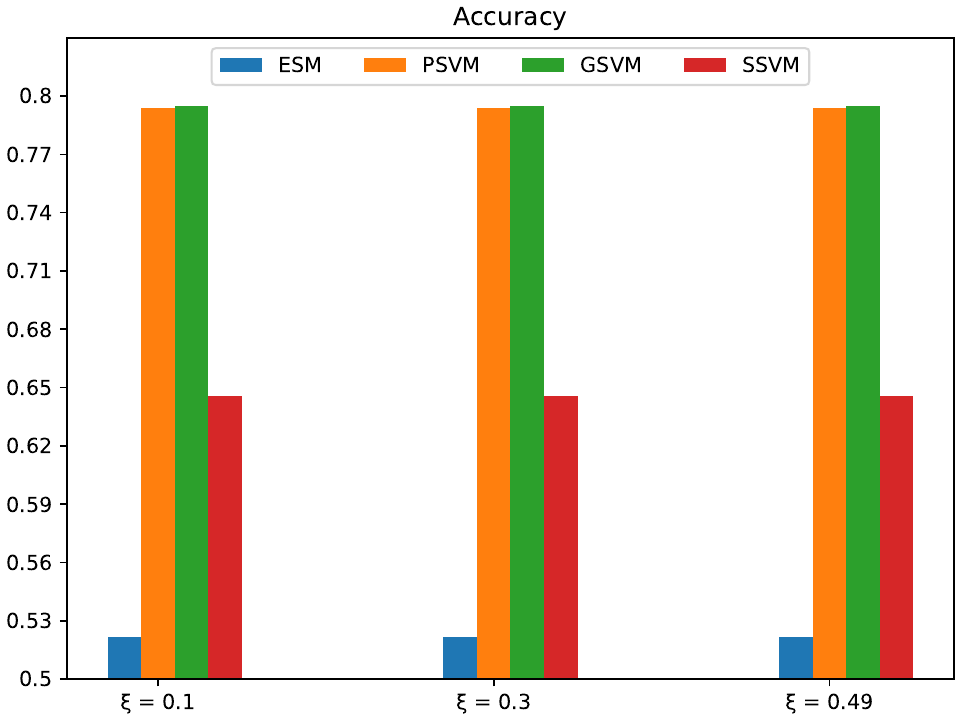}
\end{subfigure}
\vspace{0.3cm}
\begin{subfigure}{0.3\textwidth}
\includegraphics[scale=0.25]{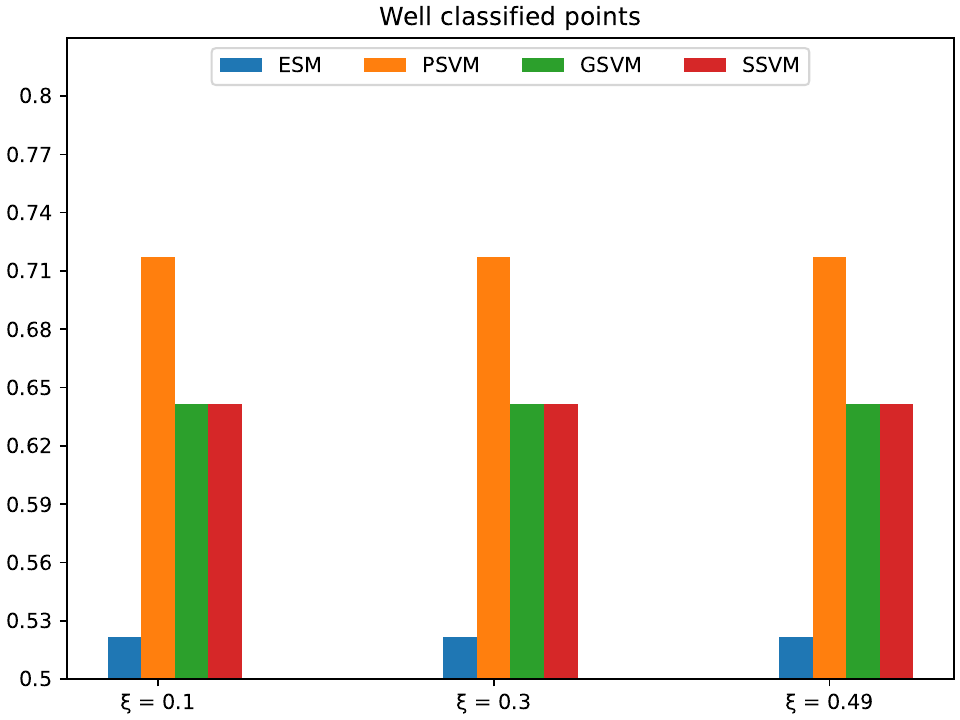}
\end{subfigure}
\hspace{0.1cm}
\begin{subfigure}{0.3\textwidth}
\includegraphics[scale=0.25]{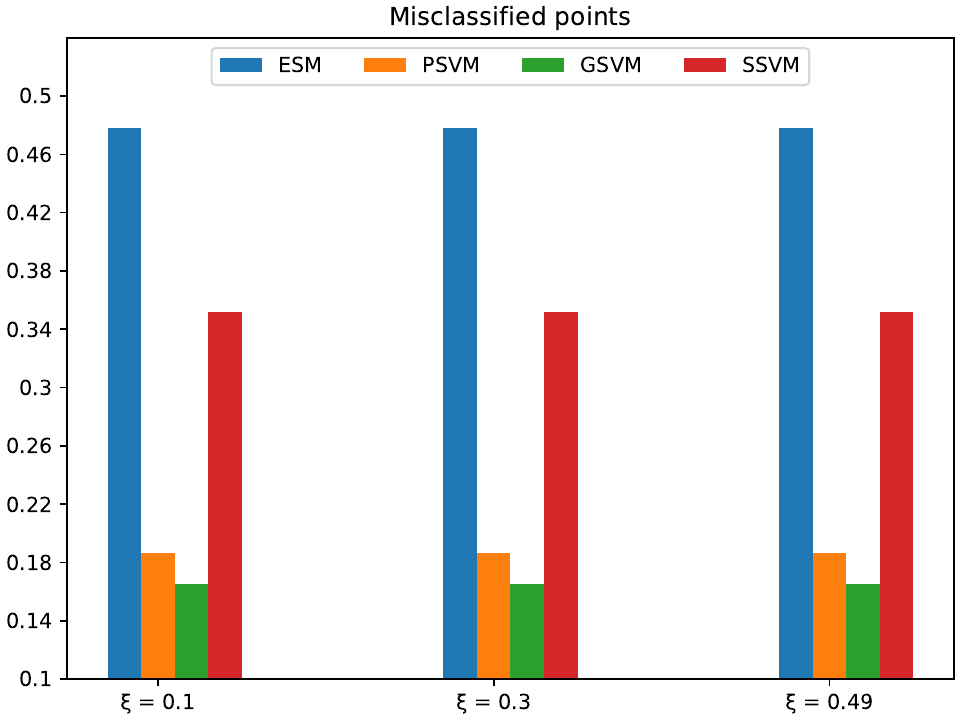}
\end{subfigure}
\hspace{0.1cm}
\begin{subfigure}{0.3\textwidth}
\includegraphics[scale=0.25]{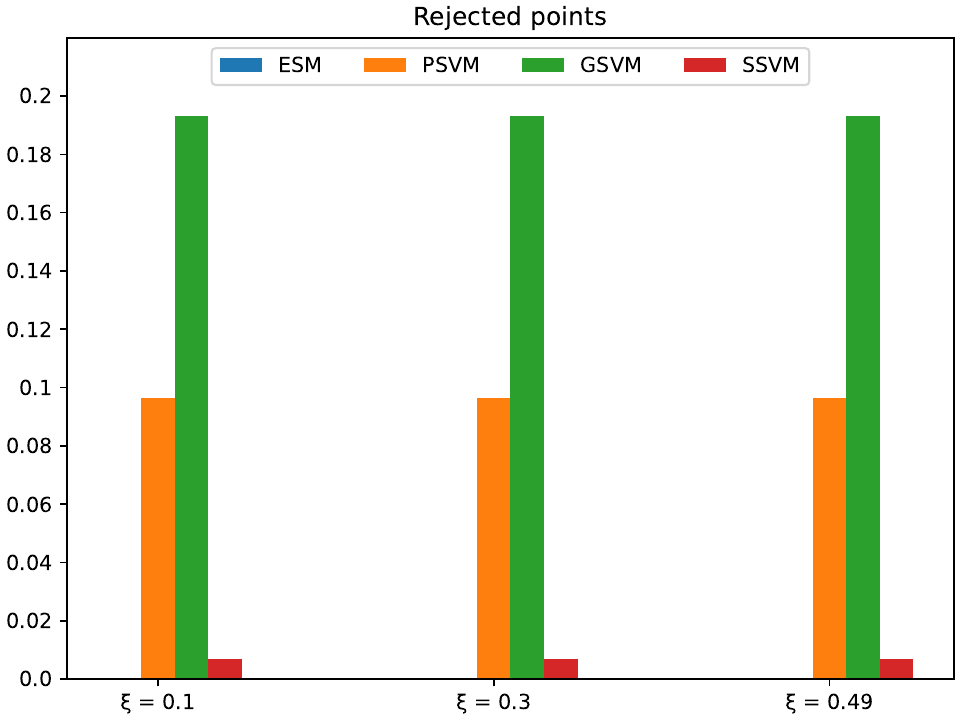}
\end{subfigure}
\end{center}
\caption{Numerical results for the real dataset Liver}
\label{fig:liv}
\end{figure}

\begin{figure}[h!]
\begin{center}
\begin{subfigure}{0.3\textwidth}
\includegraphics[scale=0.25]{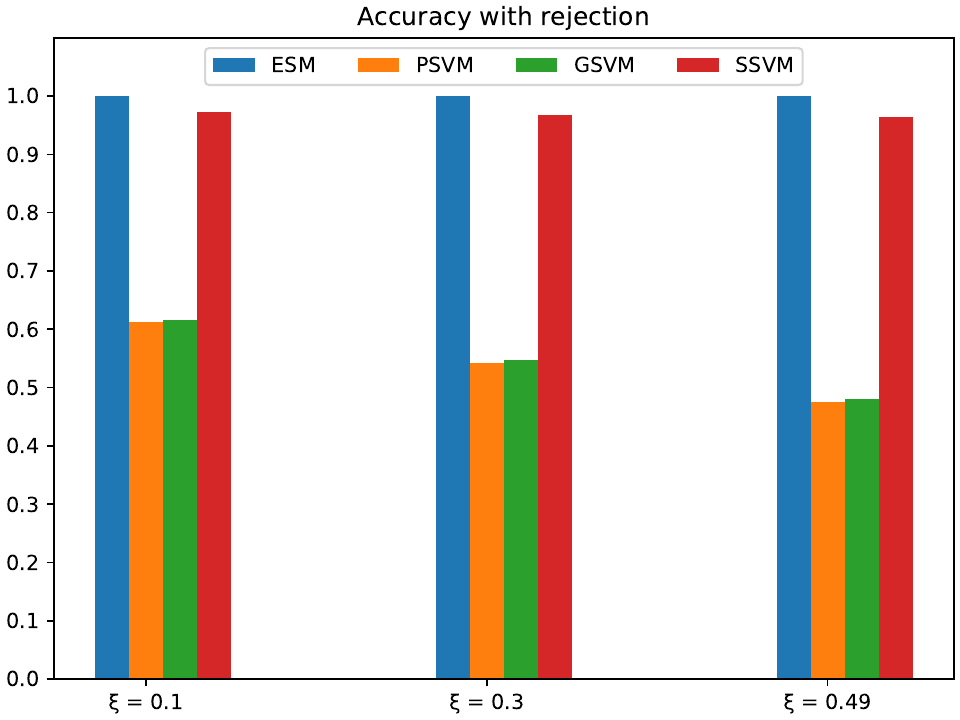}
\end{subfigure}
\hspace{1cm}
\begin{subfigure}{0.3\textwidth}
\includegraphics[scale=0.25]{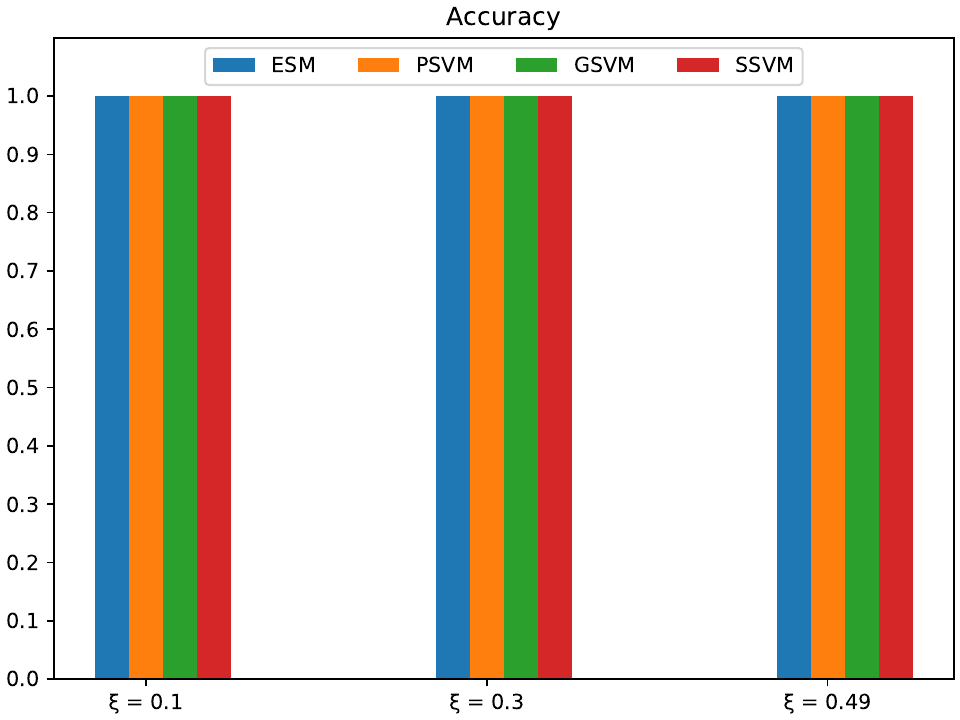}
\end{subfigure}
\vspace{0.3cm}
\begin{subfigure}{0.3\textwidth}
\includegraphics[scale=0.25]{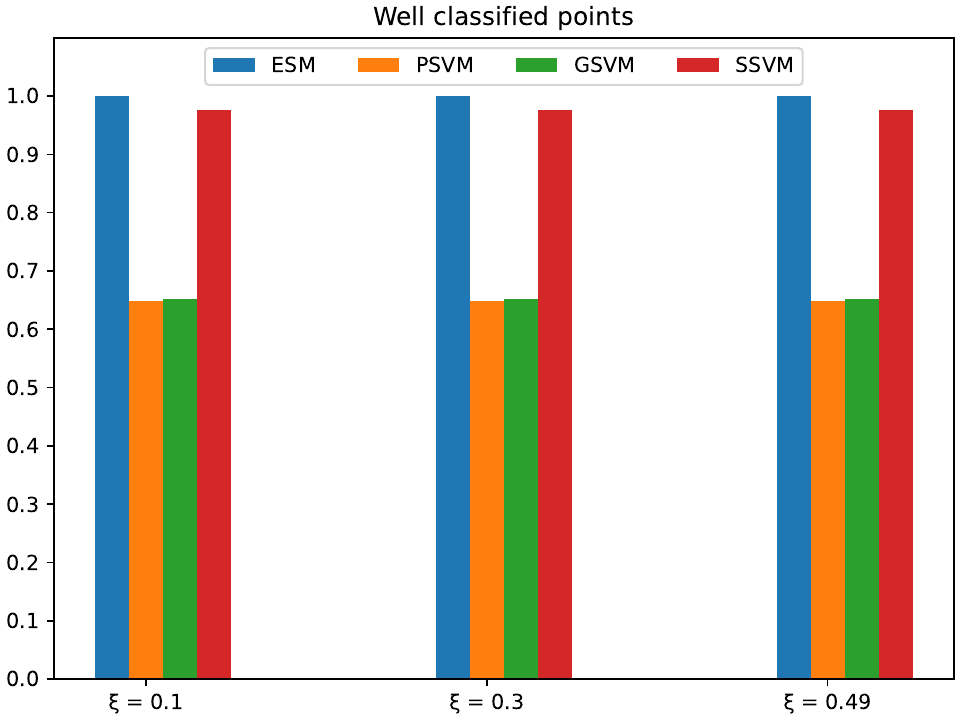}
\end{subfigure}
\hspace{0.1cm}
\begin{subfigure}{0.3\textwidth}
\includegraphics[scale=0.25]{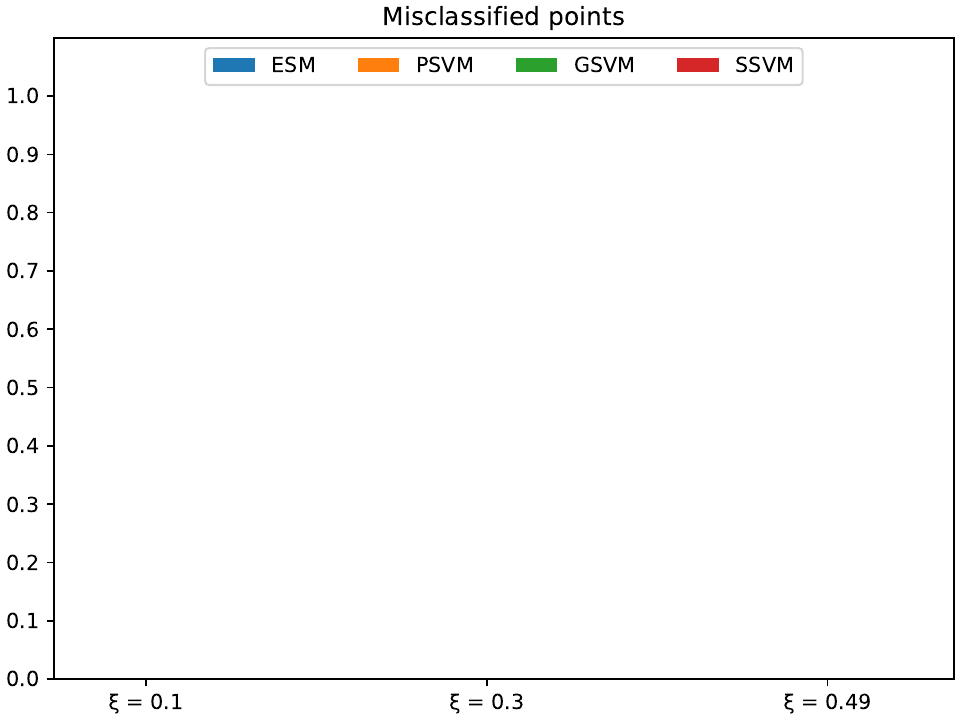}
\end{subfigure}
\hspace{0.1cm}
\begin{subfigure}{0.3\textwidth}
\includegraphics[scale=0.25]{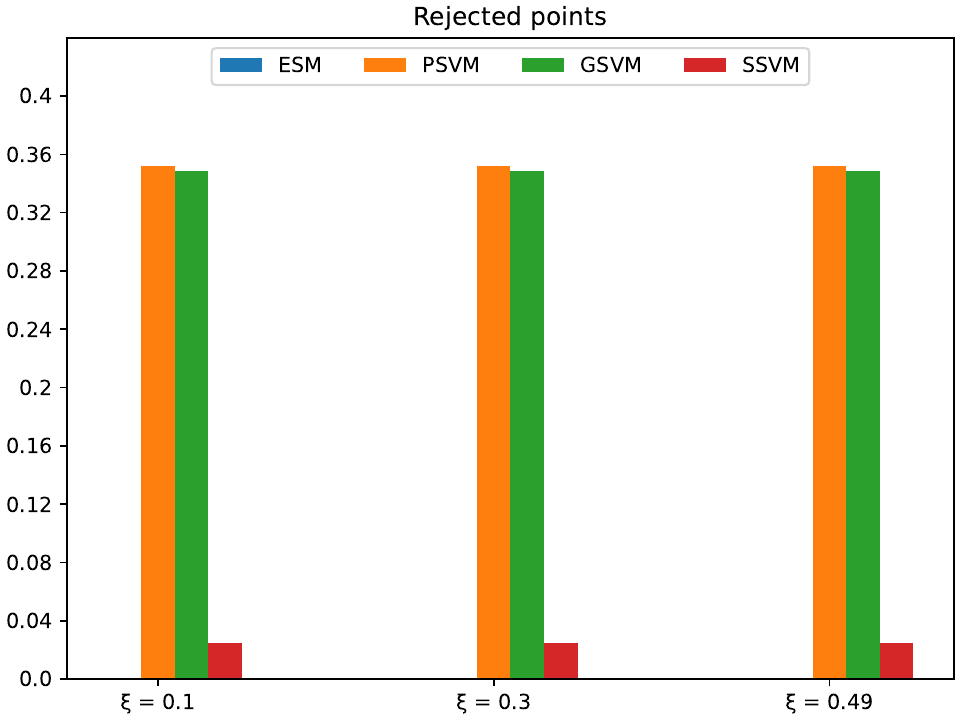}
\end{subfigure}
\end{center}
\caption{Numerical results for the real dataset Mesothelioma}
\label{fig:meso}
\end{figure}

\begin{figure}[h!]
\begin{center}
\begin{subfigure}{0.3\textwidth}
\includegraphics[scale=0.25]{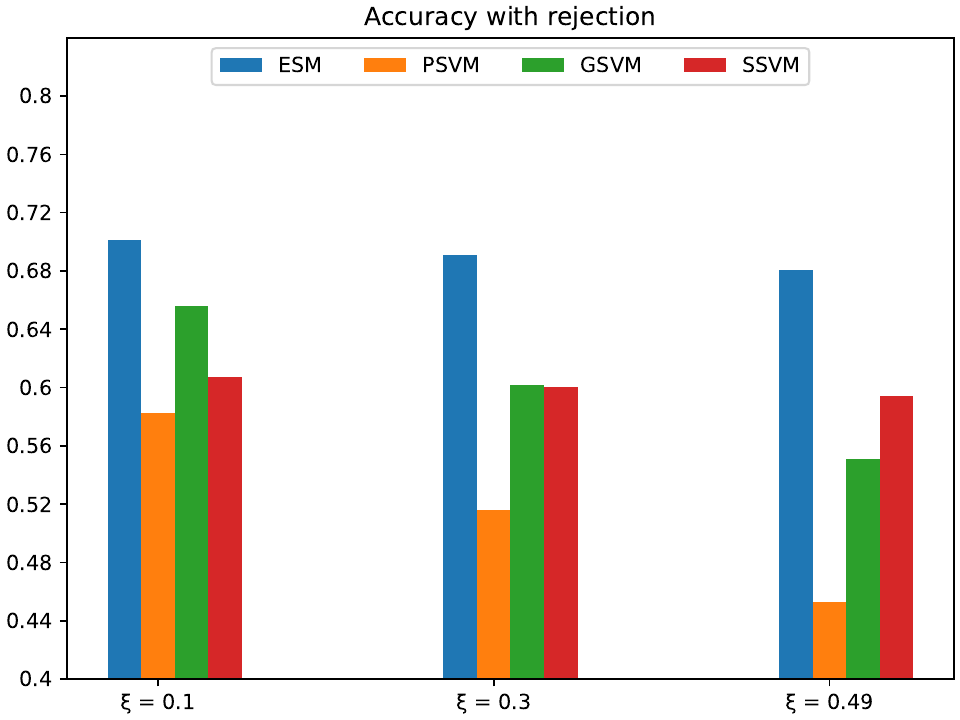}
\end{subfigure}
\hspace{1cm}
\begin{subfigure}{0.3\textwidth}
\includegraphics[scale=0.25]{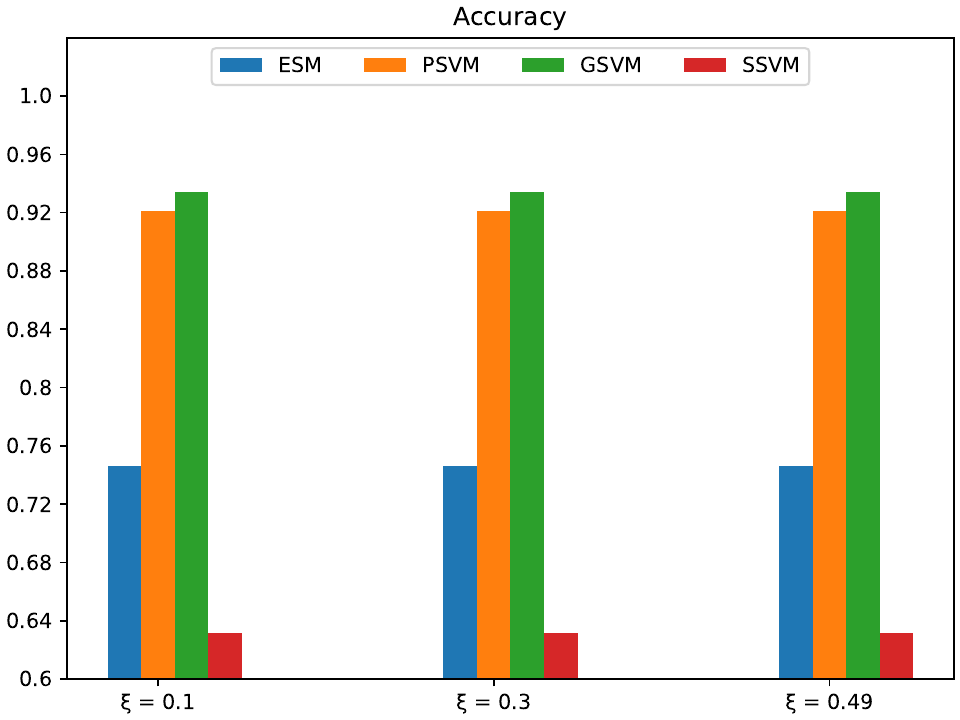}
\end{subfigure}
\vspace{0.3cm}
\begin{subfigure}{0.3\textwidth}
\includegraphics[scale=0.25]{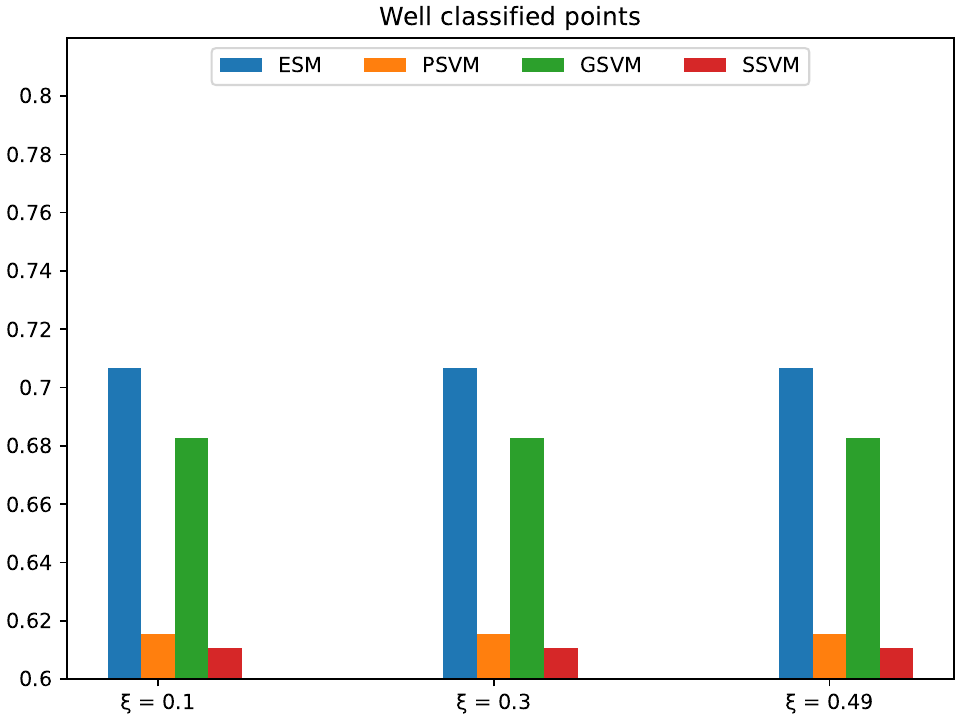}
\end{subfigure}
\hspace{0.1cm}
\begin{subfigure}{0.3\textwidth}
\includegraphics[scale=0.25]{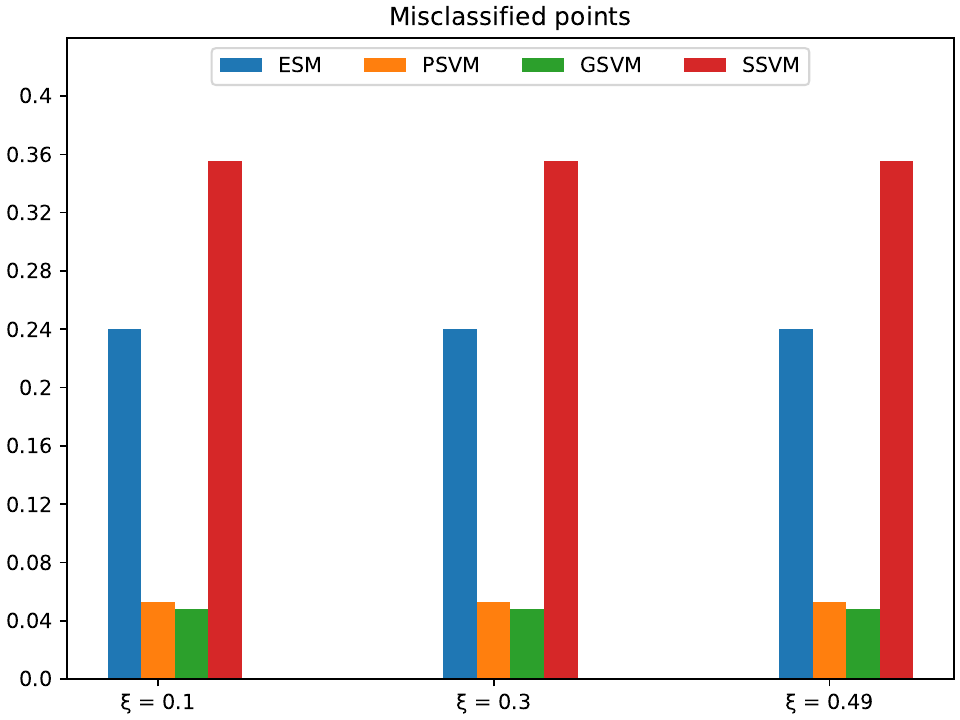}
\end{subfigure}
\hspace{0.1cm}
\begin{subfigure}{0.3\textwidth}
\includegraphics[scale=0.25]{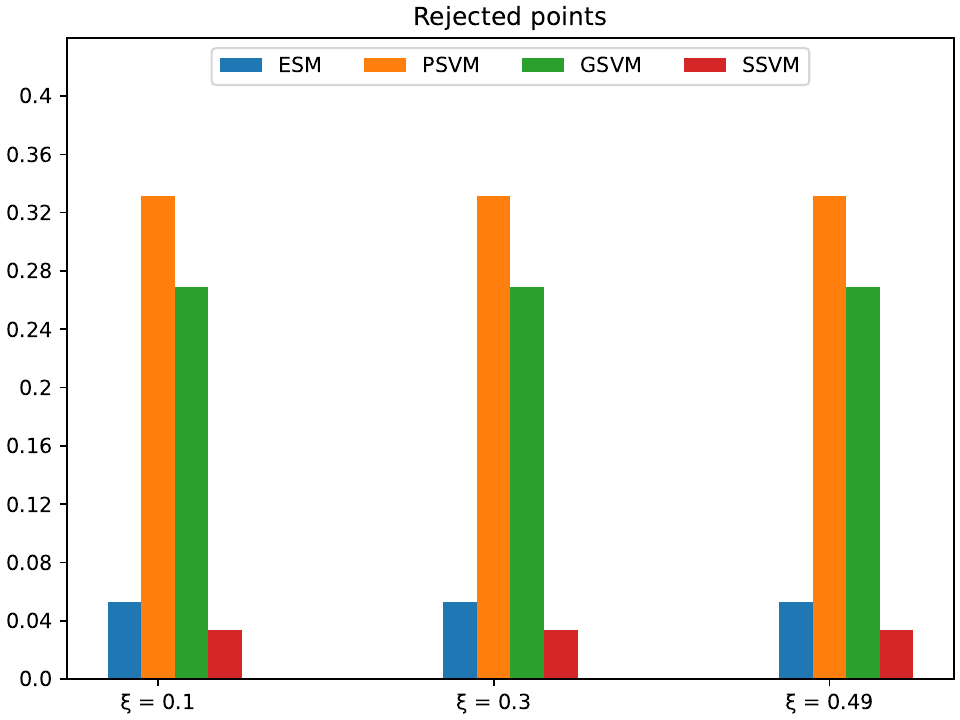}
\end{subfigure}
\end{center}
\caption{Numerical results for the real dataset Sonar}
\label{fig:son}
\end{figure}

\begin{figure}[h!]
\begin{center}
\begin{subfigure}{0.3\textwidth}
\includegraphics[scale=0.25]{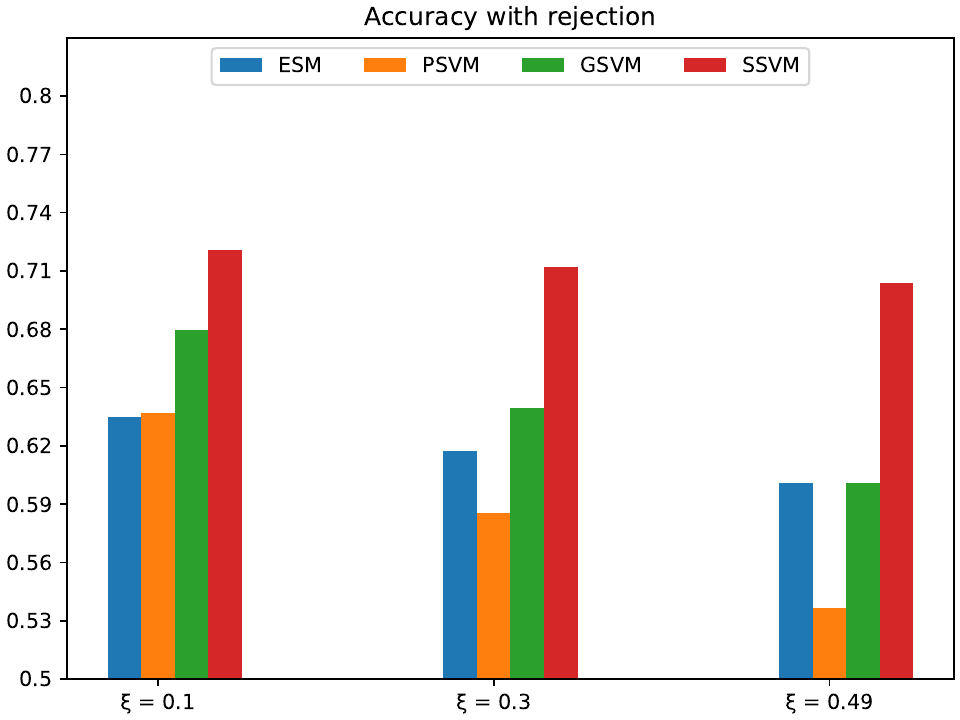}
\end{subfigure}
\hspace{1cm}
\begin{subfigure}{0.3\textwidth}
\includegraphics[scale=0.25]{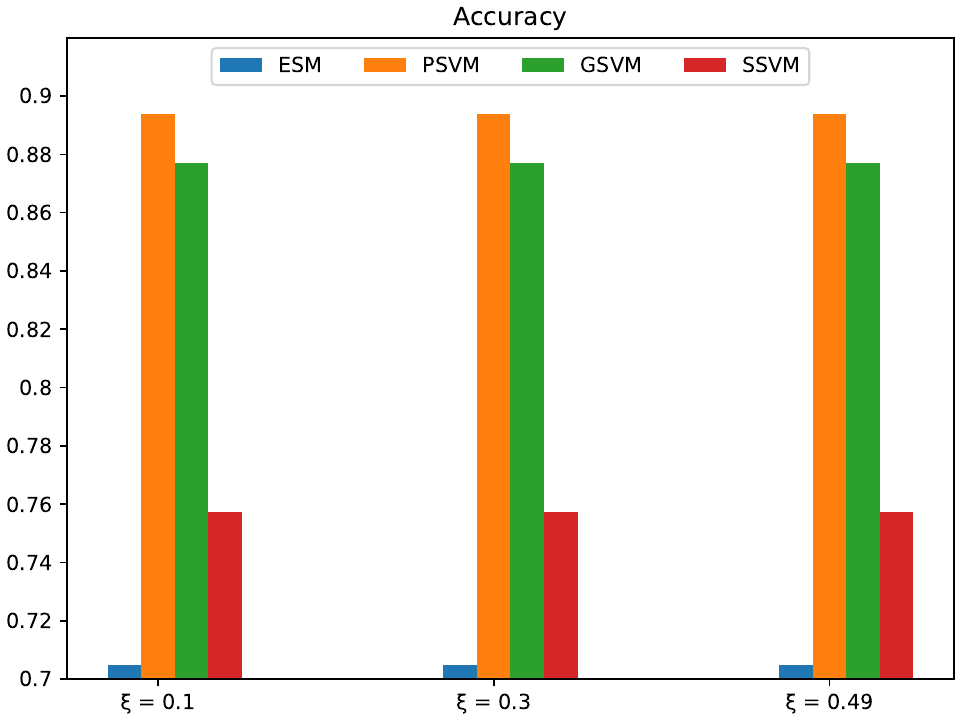}
\end{subfigure}
\vspace{0.3cm}
\begin{subfigure}{0.3\textwidth}
\includegraphics[scale=0.25]{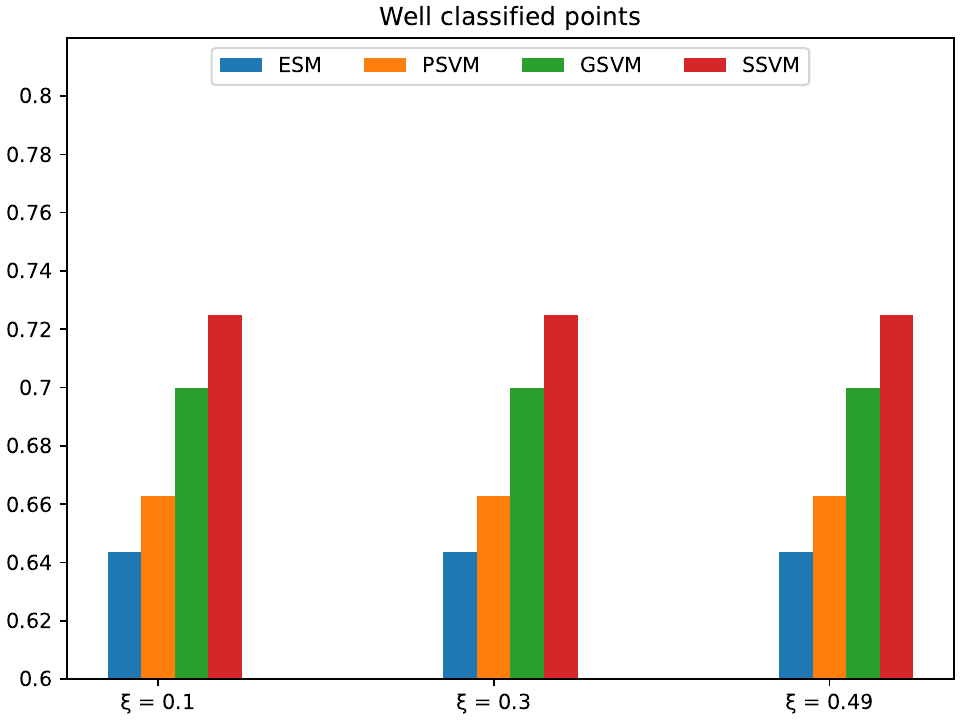}
\end{subfigure}
\hspace{0.1cm}
\begin{subfigure}{0.3\textwidth}
\includegraphics[scale=0.25]{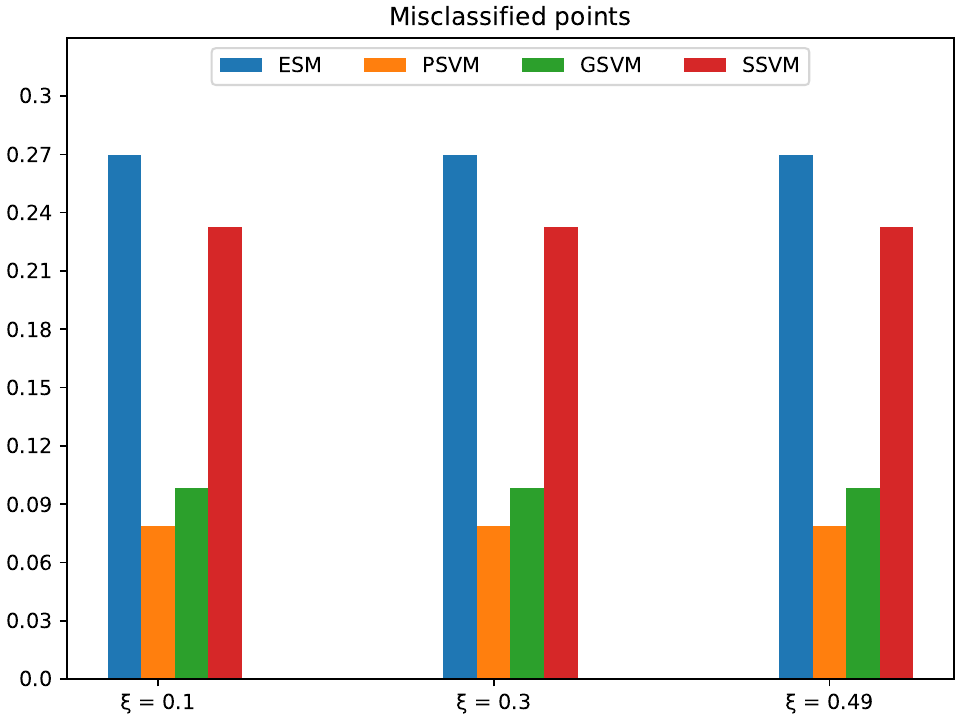}
\end{subfigure}
\hspace{0.1cm}
\begin{subfigure}{0.3\textwidth}
\includegraphics[scale=0.25]{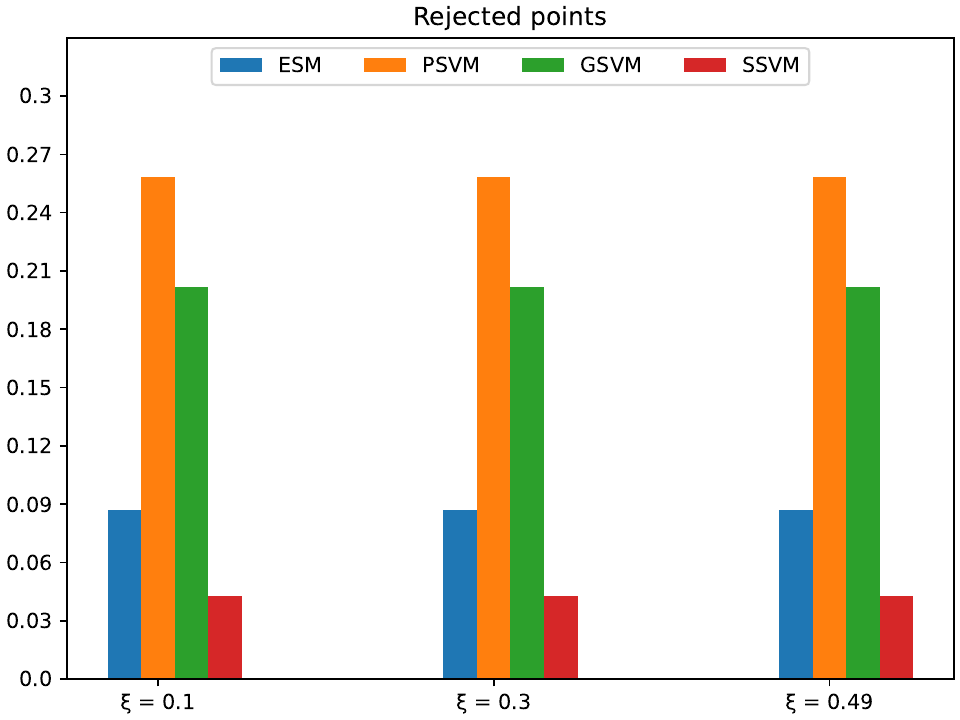}
\end{subfigure}
\end{center}
\caption{Numerical results for the real dataset svmguide3}
\label{fig:svmguide3}
\end{figure}

\begin{figure}[h!]
\begin{center}
\begin{subfigure}{0.3\textwidth}
\includegraphics[scale=0.25]{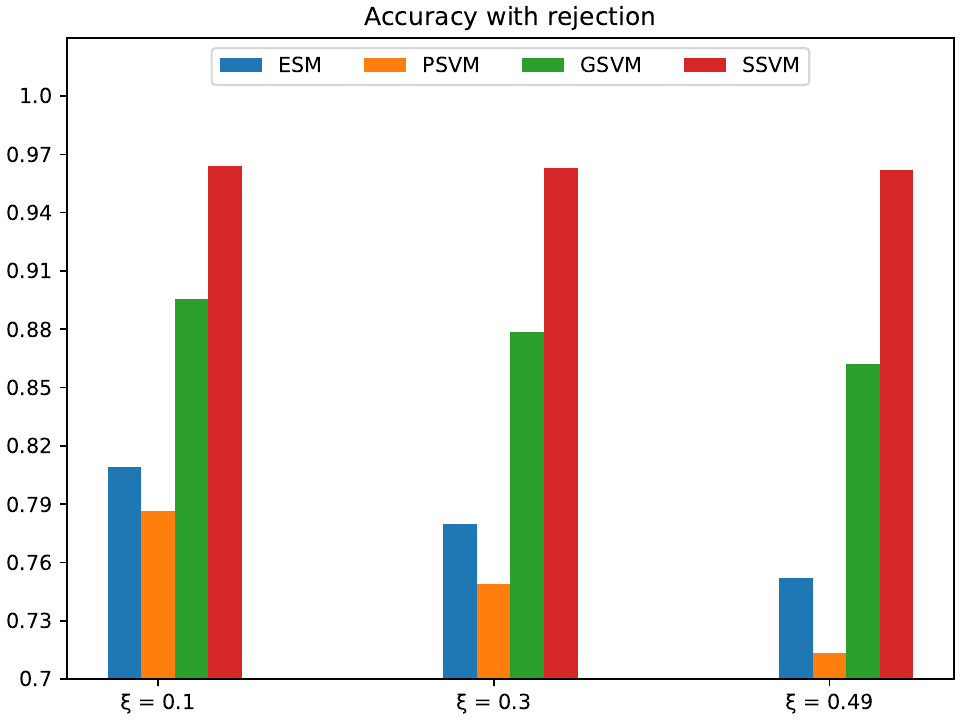}
\end{subfigure}
\hspace{1cm}
\begin{subfigure}{0.3\textwidth}
\includegraphics[scale=0.25]{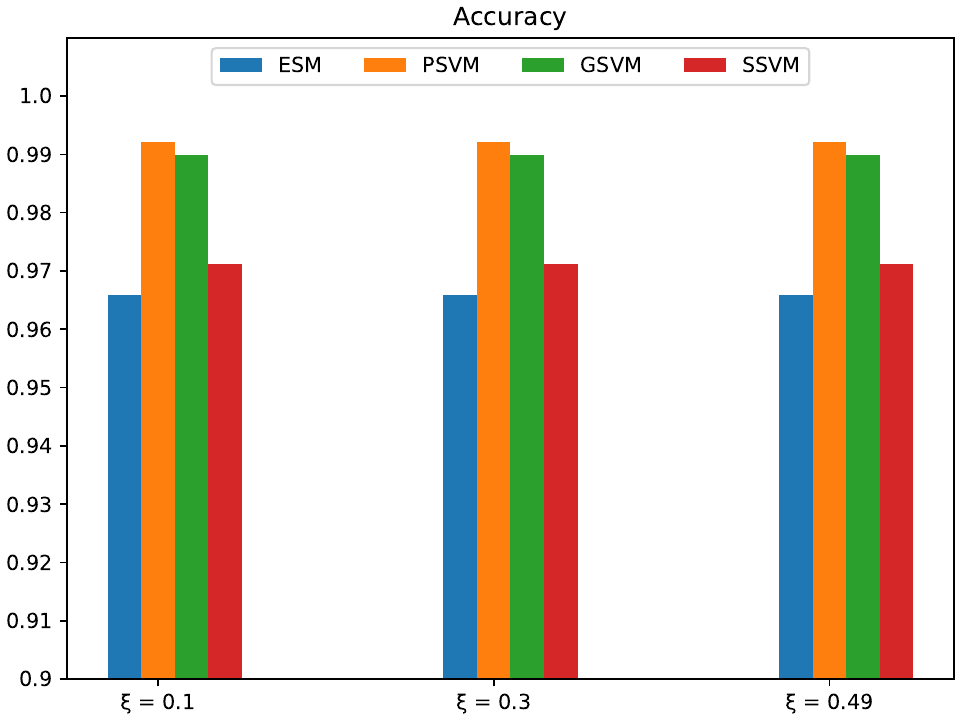}
\end{subfigure}
\vspace{0.3cm}
\begin{subfigure}{0.3\textwidth}
\includegraphics[scale=0.25]{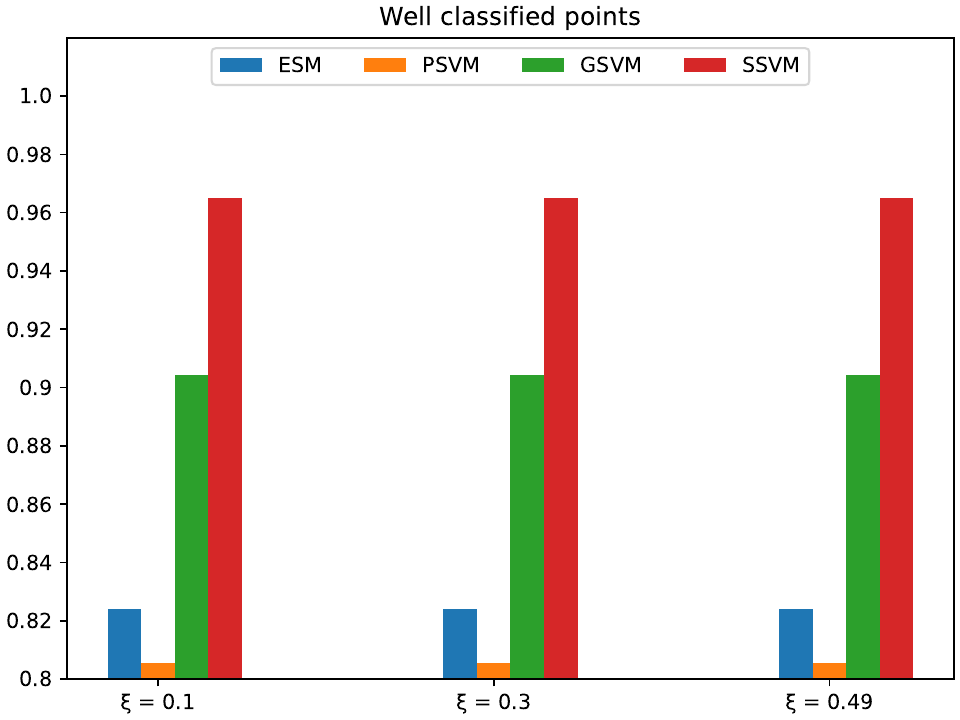}
\end{subfigure}
\hspace{0.1cm}
\begin{subfigure}{0.3\textwidth}
\includegraphics[scale=0.25]{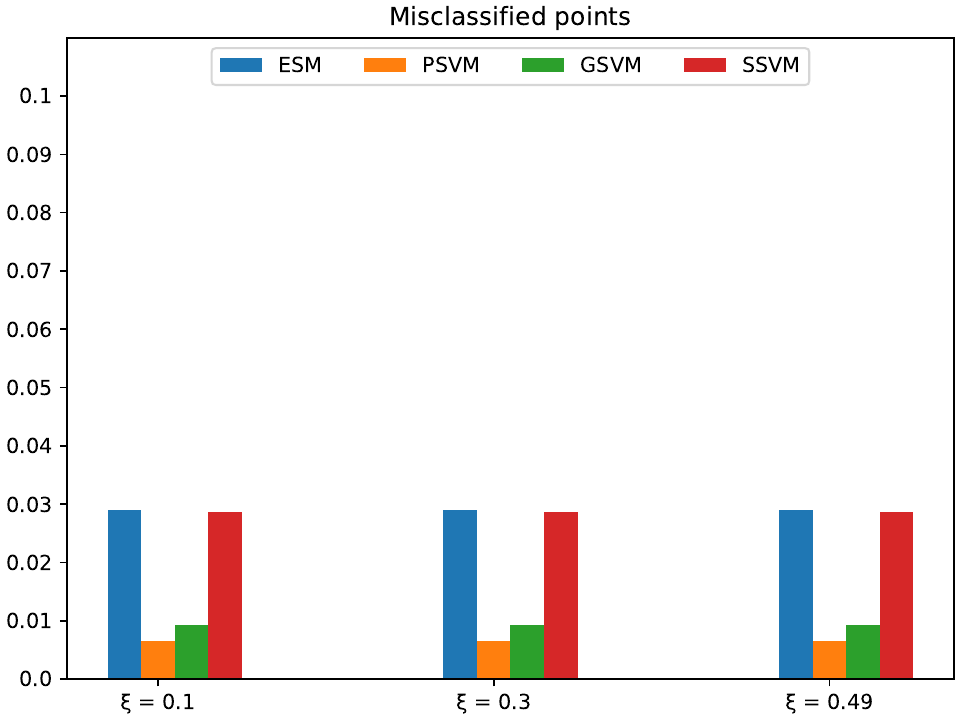}
\end{subfigure}
\hspace{0.1cm}
\begin{subfigure}{0.3\textwidth}
\includegraphics[scale=0.25]{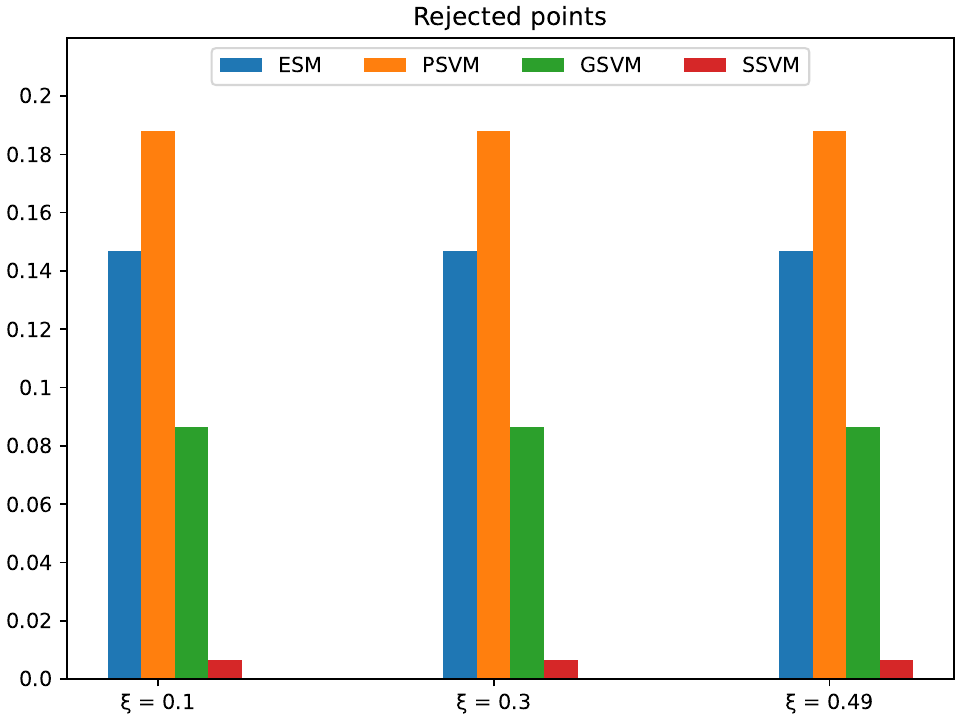}
\end{subfigure}
\end{center}
\caption{Numerical results for the real dataset w1a}
\label{fig:w1a}
\end{figure}

\noindent
Table \ref{tab:bestparamsreal} shows the best hyper-parameters found for each dataset. We observe that our model is robust with respect to the rejection cost $\xi$. The results are reported in Table \ref{tab:realresults} together with their graphic representation in Figures \ref{fig:a1a}--\ref{fig:w1a}.

\begin{table}[h!]
\begin{tabular}{|c|c|c|c|c|c|c|}
\hline
Dataset & $\xi$ & $C_1$ & $C_2$ & $C_3$ & $r$& \textsf{max\_bundle\_it}\\
\hline
\multirow{ 3}{*}{a1a} &0.1 &0.05    & 1e-12  & 5e-9  &8   &40\\
\cmidrule{2-7}
&0.3 &0.05    & 1e-12  & 5e-9  &8   &40\\
\cmidrule{2-7}
&0.49 &0.05    & 1e-12  & 5e-9  &8   &40\\
\hline
\multirow{ 3}{*}{Australian} &0.1 &0.0008&1e-12&10000&8&40\\
\cmidrule{2-7}
&0.3 &0.0008&1e-12&10000&8&40\\
\cmidrule{2-7}
&0.49 &0.0008&1e-12&10000&8&40\\
\hline
\multirow{ 3}{*}{Breast} &0.1 &1000.0&1e-12&8000.0&5&20\\
\cmidrule{2-7}
&0.3 &1000.0&1e-12&8000.0&5&20\\
\cmidrule{2-7}
& 0.49&1000.0&1e-12&8000.0&5&20\\
\hline
\multirow{ 3}{*}{Breast Wisconsin} &0.1 &1000.0&1e-12&8e-09&5&20\\
\cmidrule{2-7}
&0.3 &1000.0&1e-12&8e-09&5&20\\
\cmidrule{2-7}
&0.49 &1000.0&1e-12&8e-09&5&20\\
\hline
\multirow{ 3}{*}{Diabetes} &0.1 &3e-10&1e-12&30000.0&8&40\\
\cmidrule{2-7}
&0.3 &3e-10&1e-12&30000.0&8&40\\
\cmidrule{2-7}
&0.49 &3e-10&1e-12&30000.0&8&40\\
\hline
\multirow{ 3}{*}{Divorce} &0.1 &1e-12&1e-12&5e-12&5&20\\
\cmidrule{2-7}
&0.3 &1e-12&1e-12&5e-12&5&20\\
\cmidrule{2-7}
&0.49 &1e-12&1e-12&5e-12&5&20\\
\hline
\multirow{ 3}{*}{Fertility} &0.1 &0.001&1e-12&0.1&5&20\\
\cmidrule{2-7}
&0.3 &0.001&1e-12&0.1&5&20\\
\cmidrule{2-7}
&0.49 &0.001&1e-12&0.1&5&20\\
\hline
\multirow{ 3}{*}{Flowmeters} &0.1 &5e-9&1e-12&10&8&40\\
\cmidrule{2-7}
&0.3 &5e-9&1e-12&10&8&40\\
\cmidrule{2-7}
&0.49 &5e-9&1e-12&106&8&40\\
\hline
\multirow{ 3}{*}{Gallstone} &0.1 &1000&1e-12&500&5&20\\
\cmidrule{2-7}
&0.3 &1000&1e-12&500&5&20\\
\cmidrule{2-7}
&0.49 &1000&1e-12&500&5&20\\
\hline
\multirow{ 3}{*}{Heart} &0.1 &3e-08&1e-12&3000.0&8&40\\
\cmidrule{2-7}
&0.3 &3e-08&1e-12&3000.0&8&40\\
\cmidrule{2-7}
&0.49 &3e-08&1e-12&3000.0&8&40\\
\hline
\multirow{ 3}{*}{Ionosphere} &0.1 &8e-11&1e-12&1e-12&5&20\\
\cmidrule{2-7}
&0.3 &8e-11&1e-12&1e-12&5&20\\
\cmidrule{2-7}
&0.49 &8e-11&1e-12&1e-12&5&20\\
\hline
\multirow{ 3}{*}{Liver} &0.1 &3e-08&1e-12&3000.0&8&40\\
\cmidrule{2-7}
&0.3 &3e-08&1e-12&3000.0&8&40\\
\cmidrule{2-7}
&0.49 &3e-08&1e-12&3000.0&8&40\\
\hline
\multirow{ 3}{*}{Mesothelioma} &0.1 &100&1e-12&0.10&5&20\\
\cmidrule{2-7}
&0.3 &100&1e-12&0.1&5&20\\
\cmidrule{2-7}
&0.49 &100&1e-12&0.1&5&20\\
\hline
\multirow{ 3}{*}{Sonar} &0.1 &1000.0&1e-12&30.0&5&20\\
\cmidrule{2-7}
&0.3 &1000.0&1e-12&30.0&5&20\\
\cmidrule{2-7}
&0.49 &1000.0&1e-12&30.0&5&20\\
\hline
\multirow{ 3}{*}{svmguide3} &0.1 &3.0    &1e-12   &5000   & 8  &40\\
\cmidrule{2-7}
&0.3 &3.0    &1e-12   &5000   & 8  &40\\
\cmidrule{2-7}
&0.49 &3.0    &1e-12   &5000   & 8  &40\\
\hline
\multirow{ 3}{*}{w1a} &0.1 &3e-8    & 1e-12  & 5e-6  &5   &20\\
\cmidrule{2-7}
&0.3 &3e-8    & 1e-12  & 5e-6  &5   &20\\
\cmidrule{2-7}
&0.49 &3e-8    & 1e-12  & 5e-6  &5   &20\\
\hline
\end{tabular}
\caption{Best ESM hyper-parameters for the datasets}
\label{tab:bestparamsreal}
\end{table}

\begin{table}[h!]
\renewcommand\tabcolsep{1.0pt}
{\footnotesize{
\begin{tabular}{c|c|ccccc|ccccc|ccccc|ccccc}
  & $\xi$  & \multicolumn{5}{c|}{ESM} & \multicolumn{5}{c|}{PSVM}& \multicolumn{5}{c|}{GSVM}& \multicolumn{5}{c}{SSVM}\\
&  & AR & A & WC & MC & R & AR & A & WC & MC & R & AR & A & WC & MC & R &AR & A & WC & MC & R \\
\hline
 \multirow{3}*{\rotatebox[origin=c]{270}{a1a}}&0.1 &0.7&0.74&0.71&0.24&0.05&0.63&0.91&0.66&0.06&0.28&0.06&0.96&0.15&0.01&0.85&\bf0.77&0.78&0.77&0.22&0.01\\
\cmidrule{2-22}
  &0.3 &0.69&0.74&0.71&0.24&0.05&0.58&0.91&0.66&0.06&0.28&-0.1&0.96&0.15&0.01&0.85&\bf0.77&0.78&0.77&0.22&0.01\\
\cmidrule{2-22}
  &0.49 &0.68&0.74&0.71&0.24&0.05&0.53&0.91&0.66&0.06&0.28&-0.27&0.96&0.15&0.01&0.85&\bf0.77&0.78&0.77&0.22&0.01\\
 \hline
 \multirow{3}*{\rotatebox[origin=c]{270}{Aust.}}&0.1 & 0.53&0.67&0.55&0.27&0.18&\bf{0.57}&0.61&0.58&0.37&0.05&0.54&0.9&0.58&0.07&0.35&0.48&0.48&0.48&0.52&0.0\\
\cmidrule{2-22}
  &0.3 & 0.5&0.67&0.55&0.27&0.18&\bf{0.56}&0.61&0.58&0.37&0.05&0.47&0.9&0.58&0.07&0.35&0.48&0.48&0.48&0.52&0.0\\
\cmidrule{2-22}
  &0.49 & 0.46&0.67&0.55&0.27&0.18&\bf{0.55}&0.61&0.58&0.37&0.05&0.41&0.9&0.58&0.07&0.35&0.48&0.48&0.48&0.52&0.0\\
 \hline
 \multirow{3}*{\rotatebox[origin=c]{270}{Breast}}&0.1 & \bf{0.95}&0.95&0.95&0.05&0.0&0.88&0.99&0.89&0.01&0.1&0.9&0.99&0.9&0.01&0.09&0.94&0.95&0.94&0.05&0.02\\
\cmidrule{2-22}
  &0.3 &\bf{0.95}&0.95&f0.95&0.05&0.0&0.86&0.99&0.89&0.01&0.1&0.88&0.99&0.9&0.01&0.09&0.93&0.95&0.94&0.05&0.02\\
\cmidrule{2-22}
  &0.49 & \bf{0.95}&0.95&0.95&0.05&0.0&0.84&0.99&0.89&0.01&0.1&0.86&0.99&0.9&0.01&0.09&0.93&0.95&0.94&0.05&0.02\\
 \hline
 \multirow{3}*{\rotatebox[origin=c]{270}{Bre.W.}}&0.1 & \bf0.88&0.93&0.88&0.07&0.05&0.8&0.99&0.82&0.01&0.17&0.75&0.99&0.77&0.01&0.22&0.63&0.63&0.63&0.37&0.0\\
\cmidrule{2-22}
  &0.3 & \bf0.86&0.93&0.88&0.07&0.05&0.77&0.99&0.82&0.01&0.17&0.71&0.99&0.77&0.01&0.22&0.63&0.63&0.63&0.37&0.0\\
\cmidrule{2-22}
  &0.49 & \bf0.85&0.93&0.88&0.07&0.05&0.73&0.99&0.82&0.01&0.17&0.66&0.99&0.77&0.01&0.22&0.63&0.63&0.63&0.37&0.0\\
 \hline
 \multirow{3}*{\rotatebox[origin=c]{270}{Diab.}}&0.1 & 0.5&0.73&0.53&0.2&0.27&0.65&0.84&0.67&0.12&0.21&0.6&0.87&0.63&0.09&0.28&\bf0.66&0.66&0.66&0.34&0.01\\
\cmidrule{2-22}
  &0.3 & 0.45&0.73&0.53&0.2&0.27&0.6&0.84&0.67&0.12&0.21&0.55&0.87&0.63&0.09&0.28&\bf0.65&0.66&0.66&0.34&0.01\\
\cmidrule{2-22}
  &0.49 & 0.4&0.73&0.53&0.2&0.27&0.56&0.84&0.67&0.12&0.21&0.49&0.87&0.63&0.09&0.28&\bf0.65&0.66&0.66&0.34&0.01\\
 \hline
 \multirow{3}*{\rotatebox[origin=c]{270}{Divor.}}&0.1 & \bf0.97&0.97&0.97&0.03&0.0&0.93&1.0&0.94&0.0&0.06&0.9&1.0&0.91&0.0&0.09&0.9&1.0&0.91&0.0&0.09\\
\cmidrule{2-22}
  &0.3 & \bf0.97&0.97&0.97&0.03&0.0&0.92&1.0&0.94&0.0&0.06&0.88&1.0&0.91&0.0&0.09&0.88&1.0&0.91&0.0&0.09\\
\cmidrule{2-22}
  &0.49 & \bf0.97&0.97&0.97&0.03&0.0&0.9&1.0&0.94&0.0&0.06&0.86&1.0&0.91&0.0&0.09&0.86&1.0&0.91&0.0&0.09\\
 \hline
 \multirow{3}*{\rotatebox[origin=c]{270}{Fert.}}&0.1 & 0.45&0.89&0.49&0.06&0.45&0.6&0.86&0.63&0.1&0.27&0.61&0.94&0.64&0.04&0.32&\bf0.76&0.83&0.77&0.16&0.07\\
\cmidrule{2-22}
  &0.3 & 0.35&0.89&0.49&0.06&0.45&0.55&0.86&0.63&0.1&0.27&0.54&0.94&0.64&0.04&0.32&\bf0.75&0.83&0.77&0.16&0.07\\
\cmidrule{2-22}
  &0.49 & 0.27&0.89&0.49&0.06&0.45&0.5&0.86&0.63&0.1&0.27&0.48&0.94&0.64&0.04&0.32&\bf0.74&0.83&0.77&0.16&0.07\\
 \hline
 \multirow{3}*{\rotatebox[origin=c]{270}{Flowm.}}&0.1 & 0.38&0.93&0.44&0.03&0.53&0.6&0.9&0.63&0.07&0.3&\bf0.71&0.97&0.74&0.02&0.24&0.63&0.63&0.63&0.37&0.0\\
\cmidrule{2-22}
  &0.3 & 0.28&0.93&0.44&0.03&0.53&0.54&0.9&0.63&0.07&0.3&\bf0.66&0.97&0.74&0.02&0.24&0.63&0.63&0.63&0.37&0.0\\
\cmidrule{2-22}
  &0.49 & 0.18&0.93&0.44&0.03&0.53&0.49&0.9&0.63&0.07&0.3&0.62&0.97&0.74&0.02&0.24&\bf0.63&0.63&0.63&0.37&0.0\\
 \hline
 \multirow{3}*{\rotatebox[origin=c]{270}{Gallst.}}&0.1 & 0.33&0.54&0.36&0.31&0.33&\bf0.55&0.87&0.58&0.09&0.33&0.29&0.89&0.35&0.04&0.61&0.5&0.5&0.5&0.5&0.0\\
\cmidrule{2-22}
  &0.3 & 0.27&0.54&0.36&0.31&0.33&0.48&0.87&0.58&0.09&0.33&0.17&0.89&0.35&0.04&0.61&\bf0.5&0.5&0.5&0.5&0.0\\
\cmidrule{2-22}
  &0.49 & 0.2&0.54&0.36&0.31&0.33&0.42&0.87&0.58&0.09&0.33&0.05&0.89&0.35&0.04&0.61&\bf0.5&0.5&0.5&0.5&0.0\\
 \hline
 \multirow{3}*{\rotatebox[origin=c]{270}{Heart}}&0.1 & 0.76&0.81&0.77&0.17&0.06&0.66&0.91&0.68&0.07&0.25&0.6&0.9&0.63&0.07&0.3&\bf0.79&0.8&0.79&0.2&0.01\\
\cmidrule{2-22}
  &0.3 & 0.75&0.81&0.77&0.17&0.06&0.61&0.91&0.68&0.07&0.25&0.54&0.9&0.63&0.07&0.3&\bf0.79&0.8&0.79&0.2&0.01\\
\cmidrule{2-22}
  &0.49 & 0.74&0.81&0.77&0.17&0.06&0.56&0.91&0.68&0.07&0.25&0.48&0.9&0.63&0.07&0.3&\bf0.78&0.8&0.79&0.2&0.01\\
 \hline
 \multirow{3}*{\rotatebox[origin=c]{270}{Iono.}}&0.1 & 0.61&0.81&0.63&0.15&0.22&0.79&0.96&0.81&0.04&0.15&\bf0.81&0.93&0.82&0.06&0.11&0.78&0.79&0.78&0.21&0.01\\
\cmidrule{2-22}
  &0.3 & 0.57&0.81&0.63&0.15&0.22&0.76&0.96&0.81&0.04&0.15&\bf0.79&0.93&0.82&0.06&0.11&0.78&0.79&0.78&0.21&0.01\\
\cmidrule{2-22}
  &0.49 & 0.52&0.81&0.63&0.15&0.22&0.73&0.96&0.81&0.04&0.15&0.77&0.93&0.82&0.06&0.11&\bf0.78&0.79&0.78&0.21&0.01\\
 \hline
 \multirow{3}*{\rotatebox[origin=c]{270}{Liver}}&0.1 & 0.52&0.52&0.52&0.48&0.0&\bf0.71&0.79&0.72&0.19&0.1&0.62&0.79&0.64&0.17&0.19&0.64&0.65&0.64&0.35&0.01\\
\cmidrule{2-22}
  &0.3 & 0.52&0.52&0.52&0.48&0.0&\bf0.69&0.79&0.72&0.19&0.1&0.58&0.79&0.64&0.17&0.19&0.64&0.65&0.64&0.35&0.01\\
\cmidrule{2-22}
  &0.49 & 0.52&0.52&0.52&0.48&0.0&\bf0.67&0.79&0.72&0.19&0.1&0.55&0.79&0.64&0.17&0.19&0.64&0.65&0.64&0.35&0.01\\
 \hline
 \multirow{3}*{\rotatebox[origin=c]{270}{Meso.}}&0.1 & \bf1.0&1.0&1.0&0.0&0.0&0.61&1.0&0.65&0.0&0.35&0.62&1.0&0.65&0.0&0.35&0.97&1.0&0.98&0.0&0.02\\
\cmidrule{2-22}
  &0.3 & \bf1.0&1.0&1.0&0.0&0.0&0.54&1.0&0.65&0.0&0.35&0.55&1.0&0.65&0.0&0.35&0.97&1.0&0.98&0.0&0.02\\
\cmidrule{2-22}
  &0.49 & \bf1.0&1.0&1.0&0.0&0.0&0.48&1.0&0.65&0.0&0.35&0.48&1.0&0.65&0.0&0.35&0.96&1.0&0.98&0.0&0.02\\
 \hline
 \multirow{3}*{\rotatebox[origin=c]{270}{Son.}}&0.1 & \bf0.7&0.75&0.71&0.24&0.05&0.58&0.92&0.62&0.05&0.33&0.66&0.93&0.68&0.05&0.27&0.61&0.63&0.61&0.36&0.03\\
\cmidrule{2-22}
  &0.3 & \bf0.69&0.75&0.71&0.24&0.05&0.52&0.92&0.62&0.05&0.33&0.6&0.93&0.68&0.05&0.27&0.6&0.63&0.61&0.36&0.03\\
\cmidrule{2-22}
  &0.49 & \bf0.68&0.75&0.71&0.24&0.05&0.45&0.92&0.62&0.05&0.33&0.55&0.93&0.68&0.05&0.27&0.59&0.63&0.61&0.36&0.03\\
 \hline
  \multirow{3}*{\rotatebox[origin=c]{270}{svmg.3}}&0.1 &0.63&0.7&0.64&0.27&0.09&0.64&0.89&0.66&0.08&0.26&0.68&0.88&0.7&0.1&0.2&\bf0.72&0.76&0.72&0.23&0.04\\
\cmidrule{2-22}
  &0.3 &0.62&0.7&0.64&0.27&0.09&0.59&0.89&0.66&0.08&0.26&0.64&0.88&0.7&0.1&0.2&\bf0.71&0.76&0.72&0.23&0.04\\
\cmidrule{2-22}
  &0.49 &0.6&0.7&0.64&0.27&0.09&0.54&0.89&0.66&0.08&0.26&0.6&0.88&0.7&0.1&0.2&\bf0.7&0.76&0.72&0.23&0.04\\
 \hline
  \multirow{3}*{\rotatebox[origin=c]{270}{w1a}}&0.1 &0.81&0.97&0.82&0.03&0.15&0.79&0.99&0.81&0.01&0.19&0.9&0.99&0.9&0.01&0.09&\bf0.96&0.97&0.96&0.03&0.01\\
\cmidrule{2-22}
  &0.3 &0.78&0.97&0.82&0.03&0.15&0.75&0.99&0.81&0.01&0.19&0.88&0.99&0.9&0.01&0.09&\bf0.96&0.97&0.96&0.03&0.01\\
\cmidrule{2-22}
  &0.49 &0.75&0.97&0.82&0.03&0.15&0.71&0.99&0.81&0.01&0.19&0.86&0.99&0.9&0.01&0.09&\bf0.96&0.97&0.96&0.03&0.01\\
 \hline
  
\end{tabular}
}}
\caption{Results for the datasets, in bold the best AR score}
\label{tab:realresults}
\end{table}

\smallskip
\noindent
%
The figures illustrate how the different kernels have quite different behavior in terms of the different metrics, showing that---as expected---the choice of the kernel significantly impacts the classification performances of the methods. Remarkably, ESM is still quite different from all the other classifiers; while it may behave similarly to some of the kernels in some cases it is typically not to the same one, and it can be significantly different from them all. In some cases ESM is comparable with the best (but, often, not to all) SVM approaches, in a few cases it actually yields significantly worse results, but yet in other cases (and a significant number of them) it clearly outperform all the standard SVM approaches, despite not being endowed with the flexibility that the kernels provide. Thus, ESM provides a novel choice that can be more appropriate to some datasets. 

\section{Conclusions}

We have proposed a fully functional implementation of the ``Separation by Convex Bodies'' approach that moves away from the traditional hyperplane separator characterising the SVM approach to the ``next simpler'' convex set, i.e., the ellipsoid. This leads to some nontrivial issues, having primarily to do that the two areas of the space separated by an hyperplane are fully symmetric, but the same does not hold for an ellipsoid, where one of them is convex and the other reverse convex. Yet, appropriate reformulation tricks allow to write the training problem as a convex program. The natural formulation is a rather large-dimensional SDP, thus solving it in reasonable time for the size required by actual classification tasks, considering the need of performing grid search to the hyperparameters, is not trivial. Fortunately, it is well-known that highly accurate solutions are not needed in ML approaches, which allows us to propose heuristic approaches that produce working classifiers at reasonable computational cost. All in all, ESM provides yet another arrow in the quiver when constructing classification approaches which may be particularly interesting when rejecting points that should not really be classified is a primary concern. Future research may entail developing kernel tricks that enhance the flexibility of the approach by embedding the separating ellipsoid in a more appropriate feature space.

\bibliography{Biblio_2}
\end{document}